\documentclass[11pt,reqno]{amsart}

\usepackage{amsfonts}
\usepackage{eurosym}
\usepackage{amssymb}
\usepackage{amsthm}
\usepackage{amsmath}
\usepackage{bm}
\usepackage{cite}
\usepackage{mathrsfs}
\usepackage{xcolor}
\usepackage[OT1]{fontenc}
\usepackage[left=2cm, right=2cm, top=3cm]{geometry}
\usepackage{hyperref}
\hypersetup{colorlinks=true, linkcolor=blue, citecolor=red}

\usepackage{times}

\usepackage[shortlabels]{enumitem}
\usepackage{enumitem}

\numberwithin{equation}{section}
\usepackage[square,numbers,sectionbib]{natbib}
\usepackage{varioref}

\usepackage{listings}
\usepackage{color}
\definecolor{dkgreen}{rgb}{0,0.6,0}
\definecolor{gray}{rgb}{0.5,0.5,0.5}
\definecolor{mauve}{rgb}{0.58,0,0.82}

\allowdisplaybreaks
\newtheorem{theorem}{Theorem}[section]

\newtheorem{lemma}[theorem]{Lemma}

\theoremstyle{remark}
\newtheorem{remark}[theorem]{Remark}

\def\d{{\rm d}}

\def\f{\textbf{\textit{f}}}
\def\uu{\mathbf{u}}

\def\ww{\ww}

\def\ddt{\frac{\d}{\d t}}

\def\pp{\mathcal{P}_\sigma}
\newcommand{\al}{\boldsymbol{\alpha}}

\newcommand{\numberset}{\mathbb}
\newcommand{\N}{\numberset{N}}
\newcommand{\R}{\numberset{R}}

\everymath{\displaystyle}

\def \no#1#2#3 {{\bf #1} (#3), #2.}
\def \eds#1#2#3 {#1, #2, #3.}
\def\an #1{{#1}}

\def\pphi{\pmb\varphi}

\def\d{{\rm d}}

\def\f{\textbf{\textit{f}}}

\def\ddt{\frac{\d}{\d t}}

\def \vphi{\varphi}
\def\dx{\ dx}

\def\multibold #1{\def\arg{#1}%
	\ifx\arg\pto \let\next\relax
	\else
	\def\next{\expandafter
		\def\csname #1#1#1\endcsname{{\boldsymbol #1}}%
		\multibold}%
	\fi \next}

\def\pto{.}

\def\multical #1{\def\arg{#1}%
	\ifx\arg\pto \let\next\relax
	\else
	\def\next{\expandafter
		\def\csname #1#1\endcsname{{\cal #1}}%
		\multical}%
	\fi \next}

\def\multimathop #1 {\def\arg{#1}%
	\ifx\arg\pto \let\next\relax
	\else
	\def\next{\expandafter
		\def\csname #1\endcsname{\mathop{\rm #1}\nolimits}%
		\multimathop}%
	\fi \next}

\multibold
qwertyuiopasdfghjklzxcvbnmQWERTYUIOPASDFGHJKLZXCVBNM.

\multical
QWERTYUIOPASDFGHJKLZXCVBNM.

\multimathop
diag dist div dom mean meas sign supp .

\def\numberset{\mathbb}
\def\N{\numberset{N}}
\def\R{\numberset{R}}
\def\pp{\mathcal{P}_\sigma}
\def\ww{\www}

\everymath{\displaystyle}

\def \no#1#2#3 {{\bf #1} (#3), #2.}
\def \eds#1#2#3 {#1, #2, #3.}

\renewcommand{\d}{\mathrm d}



\begin{document}

\title[MULTI-COMPONENT CONSERVED ALLEN-CAHN EQUATIONS]
{MULTI-COMPONENT CONSERVED ALLEN-CAHN EQUATIONS}
\author[ M. Grasselli \& A. Poiatti]{Maurizio Grasselli \& Andrea Poiatti}

\address{Politecnico di Milano\\ Dipartimento di Matematica\\ Via E. Bonardi, 9 - Milano 20133 - Italy}
\email{maurizio.grasselli@polimi.it, andrea.poiatti@polimi.it}


\begin{abstract}
We consider a multi-component version of the conserved Allen-Cahn equation proposed by J. Rubinstein and P. Sternberg in 1992 as an alternative
model for phase separation. In our case, the free energy is characterized by a mixing entropy density which belongs to a large class of physically relevant
entropies like, e.g., the Boltzmann-Gibbs entropy. We establish the well-posedness
of the Cauchy-Neumann problem with respect to a natural notion of (finite) energy solution which is more regular under appropriate assumptions and
is strictly separated from pure phases if the initial datum \an{is}. We then prove that the energy solution becomes more regular and strictly separated instantaneously.
Also, we show that any finite energy solution converges to a unique equilibrium.
The validity of a dissipative inequality (identity for strong solutions) allows us to analyze the problem within the theory of infinite-dimensional
dissipative dynamical systems. On account of the obtained results, we can associate to our problem a dissipative dynamical system and we can
prove that it has a global attractor as well as an exponential attractor.
\end{abstract}



\maketitle

\section{Introduction}
Phase separation, namely the creation of two (or more) distinct phases from a single homogeneous mixture, is an important phenomenon which characterizes many important processes. In particular, it
has recently become a paradigm in Cell Biology (see, for instance, \cite{D1,D2} and references therein). A well-known mathematical model of phase separation for binary alloys was proposed by J.W. Cahn and J.E. Hilliard \cite{CH1,CH2}.
This model leads to the so-called Cahn-Hilliard equation (see, for instance, \cite{M} and references therein). More precisely, indicating by $\varphi$
the concentration of one species, phase separation can be modeled as a competition between the Boltzmann-Gibbs mixing entropy
$$
S(\varphi) = -\varphi\ln \varphi - (1-\varphi)\ln \varphi
$$
and the demixing effects due to the reciprocal attraction of the molecules of the same species
which can be described, for instance, as follows
$$
D(\varphi) = -\varphi (1-\varphi).
$$
Thus the free energy density is given by the so-called Flory-Huggins potential (see, for instance, \cite{BTP2015} and references therein)
\begin{equation}
\label{dwpotential}
W(\varphi) = -\Theta S(\varphi) + \Theta_c D(\varphi)
\end{equation}
where $\Theta>0$ is the absolute temperature of the mixture and $\Theta_0>0$ is its critical temperature (other constants are set equal to the unity). If $\Theta <\Theta_c$
then $W$ has a double well shape and phase separation takes place. Assuming that the mixture occupies a bounded domain $\Omega \subset \R^d$, $d=2,3$, the previous considerations
lead to the following free energy functional
$$
E(\vphi) :=\int_\Omega W(\varphi)dx+\frac \gamma 2 \int_\Omega \vert\nabla\varphi \vert^2dx,
$$
where the penalization term allows the creation of diffuse interfaces between the two species and also allows a convenient mathematical treatment of the phenomenon (see \cite{E}). Here $\gamma>0$
is related to the thickness of the diffuse interface. The Cahn-Hilliard equation can be introduced as a conserved gradient flow generated by the gradient of the chemical potential
$\mu$ defined by
$$
\mu = \frac{\delta E}{\delta\varphi}   = -\gamma\Delta\varphi + W^\prime(\varphi),
$$
namely, taking constant mobility equal to a constant $m>0$,
$$
\partial_t \varphi = m\Delta \mu.
$$
This equation subject to no-flux (or periodic) boundary conditions entails the conservation of the total mass $\int_\Omega \varphi(t) dx$. An alternative model has been proposed
by J. Rubinstein and P. Sternberg \cite{RS92} by modifying another well-known equation proposed by S.M. Allen and J.W. Cahn \cite{AC79} in order to ensure the mass conservation. The
equation has the form
\begin{equation}
\label{CAC}
\partial_t \varphi = \alpha(\overline{\mu} - \mu),
\end{equation}
where $\alpha>0$ and $\overline{\f}$ is defined by
\begin{equation*}
\overline{f}:=|\Omega |_d^{-1}\int_{\Omega }f(x)dx,
\end{equation*}%
for any integrable $f$. Here $|\Omega|_d$ stands for the $d$-dimensional Lebesgue measure of $\Omega $.
Equation \eqref{CAC} equipped with homogeneous Neumann boundary condition, preserves the total mass. In \cite{RS92} a (formal) asymptotic analysis was performed
with respect to a specific scaling in order to understand the motion of the separating interfaces (see also \cite{MHVB2018} for an important application).
More rigorous results can be found in \cite{BS97} where the authors show that, in a radially symmetric setting,
the sharp interface problem of a suitable scaling of \eqref{CAC} is a nonlocal motion by mean curvature. Moreover, they also prove that both \eqref{CAC} and the Cahn-Hilliard equation can be seen as degenerate limits of the viscous Cahn-Hilliard equation introduced in \cite{NC92}.

The corresponding motion by mean curvature is also analyzed in \cite{CHL2010} under more general assumptions on the evolving surface
In the quoted contributions, the mixing entropy is approximated with a smooth function defined all over $\mathbb{R}$. In this case, the potential is called regular. In particular, the double-well potential $W$ is usually represented by a fourth-order polynomial.
However, in presence of the nonlocal constraint, one cannot ensure that $\varphi$ takes its values in the physical range $[0,1]$ (see, however, \cite{G97} for an alternative model).
Instead, if the mixing entropy is not approximated by a smooth function, then the image of $\varphi$ is always contained in $[0,1]$. Well-posedness issues in the case of a smooth $W$ are standard. However, if $W$ given by \eqref{dwpotential} then proving the existence of sufficiently regular global solutions is less trivial because $S^\prime$ is singular at the endpoints and cannot be controlled by $S$ like a polynomial. In this
case, it would be nice to show that $\varphi$ stays uniformly away from $0$ and $1$, that is, if \an{the} strict separation property holds, then $S^\prime$ would be globally Lipschitz and the analysis would simplify a lot
(see, for instance, \cite{GGG} and references therein for the Cahn-Hilliard equation in two dimensions, see also \cite{P} for the case of three dimensions). In the non-conserved case, the strict separation is trivial for regular potentials and a bit less straightforward for logarithmic
type potentials like \eqref{dwpotential} (see \cite[Thm.2.3]{GPS2006}).
Concerning \eqref{CAC}, it has been proven its instantaneous validity in dimension two (see \cite{GGW}), while in dimension three the proof was given assuming that the initial datum is strictly separated (see \cite{GGP}).
Observe that the strict separation property combined with the uniqueness of a solution $\varphi$, allows to view the solution itself as the solution to a similar problem where $S$ is replaced by a smooth approximation, defined on the whole real line, which coincides with $S$ on the interval $[\delta, 1-\delta]$ and $\delta\in (0,1)$ is such that $\varphi \in [\delta,1-\delta]$. In other words, the validity of the strict separation can be interpreted as a rigorous justification of the regular approximation of a singular (e.g. Flory-Huggins) potential.

In this paper we want to reconsider these issues and say more for a multi-component version of \eqref{CAC}. In many applications, it is important to account for the presence of multiple interacting species
(see, for instance, \cite{BGSS2013,BW2005,KL2017,KK2006,GNSW2008,WS2021,ZL2022} and references therein, see also \an{\cite{M2023} and its references for the motion by mean curvature in the non-conserved case and} \cite{RB2023} for the importance of the Flory-Huggins potential). Nevertheless, to our knowledge, a comprehensive theoretical analysis of multi-component conserved Allen-Cahn equations is missing. Nonetheless,
it is worth recalling \cite{ST2009,Fu2017} and their references for non-conserved stationary problems with a regular potential. Moreover, we mention that a rigorous solution to the so-called Keller-Rubinstein-Sternberg problem on the motion by curvature has recently been given in \cite{FLWZ2023} (see also its references).
On the contrary, multi-component Cahn-Hilliard equations have been analyzed long ago in the pioneering paper \cite{EL} (see also \cite{GGPS} and its references for further results and recent developments).
As we shall see, one of the advantages (and our main result) is the fact that any weak solution becomes instantaneously strong and strictly separated also in dimension three, while this property is known only
in dimension two for the corresponding multi-component Cahn-Hilliard equation. This regularization allows us to investigate the longtime behavior of solutions in some details, that is, we prove the existence of a
global and an exponential attractor. Also, we can show that any weak solution converges to a single stationary state.
The present analysis can also be viewed as a first step towards the analysis of multi-component Navier-Stokes-Allen-Cahn systems (see, for instance, \cite{ATT2023,Y2021,YH2022}, see also \cite{GGW,LOCY2022} and references therein for binary fluids). We also believe that this contribution
is a significant addition to \cite[Sec.9]{KSZ}.

In the multi-component case, we denote by $\mathbf{u}: \Omega \times (0,T) \to \mathbb{R}^N$ the vector-valued function of concentration species whose components must satisfy the constraint
\begin{equation}
\label{constraint}
\sum_{i=1}^{N} u_{i}\equiv 1.
\end{equation}
The free energy density takes the form
			\begin{align}
			\Psi(\textbf{s})=\sum_{i=1}^N \psi(s_i)-\frac 1 2 \textbf{s}^T\textbf{As},
			\label{pssi}
			\end{align}
            where $\mathbf{A}$ is a constant symmetric $N\times N$ matrix with the largest
			eigenvalue $\lambda _{\mathbf{A}}>0$. Concerning $\psi$, here we are mainly interested in the Boltzmann-Gibbs mixing entropy, namely,
		\begin{equation}
        \label{potential}
		\Psi ^{1}(\mathbf{u}):=\theta \sum_{i=1}^{N}{u}_{i}\ln u_{i} =\sum_{i=1}^{N}\psi ({u}_{i}),
		 \end{equation}%
		where $\theta >0$ is the absolute temperature \an{of} the mixture. However, our framework also includes many other (physically
		relevant) entropy functions $\Psi ^{1}:\left[ 0,1\right] \rightarrow \mathbb{R}_{+}$ (see \cite{GGG,GGPS}).
	The free energy $\mathcal{E}$ is thus defined as	
	$$
	\mathcal{E}(\uu) := W(\uu) +\frac \gamma 2 \int_\Omega \vert\nabla\uu\vert^2dx,
	$$
    where
    $$
    W(\uu)=\int_\Omega \Psi(\uu)dx.
    $$
    Setting
    $$
    \mu^0_{i} = \frac{\delta W  }{\delta {u_i}} = \Psi_{, u_i}, \qquad i=1,\ldots ,N,
    $$
    the vector $\boldsymbol{\mu}^0$ is the chemical potential without capillarity and
    $$
    \boldsymbol{\mu} = -\gamma\Delta\mathbf{u} + \boldsymbol{\mu}^0
    $$
    is the chemical potential.

    Summing up, arguing as in \cite{EL} for the Cahn-Hilliard case, the goal of this work is to study the following initial and boundary value problem
			\begin{equation}
			\begin{cases}
			\partial _{t}\mathbf{u}+\boldsymbol{\alpha } (\ww-\overline{\ww})=\textbf{0},\quad \text{ in
			}\Omega \times (0,T), \\
			\ww = \mathbf{P}\boldsymbol{\mu} = -\gamma \Delta \mathbf{u} +\mathbf{P}\boldsymbol{\mu}^0  ,\quad \text{
				in }\Omega \times (0,T), \\
			\nabla {u}_{i}\cdot \mathbf{n}=0,\quad \text{ on }\partial \Omega \times (0,T),  \quad i=1,\ldots ,N, \\
			\mathbf{u}(0)=\uu_0,\quad \text{ in }\Omega .%
			\end{cases}
			\label{syst}
			\end{equation}
			The (constant) mobility matrix $\boldsymbol{\alpha }$ is
			a symmetric, positive semidefinite $N\times N$ matrix such that its
			kernel is given by $span\{\mathbf{\boldsymbol{\zeta}}\}$ (where $\zeta_{i}=1$,
			for $i=1,\ldots ,N$). Here $\mathbf{P}$ is defined as follows (see also the next section)
	\begin{equation}
    \label{proj}
	\left( \mathbf{P}\mathbf{v}\right) _{l}=v_{l} - \frac{1}{N}\sum_{m=1}^{N} v_{m}, \qquad l=1,\ldots ,N.
	\end{equation}
Then it is easy to check that, formally, a solution to the above problem with $\mathbf{P}\Delta\mathbf{u}=\Delta \mathbf{P}\mathbf{u}$ in place of $\Delta\uu$ satisfies \eqref{constraint} if the initial datum does, using in \eqref{syst}$_1$
\an{the property
$$
\sum_{l=1}^{N} \left( \mathbf{P}\mathbf{v}\right) _{l} =0
$$
and the fact that $\sum_{i=1}^N\al_{ij}=0$ for any $j=1,\ldots,N$  (recall also that $\al$ is symmetric).}
Therefore $\mathbf{P}\Delta \mathbf{u} = \Delta\mathbf{u}$.

The plan of the paper goes as follows. In the next section we introduce the notation, the functional setup, and some basic assumptions on the mobility matrix $\al$. Also, we discuss
the basic assumptions on the potential (more general than \eqref{pssi}-\eqref{potential}) and its regularization. The main results are stated in Section~3 and the last subsection contains the proof
of the convergence to a single equilibrium. The proofs of the well-posedness and regularity results, including the strict separation property, can be found in Section~4.
The existence of the global attractor and of an exponential attractor are proven in Sections~5 and 6,
respectively.

	\section{The mathematical framework}
	
	\label{setting} The (real) Sobolev spaces are denoted as usual by $W^{k,p}(\Omega )$%
	, where $k\in \mathbb{N}$ and $1\leq p\leq \infty $, with norm $\Vert \cdot
	\Vert _{W^{k,p}(\Omega )}$. The Hilbert space $W^{k,2}(\Omega )$ is denoted
	by $H^{k}(\Omega )$ with norm $\Vert \cdot \Vert _{H^{k}(\Omega )}$.
	Moreover, given a space $X$, we denote by $\mathbf{X}$ the space of vectors
	of three components, each one belonging to $X$. We then denote by $(\cdot
	,\cdot )$ the inner product in $L^{2}(\Omega )$ and by $\Vert \cdot \Vert $
	the induced norm. We indicate by $(\cdot ,\cdot )_{X}$ and $\Vert \cdot
	\Vert _{X}$ the canonical inner product and its induced norm in a generic (real) Hilbert
	space $X$, respectively. 
	Further, we introduce the
	affine hyperplane
	\begin{equation}
	\Sigma :=\left\{ \mathbf{c}^{\prime }\in \mathbb{R}^{N}:%
	\sum_{i=1}^{N}c_{i}^{\prime }=1\right\} ,  \label{sigma}
	\end{equation}%
	the Gibbs simplex
	\begin{equation}
	\mathbf{G}:=\left\{ \mathbf{c}^{\prime }\in \mathbb{R}^{N}:%
	\sum_{i=1}^{N}c_{i}^{\prime }=1,\quad c_{i}^{\prime }\geq 0,\quad i=1,\ldots
	,N\right\} ,  \label{Gibbs}
	\end{equation}%
	and the tangent space to $\Sigma $
	\begin{equation}
	T\Sigma :=\left\{ \mathbf{d}^{\prime }\in \mathbb{R}^{N}:%
	\sum_{i=1}^{N}d_{i}^{\prime }=0\right\} .  \label{sigma2}
	\end{equation}%
	We introduce the following notation:
	{\color{black}\begin{equation*}
\an{\mathbf{H}_{0}}:=\{\fff\in \mathbf{L}^2(\Omega):\ \int_{\Omega
}\fff\\dx =\mathbf{0}\text{ and } \fff(x)\in
T\Sigma \;\text{ for a.a. }x\in \Omega \},
\end{equation*}%
\begin{equation*}
\an{\widetilde{\HHH}_{0}}:=\{\fff\in \mathbf{L}^2(\Omega):\fff(x)\in T\Sigma \,\text{ for a.a. }x\in
\Omega \},
\end{equation*}%
\begin{equation*}
\VVV_{0}:=\{\fff\in \mathbf{H}^1(\Omega):\int_{\Omega }%
\fff\\dx =\mathbf{0}\text{ and }\fff(x)\in T\Sigma\, \text{ for a.a. }x\in
\Omega \},
\end{equation*}%
\begin{equation*}
\widetilde{\VVV}_{0}:=\{\fff\in \mathbf{H}^1(\Omega):\ \fff%
(x)\in T\Sigma \, \text{ for a.a. }x\in \Omega \}.
\end{equation*}
}	%
	\an{Notice that the spaces above are still Hilbert spaces with the same
	inner products given in $\mathbf{L}^2(\Omega)$ for the first two, and $\mathbf{H}^1(\Omega)$, for the others. We also have (see \cite{GGPS}) the Hilbert triplets $\VVV_0\hookrightarrow \HHH_0\hookrightarrow \VVV_0'$ and $\widetilde{\VVV}_0\hookrightarrow \widetilde{\HHH}_0\hookrightarrow \widetilde{\VVV}_0'$.}
	
	Recalling \eqref{proj}, we now define rigorously the Euclidean projection $\mathbf{P}$ of $\mathbb{R}^{N}$
	onto $T\Sigma $, which is, for $l=1,\ldots ,N$,
	\begin{equation*}
	\left( \mathbf{P}\mathbf{v}\right) _{l}=\left( \mathbf{v}-\left( \frac{1}{N}%
	\sum_{i=1}^{N}v_{i}\right) \boldsymbol{\zeta}\right) _{l},
	\end{equation*}%
	where $\boldsymbol{\zeta}:=(1,1,\ldots,1)$.
	Notice that the projector $\mathbf{P}$ is also an orthogonal $\mathbf{L}%
	^{2}(\Omega )$-projector, being symmetric and idempotent. We now assume that $\boldsymbol{\alpha }$
	is positive definite over $T\Sigma $. This will constitute the main
	assumption on the mobility matrix in this contribution, since it is enough to prove the existence of weak (and strong) solutions. Nevertheless it is not enough to show the validity of a continuous dependence estimate. Thus we need a second assumption (see (\textbf{M1})). More precisely, we assume that:
	
	\begin{itemize}
		\item[(\textbf{M0})] there exists $l_{0}>0$
		such that
		\begin{equation}
		\boldsymbol{\alpha }\boldsymbol{\eta }\cdot \boldsymbol{\eta }\geq l_{0}%
		\boldsymbol{\eta }\cdot \boldsymbol{\eta },\quad\forall\,\boldsymbol{%
			\eta }\in T\Sigma ;  \label{pos}
		\end{equation}
		\item[(\textbf{M1})] $\al\in\R^{N\times N}$ has the structure
		\begin{align}
		\al=\begin{bmatrix}
	A &B&\ldots&B\\
		B &  A&\ldots&B\\
		\vdots&\vdots&\vdots&\vdots\\
		B&\ldots&\dots&A
		\end{bmatrix},
		\label{general1}
		\end{align}
		where $A>0$ and $A+(N-1)B=0$, so that $B=-\frac A{N-1}<0$.
		\end{itemize}
\begin{remark}
	Note that assumption (\textbf{M1}) can be also rewritten as follows: there exists $\xi>0$ such that
	\begin{align}
	\al=\xi\begin{bmatrix}
	N-1 &-1&\ldots&-1\\
	-1 &  N-1&\ldots&-1\\
	\vdots&\vdots&\vdots&\vdots\\
	-1&\ldots&\dots&N-1
	\end{bmatrix}.
	\label{general2}
	\end{align}
	A matrix of this kind is the natural extension to the case $N>2$ of the admissible matrix $\al$ when $N=2$, which has necessarily the form \eqref{general2}, as one can easily verify. Observe that when $\xi=1$ the matrix $\al$ is simply the representative matrix of the projector $\textbf{P}$, i.e., the identity operator over the space $T\Sigma$.
	We also point out that $\al$ is for sure positive semidefinite and satisfies \eqref{pos}, since it has a zero simple eigenvalue corresponding to the eigenspace $T\Sigma^\perp$, whereas on $T\Sigma$ we see by Lemma \ref{lem} below (with $\mathbf{C}$ equal to the $N\times N$ identity matrix) that $\al$ is positive definite. In particular, one could show that the eigenvalues of $\al$ are $\lambda_1=0$ (corresponding to the eigenvector $(1,1,\ldots,1)$), and $\lambda_i=\xi N$, for $i=2,\ldots,N$, whose eigenspace is clearly $T\Sigma$.
\end{remark}	
Next, we define the set
	\begin{equation}
	\mathcal{K}:=\left\{ \boldsymbol{\eta }\in \mathbf{H}^{1}(\Omega
	):\sum_{i=1}^{N}\eta _{i}=1,\quad\eta _{i}\geq 0,\quad \forall i=1,\ldots
	,N\right\} .  \label{K}
	\end{equation}%
	For the sake of simplicity we will adopt the compact notation $\mathbf{v}%
	\geq k$, with $\mathbf{v}\in \mathbb{R}^{N}$ and $k\in \mathbb{R}$ to
	indicate the relations $v_{i}\geq k$, $i=1,\ldots ,N$.
	
	Recalling \eqref{potential}, we now set
	\begin{equation}
	\left( \boldsymbol{\phi }(\mathbf{u})\right) _{i}=\phi ({u}%
	_{i}):=\psi ^{\prime }({u}_{i}),\quad i=1,\ldots ,N.
	\label{boldphi}\end{equation}%
	In order to include a large admissible class of entropy functionals in (\ref{pssi}), we suppose that%
	\begin{equation*}
	\psi \in C\left[ 0,1\right] \cap C^{2}(0,1]
	\end{equation*}%
	has the following properties:
	
	\begin{itemize}
		\item[(\textbf{E0})] $\psi^{\prime \prime }(s)\geq \zeta >0,$ for all $s\in
		(0,1];$
		
		\item[(\textbf{E1})] $\lim_{s\rightarrow 0^{+}}\psi^{\prime }\left( s\right)
		=-\infty ;$
		
		\item[(\textbf{E2})] $\lim_{s\to0^+} \left(\psi^\prime(s-2s^2)-\psi^\prime(2s^2)\right)=+\infty$.
	\end{itemize}
	

 \an{As in \cite{GGPS}, we also extend $\psi(s)=+\infty ,$ for any $s\in (-\infty ,0)$, and extend $\psi$ for all $s\in \lbrack 1,\infty )$ so that $\psi$ is a $C^2$ function on $(0,+\infty)$ and  $(\textbf{E0})$ holds for any $s>0$. In particular we define
\begin{align}
\psi(s):=As^3+Bs^2+Ds,\quad\text{ for all }s\geq 1,
    \label{psiext}
\end{align}
with
\begin{align*}
    \begin{cases}
      A=\psi(1)-\psi'(1)+\frac 1 2 \psi''(1),\\
      B=-3\psi(1)+3\psi'(1)-\psi''(1),\\
      D=3\psi(1)-2\psi'(1)+\frac 1 2\psi''(1).
    \end{cases}
\end{align*}
}

 We
	refer the reader to \cite[Section 6.3]{GGG} for some other important classes
	of mixing potentials that are singular at $0$.
	Furthermore, following the general scheme developed in \cite[Section 3.1]%
	{GGG2}, by (\textbf{E0})-(\textbf{E1}) we can define an approximation of the
	potential $\psi $ by means of a sequence $\left\{ \psi _{\varepsilon
	}\right\} _{\varepsilon >0}$ of everywhere defined non-negative functions.
	More precisely, let
	\begin{equation}
	\psi _{\varepsilon }(s)=\frac{\varepsilon }{2}|\mathbb{T}_{\varepsilon
	}s|^{2}+\psi (J_{\varepsilon }(s)),\qquad s\in \mathbb{R},\ \varepsilon >0,
	\label{approx}
	\end{equation}%
	where $J_{\varepsilon }=(I+\varepsilon \mathbb{A})^{-1}: \an{\R\to (0,+\infty)}$ is the resolvent
	operator and $\mathbb{T}_{\varepsilon }=\frac{1}{\varepsilon }%
	(I-J_{\varepsilon })$ is the Yosida approximation of $\mathbb{T}\left(
	s\right) :=\psi^{^{\prime }}\left( s\right) ,$ for all $s\in \mathfrak{D}%
	\left( \mathbb{T}\right) =(0,1]$. According to the general theory of maximal
	monotone operators, as already developed in \cite[Section 2]{GGPS}, the following properties hold:
	
	\begin{itemize}
		\item[(i)] $\psi _{\varepsilon }$ is convex and $\psi _{\varepsilon
		}(s)\nearrow \psi(s)$, for all $s\in \mathbb{R}$, as $\varepsilon $ goes to $%
		0$;
		
		\item[(ii)] $\psi _{\varepsilon }^{\prime }(s)=\mathbb{A}_{\varepsilon }(s) $
		and \an{$\phi _{\varepsilon }:=\psi _{\varepsilon }^{\prime }$} is globally Lipschitz with constant $\frac{1}{\varepsilon }$;
		
		\item[(iii)] $|\psi _{\varepsilon }^{\prime }(s)|\nearrow |\psi^{\prime
		}(s)| $ for all $s\in (0,1]$ and $|\psi _{\varepsilon }^{\prime
		}(s)|\nearrow \infty ,$ for all $s\in (-\infty ,0]$, as $\varepsilon $ goes
		to $0$;
		
		
		\item[(iv)] for any $\varepsilon \in (0,1]$, there holds
		\begin{equation*}
		\psi _{\varepsilon }^{\prime \prime }(s)\geq \frac{\zeta}{1+\zeta},\quad
		\text{for all}\ s\in \mathbb{R};
		\end{equation*}
		
		\item[(v)] for any compact subset $M\subset (0,1]$, $\psi _{\varepsilon
		}^{\prime }$ converges uniformly to $\psi^{\prime }$ on $M$;
		
		\item[(vi)] for any $\varepsilon _{0}>0$ there exists $\tilde{K}=\tilde{K}%
		(\varepsilon _{0})>0$ such that
		\begin{equation*}
		\sum_{i=1}^{N}\psi _{\varepsilon }(r_{i})\geq \frac{1}{4\varepsilon _{0}}|%
		\mathbf{r}|^{2}-\tilde{K},\quad \forall \mathbf{r}\in \mathbb{R}^{N},\quad
		\forall 0<\varepsilon <\varepsilon _{0}.
		\end{equation*}
\end{itemize}

		The latter property directly follows from a straightforward adaptation of
		\cite[Lemma 3.11]{GGG2}, which entails that for any $\varepsilon _{0}>0$
		there exists $C=C(\varepsilon _{0})>0$ such that $\psi _{\varepsilon
		}(s)\geq \frac{1}{4\varepsilon _{0}}s^{2}-C$, for any $s\in \mathbb{R}$ and
		any $0<\varepsilon <\varepsilon _{0}$ (see also \cite[Section 2]{GGPS}). Let us now introduce
		\begin{equation*}
		\Psi _{\varepsilon }(\mathbf{r}):=\sum_{i=1}^{N}\psi _{\varepsilon }(r_{i})-%
		\frac{1}{2}\mathbf{r}^{T}\mathbf{Ar}=\Psi _{\varepsilon }^{1}(\mathbf{r})-%
		\frac{1}{2}\mathbf{r}^{T}\mathbf{Ar},
		\end{equation*}%
  \an{where, as presented in the Introduction, $\mathbf{A}$ is a symmetric $N\times N$ matrix with $\lambda _{\mathbf{A}}>0$ as the largest eigenvalue. }
		We thus have that for any $\varepsilon _{0}>0$ sufficiently small there
		exist $K=K(\varepsilon _{0})>0$ and $C=C(\varepsilon _{0})>0$, $%
		C(\varepsilon _{0})\nearrow +\infty $ as $\varepsilon _{0}\rightarrow 0$,
		such that
		\begin{equation*}
		\Psi _{\varepsilon }(\mathbf{r})\geq C(\varepsilon _{0})|\mathbf{r}%
		|^{2}-K,\quad \forall \mathbf{r}\in \mathbb{R}^{N},\quad \forall
		\varepsilon \in(0,\varepsilon _{0}).
		\end{equation*}%
		In particular, this comes from the fact that $-\frac{1}{2}\mathbf{r}\cdot
		\mathbf{Ar}\geq -\frac{\lambda _{\mathbf{A}}}{2}|\mathbf{r}|^{2}$ and $%
		\varepsilon _{0}$ has to be small enough so that, e.g., $C(\varepsilon _{0})=%
		\frac{1}{4\varepsilon _{0}}-\frac{\lambda _{\mathbf{A}}}{2}>0$.
			
\begin{remark}
	We point out that, differently from the standard assumptions on $\psi$ (see, e.g., \cite{GGG,GGW}) here we do not need the assumption \begin{equation*}
	\phi^{{\prime }}\left( s\right) =\psi^{\prime \prime }(s)\leq C\mathrm{e}%
	^{C|\psi^{\prime }(s)|^{\beta}},\text{ for all }s\in (0,1], \; \beta \in [1,2),
	\end{equation*}
since to deduce the validity of the instantaneous strict separation property we will make use only of assumptions (\textbf{E0})-(\textbf{E2}). Clearly the logarithmic potential \eqref{pssi}-\eqref{potential} satisfies assumptions (\textbf{E0})-(\textbf{E2}) and is then included in our analysis. Indeed also assumption (\textbf{E2}) certainly holds for the logarithmic potential since $\psi^\prime(s)=\theta(\ln(s)+1)$ and thus $ \psi^\prime(s-2s^2)-\psi^\prime(2s^2)=\theta(\ln(s-2s^2)-\ln(2s^2))=\theta\ln\left(\frac 1 {2s}-1\right)\to +\infty$ as $s\to 0^+$. Moreover, it seems that if we consider potentials exploding at infinity more slowly than the logarithm then (\textbf{E2}) is not satisfied. Indeed, if, for instance, we consider $\psi^\prime(s)=-\ln(\vert\ln(s)\vert)$ then we get
$\psi^\prime(s-2s^2)-\psi^\prime(2s^2)=-\ln(\vert\ln(s-2s^2)\vert)+\ln(\vert\ln(2s^2)\vert)\to \ln(2)$ as $s\to0$.
\end{remark}	

\section{Main results}
This section is divided into several subsections according to the nature of the results.

\subsection{Well-posedness and regularity}

\label{well-posedness}
We first deal with well-posedness and regularity (see \cite{GGPS} for the multi-component Cahn-Hilliard system). We have

\begin{theorem}\label{thm}
	\begin{enumerate}
\item 	Assume (\textbf{M0}) and (\textbf{E0})-(\textbf{E1}), and let $%
	\uu_0\in \mathcal{K}$. Suppose that
	\begin{equation}
	\label{avesep}
	\delta _{0}<\overline{\mathbf{u}}_{0},
	\end{equation}
	for some $0<\delta _{0}<\frac{1}{N}$.
	Then there exists a pair $(\mathbf{u},\ww)$ defined on $[0,\infty)$, called global finite energy solution to \eqref{syst},
	such that, for any $T>0$, it has the following properties
	\begin{align*}
	& \mathbf{u}\in C([0,T];\textbf{L}^2(\Omega))\cap L^{\infty }(0,T;\mathbf{H}^{1}(\Omega ))\cap L^2(0,T;\mathbf{H}^2(\Omega)), \\
	& \partial _{t}\mathbf{u}\in L^{2}(0,T;\mathbf{L}^{2}(\Omega ) ),\\
	& \ww\in L^{2}(0,T;\mathbf{L}^{2}(\Omega )),
	\\
	& \phi (u_i)\in L^{2}(0,T;L^{2}(\Omega )),\quad \an{i=1,\ldots, N},
	\end{align*}
    and satisfies
	\begin{align}
	& \mathbf{u}(\cdot ,t)\in \mathcal{K},\quad \overline{\mathbf{u}}(\cdot ,t)\equiv \overline{\mathbf{u}}_{0}, \quad \text{ for a.a. }t\in
	(0,T),\label{eq1bis}
	\\
	& 0<\mathbf{u}(x,t)<1 \quad \text{ for a.a. } (x,t)\in \Omega \times (0,T),  \label{eq2}
	\\
	& \partial_t\mathbf{u} + \boldsymbol{\alpha }( \ww-\overline{\ww})=\mathbf{0},\quad
   \text{ a.e. in } \Omega\times (0,T),  \label{phi1} \\
	& \ww=\mathbf{P}(-\mathbf{Au}+\boldsymbol{\phi }(\uu))-\gamma\Delta\uu \quad\text{ a.e. in }\Omega\times(0,T),\label{mu1}\\
    &\partial_\textbf{n}\uu=\mathbf{0} \quad \text{ a.e. on } \partial\Omega\times (0,T), \label{nbc}\\
    & \mathbf{u}(0)=\uu_0 \quad \text{ a.e. in } \Omega.  \label{eq1}
 	\end{align}
    Moreover, the following energy inequality holds
	\begin{equation}
	\label{ident}
	\mathcal{E}(t)+\int_0^t(\pmb\alpha  (\ww(s)-\overline{\ww}(s)), \ww(s)-\overline{\ww}(s)
	)ds\leq 	\mathcal{E}(0),\quad \forall\, t\in \lbrack 0,\infty].
	\end{equation}%
	If, in addition, (\textbf{M1}) holds and  $\overline{\uu}_0^1=\overline{\uu}_0^2$, then two solutions ${\uu}_1$, ${\uu}_2$ are such that
	\begin{align}
	\Vert \uu_1(t)-\uu_2(t)\Vert\leq C\Vert \uu_0^1-\uu_0^2\Vert,\quad \forall\, t\in[0,T],
	\label{contdep}
	\end{align}
    for some $C=C(T)>0$ and uniqueness follows.
	
\item Assume (\textbf{M0}) and (\textbf{E0})-(\textbf{E1}) and let $\uu_0\in \mathcal{K}\cap \mathbf{H}^2(\Omega)$ is such that $\partial_\textbf{n}\uu_0=\mathbf{0}$ almost everywhere on $\partial\Omega$, and $\an{\phi}({u}_{0,i})\in L^2(\Omega)$ for any $i=1,\ldots,N$, then there is a finite energy solution $\uu$ such that, for any $T>0$,
\begin{align}
& \notag\mathbf{u}\in C([0,T];\textbf{H}^1(\Omega))\cap L^{\infty }(0,T;%
\mathbf{H}^{2}(\Omega )), \\
& \label{ddt}\partial _{t}\mathbf{u}\in L^{2}(0,T;\mathbf{H}^{1}(\Omega )),\\
& \ww\in L^{\infty}(0,T;\mathbf{L}^{2}(\Omega )) \cap L^{2}(0,T;\mathbf{H}^{1}(\Omega )),\label{mu}
\\
& \phi (u_i)\in L^{\infty}(0,T;L^{2}(\Omega )),\quad i=1,\ldots,N.\label{phi0}
\end{align}
Moreover, $\uu$ satisfies the energy identity
	\begin{equation}
\label{ident0}
\ddt\mathcal{E}+(\pmb\alpha  (\ww-\overline{\ww}), \ww-\overline{\ww}
)=0,\quad \text{ for a.a. }\, t\in \lbrack 0,\infty).
\end{equation}%
\item Let all the above assumptions hold along with (\textbf{E2}) and suppose that $\uu_0$ is strictly separated, i.e., there exists ${\delta}_0\in (0,\tfrac 1 N)$
such that $\delta_0<\uu_0$ everywhere in $\overline{\Omega}$, then the (unique) strong solution $\uu$ is strictly separated as well, i.e., for any $T>0$ there exists $\delta=\delta(\tau,T)\in(0,\tfrac 1 N]$
such that \begin{align}
\delta\leq\mathbf{u}(x,t),\quad \forall (x,t)\in\overline{\Omega}\times[0,T].
\label{sep}
\end{align}
\end{enumerate}
\end{theorem}

\begin{remark}
On account of \eqref{mu}, one could also prove that $\boldsymbol{\phi}(\uu)\in L^2(0,T;\textbf{L}^p(\Omega))$, \an{where $\boldsymbol{\phi}$ is defined in \eqref{boldphi}}, and $\mathbf{u}\in L^{2}(0,T;\WWW^{2,p}(\Omega ))$ where $p=6$ if $d=3$, while $p\in [2,\infty)$ if $d=2$, by slightly adapting part of the proof of \cite[Thm.3.1]{GGPS} (which is performed for the $L^\infty$-in-time case). Again the main issue is the presence of the projector $\mathbf{P}$ in the definition of $\ww$ (cf. \cite[Cor.1]{Abels2009} for the scalar case).
\end{remark}

\begin{remark}
	\label{control} Notice that \eqref{avesep}
	implies that there exists $\rho>0$ such that $\rho<\overline{u}_{0,i}<1-\rho$
	for any $i=1,\ldots,N$. Indeed, we have, for any $i=1,\ldots,N$,
	\begin{equation*}
	\delta_0<\min_{j=1,\ldots,N}{\overline{u}_{0,j}}\leq \overline{u}_{0,i}=1-\sum_{j\not= i}{\overline{u}}_{0,j}\leq 1-(N-1)\min_{j=1,\ldots,N}{\overline{u}_{0,j}}<1-(N-1)\delta_0,
	\end{equation*}
	and thus we can choose, e.g., $\rho=\delta_0$, being $N\geq 2$.
\end{remark}

\begin{remark}
	Arguing as in \cite[Prop.2.1]{EL}, we easily obtain that $%
	\sum_{i=1}^{N}u_i=1$ and $\sum_{i=1}^{N}w_{i}=0$.
	Moreover, by choosing $\pmb\eta \equiv \boldsymbol{\eta}_{i}$, $\boldsymbol{\eta}_{i}$
	being the $i$-th unit vector, we get that the total mass of each component $u_{i}$ is preserved, i.e.,
	\begin{equation*}
	\overline{\mathbf{u}}(t)\equiv \overline{\mathbf{u}}_{0},\quad \forall\,t\geq 0.
	\end{equation*}
\end{remark}

\begin{remark}
	From point (2) of Theorem \ref{thm}, we deduce that $\nabla u_i\cdot \textbf{n}\in C_w([0,T];H^\frac 1 2(\partial\Omega))$.  Thus
    $\nabla u_i\cdot \textbf{n}=0$ for any $t\in[0,T]$ almost everywhere on $\partial\Omega$, for $i=1,\ldots,N$. Furthermore, since we also have
	$$
	\Vert \uu\Vert_{L^\infty(0,T;\mathbf{H}^2(\Omega))}\leq C(T),
	$$
	being $\Vert \uu(\cdot)\Vert_{\mathbf{H}^2(\Omega)}$ lower semicontinuous, we get
	\begin{equation}
    \label{H2bound}
	\Vert \uu(t)\Vert_{\mathbf{H}^2(\Omega)}\leq C(T), \quad\forall\,t\in [0,T].
	\end{equation}
	 \label{regularity}
\end{remark}

\begin{remark}
	Recalling point (3) of Theorem \ref{thm}, observe that \eqref{constraint} and \eqref{sep} imply the existence of $%
	\delta_1:=(N-1)\delta>0$ such that $\mathbf{u}\leq1-\delta_1$
	almost everywhere in ${\Omega}\times[0,T]$, i.e., each component is
	strictly separated from the pure phases $0$ and $1$. Moreover, property \eqref{sep} holds on $\overline{\Omega}\times[0,T]$ since from its proof (see Section \ref{exi}) we deduce that, for any $t\in[0,T]$,
	$$\uu(t)\geq \delta\quad\text{ a.e. in }\Omega.$$
	Then, by Remark \ref{regularity} we know that $\uu(t)\in \mathbf{H}^2(\Omega)\hookrightarrow C(\overline{\Omega})$ for any $t\in[0,T]$, implying that
	$$\uu(x,t)\geq \delta,\quad\forall (x,t)\in\overline{\Omega}\times[0,T].$$
	
\end{remark}
\begin{remark}
The quantity $%
	\delta >0$ in the separation property only depends on the initial data through the initial data energy $\mathcal{E}(0)$, $\overline{\uu}_0$, $\delta_0$ and $\Vert \uu_0\Vert_{\mathbf{H}^2(\Omega)}$. The same goes, but  $\delta_0$, for all the constants involved in the regularity estimates of point (2) of the Theorem \ref{thm}. 
	\label{dependence2}
\end{remark}
\an{\begin{remark}
    As it will be clear from the proof (see also Remark \ref{N2}), in the case $N=2$ assumption (\textbf{E2}) is not needed to prove \eqref{sep}. This agrees with the result obtained in \cite{GGP} for binary mixtures.
    \label{N2bis}
\end{remark}}
On account of the dissipative nature of the system, we have the following uniform control of the energy $\mathcal{E}$.

\begin{theorem}
	\label{exp} Let the assumptions of Theorem \ref{thm}, point (1), hold. Then the energy of solution $\mathbf{u}$ satisfies the following inequality
	\begin{equation}
	\label{dissineq}
	\mathcal{E}(t)\leq C_{1}e^{-\omega t}\mathcal{E}(0)+C_{2},\quad \forall
	t\geq 0,
	\end{equation}%
	where $C_{1},C_{2}>0$   depend on $\Omega$, $\al$, ${\Psi}$, and $\overline{\mathbf{u}}_{0}$, while $\omega >0$ is a universal
	constant.
\end{theorem}

We can prove that any weak solution given by Theorem \ref{thm} instantaneously regularizes.
Thanks to this, we can show the instantaneous strict separation property in both dimensions two and three.
This means that, for any $\tau>0$, there exists $0<{\delta}=\delta(\tau)<\tfrac{1}{N}$
such that $\delta\leq\mathbf{u}$ almost everywhere in $\Omega%
\times[\tau,\infty)$. Again, notice that this implies the existence of $%
\delta_1:=(N-1)\delta>0$ such that $\mathbf{u}\geq1-\delta_1$
almost everywhere in ${\Omega}\times[\tau,\infty)$, i.e., each component is
strictly separated from both the pure phases $0$ and $1$. More precisely, the following result holds
\begin{theorem}
	\label{thm2} Let the assumptions of Theorem \ref{exp} hold, together with (\textbf{M1}) and (\textbf{E2}). Then the energy solution
	$(\mathbf{u},\ww)$, defined for all $t\geq 0$, is such that, for any $\tau >0$,
	\begin{align}
	& \mathbf{u}\in C([\tau,\infty );\mathbf{H}^{1}(\Omega ))\cap L^{\infty }(\tau,\infty;\mathbf{H}^{2}(\Omega )),\label{regg0} \\
	& \partial _{t}\mathbf{u}\in L^{2}(t,t+1;\mathbf{H}^{1}(\Omega )),\quad
	\forall\, t\geq \tau ,  \label{regg} \\
	& \ww\in L^{\infty }(\tau ,\infty ;\mathbf{L}^{2}(\Omega )), \quad
	\forall\, t\geq \tau, \label{reggbis}
	\\
	& \phi (u_i)\in L^{\infty }(\tau ,\infty ;L^{2}(\Omega )),\quad
	i=1,\ldots ,N. \label{reggtris}
	\end{align}
    Moreover, $\uu$ and $\ww$ are uniformly bounded in the above spaces by positive constants only depending on  $\Omega$, $\al$, $\Psi$, $\overline{\mathbf{u}}_{0}$, and $\mathcal{E}(0)$.
	In particular, the energy identity \eqref{ident0} holds for almost any $t\geq \tau$.
	Moreover, there exists $0<{\delta }=\delta (\tau )\leq \tfrac{1}{N}$
	such that
	\begin{equation}
	\delta \leq\mathbf{u}(x,t), \quad \forall\, (x,t)\in \overline{\Omega}\times \lbrack \tau ,\infty),
	\label{delt}
	\end{equation}%
	i.e., the instantaneous strict separation property holds.
\end{theorem}

\begin{remark}
It is standard to notice that $\uu_\delta:=(\uu-\delta)^-\in C([\tau,\infty);\mathbf{L}^2(\Omega))$. In the proof of the strict separation property (see Section \ref{pr2}) we obtain
 	\begin{equation*}
 \delta \leq\mathbf{u}(x,t), \quad \text{ for a.a. }\, (x,t)\in {\Omega}\times \lbrack \tau ,\infty),
 \end{equation*}%
which then implies $\Vert\uu_\delta(t)\Vert\equiv 0$ for almost any $t\in[\tau,\infty)$, and thus it holds \textit{for any }$t\in[\tau,\infty)$ by continuity. This means that we have
 	\begin{equation}
\delta \leq\mathbf{u}(t), \quad \forall\, t\in[\tau,\infty),\text{ a.e. in }\Omega.
\label{p3}
\end{equation}%
By \eqref{H2bound} and its global nature ensured by Theorem \ref{thm2}, we have that $\uu(t)\in \mathbf{H}^2(\Omega)$, for any $t\in[\tau,\infty)$, entailing that \eqref{p3} holds for any $(x,t)\in \overline{\Omega}\times[\tau,\infty)$.
	 \label{sepimp}
\end{remark}
\an{\begin{remark}
    We point out that, as observed in Remark \ref{N2bis}, assumption (\textbf{E2}) is not needed to prove \eqref{delt} when $N=2$, i.e., for binary mixtures.
\end{remark}}
\subsection{Existence of the regular global attractor}

\label{dynamical} We now define a complete metric space which will be the phase space of
the dissipative dynamical system (see, for instance, \cite{Temam}) associated with \eqref{syst}.
For a given $\mathbf{M}\in \Sigma$,
such that ${M}_{i}\in (0,1),$ for any $i=1,\ldots ,N$, we set
\begin{equation*}
\mathcal{V}_{\mathbf{M}}:=\{\mathbf{u}\in \mathbf{H}^{1}(\Omega ):\ 0\leq
\mathbf{u}(x)\leq 1, \quad \text{ for a.a. }x\in \Omega ,\quad \overline{%
	\mathbf{u}}=\mathbf{M},\quad \sum_{i=1}^{N}u_i = 1\},
\end{equation*}%
endowed with the $\mathbf{H}^{1}$-topology. In particular we consider the one induced by the
equivalent norm $\Vert \mathbf{u}\Vert _{\mathcal{V}_{\mathbf{M}}}=\Vert
\nabla \mathbf{u}\Vert +|\overline{\mathbf{u}}|$. This is a complete metric space.
Thus, under the same assumptions of Theorem \ref{thm2}, we can define a
dynamical system $(\mathcal{V}_{\mathbf{M}},S(t))$ where
\begin{equation*}
S(t):\mathcal{V}_{\mathbf{M}}\rightarrow \mathcal{V}_{\mathbf{M}}, \quad
S(t)\mathbf{u}_{0}=\mathbf{u}(t),\quad \forall \, t\geq 0.
\end{equation*}%
Observe that $S(t)$ satisfies the following properties:

\begin{itemize}
	\item $S(0)=Id_{\mathcal{V}_\textbf{M}}$;
	
	\item $S(t+\tau)=S(t)S(\tau)$, for every $\mathbf{u}_0\in \mathcal{V}_%
	\textbf{M}$;
	
	\item $t\mapsto S(t)\mathbf{u}_0\in C((0,\infty);\mathcal{V}_\textbf{M})$,
	for every $\mathbf{u}_0\in\mathcal{V}_\textbf{M}$;
	
	\item $\mathbf{u}_0\mapsto S(t)\mathbf{u}_0\in C(\mathcal{V}_\textbf{M};%
	\mathcal{V}_\textbf{M})$, for any $t\in[0,\infty)$.
\end{itemize}
In particular, $t\mapsto S(t)\mathbf{u}_0\in C((0,\infty);\mathcal{V}_\textbf{M})$ comes from the instantaneous regularization, so that, for any $\tau>0$, $\uu\in C([\tau,\infty);\mathcal{V}_{\mathbf{M}})$,
whereas the last property can be proved as follows. From \eqref{contdep}
together with the $\mathbf{H}^{2}$-regularity (for any $t>0$) and the
interpolation estimate%
\begin{equation*}
\Vert \cdot \Vert _{\mathbf{H}^{1}(\Omega )}\leq C\Vert \cdot \Vert _{%
	\mathbf{H}^{2}(\Omega )}^{\frac{1}{2}}\Vert \cdot \Vert^{\frac{1}{2}}
\end{equation*}%
we deduce that $\mathbf{u}_{0}\mapsto S(t)\mathbf{u}_{0}\in C(\mathcal{V}_{%
	\mathbf{M}};\mathcal{V}_{\mathbf{M}})$, for any $t\in (0,\infty )$. \an{This is indeed a consequence of \eqref{regg0}, since $\uu\in L^\infty(\tau,\infty;\mathbf{H}^2(\Omega))$ for any $\tau>0$ entails that, given two intial data $\uu_{0,1},\uu_{0,2}\in \mathcal{V}_\mathbf{M}$, for any $t>0$,
 \begin{align*}
     &\Vert S(t)\mathbf{u}_{0,1}-S(t)\mathbf{u}_{0,2}\Vert_{\mathbf{H}^1(\Omega)}\\&\leq C\Vert S(t)\mathbf{u}_{0,1}-S(t)\mathbf{u}_{0,2}\Vert_{\mathbf{H}^2(\Omega)}^{\frac{1}{2}}\Vert S(t)\mathbf{u}_{0,1}-S(t)\mathbf{u}_{0,2}\Vert^{\frac{1}{2}}\\&\leq C(t)\Vert S(t)\mathbf{u}_{0,1}-S(t)\mathbf{u}_{0,2}\Vert^{\frac{1}{2}}\leq C(t)\Vert \mathbf{u}_{0,1}-\mathbf{u}_{0,2}\Vert^{\frac{1}{2}},
 \end{align*}
 where in the last step we also used \eqref{contdep}.
 }
  The case
$t=0$ is trivial.

Furthermore, we recall that the \textit{global attractor} is the unique compact set $%
\mathcal{A}\subset \mathcal{V}_\textbf{M}$ such that

\begin{itemize}
	\item $\mathcal{A}$ is fully invariant, i.e., $S(t)\mathcal{A}=\mathcal{A}$
	for every $t\geq0$;
	
	\item $\mathcal{A}$ is attracting for the semigroup, i.e.,
	\begin{equation*}
	\lim_{t\rightarrow\infty}[\text{dist}_{\mathcal{V}_\textbf{M}}(S(t)B,%
	\mathcal{A})]=0
	\end{equation*}
	for every bounded set $B\subset \mathcal{V}_\textbf{M}$.
\end{itemize}

The dissipative inequality \eqref{dissineq} and the instantaneous regularization of the energy solution allow us to prove

\begin{theorem}
	\label{global} Let the assumptions of Theorem \ref{thm2} hold.  Then the dynamical system $(\mathcal{V}_{\mathbf{M}},S(t))$ admits
	a (unique) connected global attractor $\mathcal{A}\subset \mathcal{V}_{%
		\mathbf{M}}$ which is bounded in $\mathbf{H}^2(\Omega)$.
\end{theorem}

\begin{remark}
	\label{trick}The proof of this result is based on showing that the dynamical
	system $(\mathcal{V}_{\mathbf{M}},S(t))$ admits a compact absorbing set $%
	\mathcal{B}_{0}$ (see Section \ref{pr3} below).
\end{remark}
\subsection{Existence of an exponential attractor}
Thanks to the validity of the strict separation property in dimensions two and three, we can prove the existence of an exponential attractor in dimensions two and three.
We first recall (see, e.g., \cite{Zelik}) that a compact set $\mathcal{M}\subset \mathcal{V}_{%
	\mathbf{M}}$ is an exponential attractor for $( \mathcal{V}_{%
	\mathbf{M}},S(t))$ if
\begin{itemize}
	\item $\mathcal{M}$ is positively invariant, i.e., $S(t)\mathcal{M}\subset \mathcal{M}$ for every $t\geq0$;
	\item $\mathcal{M}$ is exponentially attracting, i.e, there exists $\overline{\omega}>0$ such that
	$$
	\text{dist}_{\mathcal{V}_{%
			\mathbf{M}}}(S(t)\mathcal{B},\mathcal{M})\leq \mathcal{Q}(\Vert\mathcal{B}\Vert_{\mathcal{V}_{%
				\mathbf{M}}})e^{-{\overline{\omega}} t}
	$$
	for every bounded $\mathcal{B}\subset \mathcal{V}_{%
		\mathbf{M}}$, where $\mathcal{Q}(\cdot)$ denotes a generic increasing positive function;
	\item $\mathcal{M}$ has finite fractal dimension in $ \mathcal{V}_{%
		\mathbf{M}}$, where the fractal dimension is defined as $$
	\text{dim}_{ \mathcal{V}_{%
			\mathbf{M}}}(\mathcal{M})=\limsup_{\epsilon\to0}\frac{\log N(\epsilon)}{-\log \epsilon},
	$$
	and $N(\epsilon)$ is the minimum number of $\epsilon$-balls of $ \mathcal{V}_{%
		\mathbf{M}}$ necessary to cover $\mathcal{M}$.
\end{itemize}
Observe that the exponential attractor is not unique, and that, by definition, $\mathcal{A}\subset\mathcal{M}$, so that from the existence result of an exponential attractor we deduce that the global attractor $\mathcal{A}$ is of finite fractal dimension.
We thus have the following
\begin{theorem}
Let the assumptions of Theorem \ref{thm2} hold. Moreover, assume that $\psi\in C^3(0,1]$. Then the dynamical system $(\mathcal{V}_{%
		\mathbf{M}}, S(t))$ possesses an exponential attractor $\mathcal{M}$ which is bounded in $\mathbf{H}^2(\Omega)$. Besides, $\mathcal{A}\subset\mathcal{M}$ has finite fractal dimension in $\mathcal{V}_\mathbf{M}$.
	\label{expatt}
\end{theorem}

\subsection{Convergence to equilibrium}

\label{convec} In this section we discuss the convergence of any weak solution to a single equilibrium. We have all the ingredients to state and prove the result.

We consider the phase space $\mathcal{V}_{\textbf{M}}$ as in the previous section.
Under the assumptions of Theorem \ref{thm2}, we define the $\omega $-limit set $\omega (\uu_0)$ of a given $\uu_0\in \mathcal{V}_{\mathbf{M}}$
\begin{equation*}
\omega (\uu_0)=\{\mathbf{z}\in \mathbf{H}^{2r}(\Omega )\cap
\mathcal{V}_{\mathbf{M}}:\exists\,t_{n} \nearrow \infty \text{ s.t. }\mathbf{%
	u}(t_{n})\rightarrow \mathbf{z}\text{ in }\mathbf{H}^{2r}(\Omega )\},
\end{equation*}%
where $r\in \lbrack \tfrac{1}{2},1)$. In particular, we
fix $r\in (\tfrac{d}{4},1)$. We thus have

\begin{theorem}
	\label{main}  Let the assumptions of Theorem \ref{thm2} hold and suppose, in addition, that $\psi$ is (real) analytic in $(0,1)$.
Then, for any $\uu_0\in \mathcal{V}_\textbf{M}$, it holds $\omega(\mathbf{u%
	}_0)=\{\uu_\infty\}$, where $\uu_\infty
	\in\mathcal{V}_\textbf{M}$ is a solution to
	\begin{equation*}
	\begin{cases}
	-\Delta {\uu}_{\infty }+\mathbf{P}\Psi^1_{,\mathbf{u}}(\mathbf{%
		\uu}_{\infty})=\pmb f, \quad \text{ a.e. in }\Omega , \\
	\partial _{\mathbf{n}}{\uu}_{\infty }=0, \quad \text{ a.e. on }\partial
	\Omega , \\
	\sum_{i=1}^{N} u_{\infty ,i}= 1, \quad \text{ in } \Omega.
	\end{cases}
	\end{equation*}
	with $\pmb{f}=\mathbf{P}\mathbf{A}\uu_\infty+ \overline{\mathbf{P}\Psi_{,\mathbf{u}}({\uu_\infty})}$.
	Moreover, $\overline{\mathbf{u}}_\infty=\mathbf{M}$, there exists $\delta>0$ so
	that
	\begin{equation*}
	\delta<\uu_\infty(x),\quad \forall\,x\in \overline{\Omega },
	\end{equation*}
	and the (unique) weak solution $\mathbf{u}(t)$ is such that
	\begin{align}
	\mathbf{u}(t)=S(t)\uu_0\underset{ t\to \infty}{\longrightarrow}
	\uu_\infty\quad \text{in }\mathbf{H}^{2r}(\Omega),  \quad \forall\,r\in(0,1).\label{sing}
	\end{align}
\end{theorem}
\begin{proof}
	The proof of this Theorem is exactly the same as the one of \cite[Thm. 3.16]{GGPS}. Indeed, the only difference is in the energy estimate given by the application of {\L}ojasiewicz-Simon inequality (see \cite[Sec.7.3]{GGPS}), in which we need to substitute $\nabla\ww$ with $\ww-\overline{\ww}$ (basically we do not need to apply Poincar\'{e}'s inequality, but we keep $\Vert \ww-\overline{\ww}\Vert$ in the inequality for $\mathcal{E}^\prime$).
\end{proof}
Theorem \ref{main} is still valid without assumption (\textbf{E2}). Indeed, in the proof we do not need the instantaneous strict separation property, for which that assumption is essential. It is also worth noticing that, without assuming (\textbf{E2}), by the same proof of \cite[Thm.3.13]{GGPS}, we can show that the asymptotic strict separation property holds, i.e.,
\begin{theorem}
	Let the assumptions of Theorem \ref{thm2} hold \an{except for} (\textbf{E2}). Then, for any $%
	\mathbf{M}\in (0,1)$, $\mathbf{M}\in \Sigma$, and for any initial datum $\uu_0\in \mathcal{V}_{%
		\mathbf{M}}$, there exists $\delta >0$ and $t^{\ast }=t^{\ast }(\mathbf{u}%
	_{0})$ such that the corresponding (unique) solution $\mathbf{u}$ satisfies:
	\begin{equation}
	\delta <\mathbf{u}(x,t),\quad \text{ for any }(x,t)\in \overline{\Omega} \times
	(t^{\ast },+\infty ).  \label{sep1}
	\end{equation}
\end{theorem}

\section{Proofs of Subsection \protect\ref{well-posedness}}

\label{pr1}

Here we collect the proofs of Theorems \ref{thm}, \ref{exp}, and \ref{thm2}.

\subsection{Proof of Theorem \protect\ref{thm}}
\label{exi}
This proof is divided into three parts. We first prove \eqref{contdep} which seems to require (\textbf{M1}). The reason is related to the following
\begin{lemma}
	\label{lem}
	Let (\textbf{M1}) hold. Then there exists $\gamma_N>0$ such that, given any $\mathbf{C}=\text{diag}(c_1,\ldots,c_N)$, with $c_i\geq0$ for any $i=1,\ldots,N$,
	\begin{align}
	\boldsymbol{\zeta }^T({\mathbf{C}}\al)\boldsymbol{\zeta}\geq \gamma_N\left(\min_{i=1,\ldots,N,\ c_i\al_{ii}>0}c_i\al_{ii}\right)\vert \boldsymbol{\zeta }\vert^2\geq 0,
	\label{below}
	\end{align}
	for any $\boldsymbol{\zeta }\in T\Sigma$. In particular, for any $N\geq2$, considering the equivalent structure \eqref{general2}, we have \an{$$\gamma_{N}:=\frac{ N}{N-1}.$$}
\end{lemma}
\begin{remark}
	Notice that, being $\al$ positive semidefinite, $\al_{ii}\geq 0$, for any $i=1,\ldots,N$.
\end{remark}
\begin{remark}
What \an{is }needed to prove \eqref{contdep} is actually \eqref{below}. Nevertheless, the matrix structure \eqref{general1} is the only example case we know that implies \eqref{below}.
	\end{remark}
\begin{proof}
	Note that, for $\boldsymbol{\zeta}\in T\Sigma$, we have $\boldsymbol{\zeta}_N=-\sum_{i=1}^{N-1}\boldsymbol{\zeta}_i $, thus, \an{exploiting the form \eqref{general2} of the matrix $\al$,}
	\begin{align*}
	&\boldsymbol{\zeta}^T{\mathbf{C}}\al\boldsymbol{\zeta}=\xi\left[\boldsymbol{\zeta}_1,\ldots,-\sum_{i=1}^{N-1}\boldsymbol{\zeta}_i \right]\begin{bmatrix}
	c_1\left(\boldsymbol{\zeta}_1(N-1)-\sum_{j\not=1}^{N-1}\boldsymbol{\zeta}_j+\sum_{j=1}^{N-1}\boldsymbol{\zeta}_j\right)\\\vdots\\c_i\left(\boldsymbol{\zeta}_i(N-1)-\sum_{j\not=i}^{N-1}\boldsymbol{\zeta}_j+\sum_{j=1}^{N-1}\boldsymbol{\zeta}_j\right)\\\vdots\\c_N\left(-\sum_{j=1}^{N-1}\boldsymbol{\zeta}_j-(N-1)\sum_{j=1}^{N-1}\boldsymbol{\zeta}_j\right)
	\end{bmatrix}\\&\\&=\xi\left[\boldsymbol{\zeta}_1,\ldots,-\sum_{i=1}^{N-1}\boldsymbol{\zeta}_i \right]\begin{bmatrix}
	c_1N\boldsymbol{\zeta}_1\\\vdots\\c_iN\boldsymbol{\zeta}_i\\\vdots\\-c_NN\sum_{j=1}^{N-1}\boldsymbol{\zeta}_j
	\end{bmatrix}=N\xi\sum_{i=1}^{N-1}c_i\vert\boldsymbol{\zeta}_i\vert^2+\xi Nc_N\left\vert\sum_{j=1}^{N-1}\boldsymbol{\zeta}_j\right\vert^2\\&\geq \xi N \left(\min_{i=1,\ldots,N,\ c_i>0}c_i\right)\vert\boldsymbol{\zeta}\vert^2\geq \gamma_N\left(\min_{i=1,\ldots,N,c_i\al_{ii}>0}c_i\al_{ii}\right)\vert\boldsymbol{\zeta}\vert^2
	\end{align*}
	with \an{$\gamma_N=\frac{ N}{N-1}$}. Thus \eqref{below} holds.
\end{proof}

\textbf{Continuous dependence estimate.} We can now prove \eqref{contdep}.
Let us consider two solutions $\uu^1$ and $\uu^2$ and take the difference between the equations they solve. Taking  $\uu=\uu^1-\uu^2$ as a test function in the resulting equation and
recalling, by mass conservation, that $\overline{\uu}\equiv \textbf{0}$, we deduce (note that $\al \mathbf{P}\Psi^1_{,\uu}(\uu^k)=\al\Psi^1_{,\uu}(\uu^k)$ for $k=1,2$, \an{since $\al\left(\tfrac 1 N\sum_{i=1}^N\psi'(u_i^k)\boldsymbol{\zeta}\right)=\left(\tfrac 1 N\sum_{i=1}^N\psi'(u_i^k)\right)\al\boldsymbol{\zeta}=0$, where $\boldsymbol{\zeta}=(1,\ldots,1)$})
\begin{align}
\dfrac 1 2 \ddt\Vert \uu\Vert^2+\gamma\left(\nabla \uu,\al\nabla \uu\right)+\sum_{i,j=1}^N\left(\al_{ij}(\psi^\prime(u^1_j)-\psi^\prime(u^2_j)),{u}_i\right)=(\al{\mathbf{A}}\uu,\uu).
\label{est4}
\end{align}
\an{Notice that $\overline{\ww}_1-\overline{\ww}_2$ does not appear in \eqref{est4}, since we have
\begin{align*}
&(\al \left((\ww_1- \ww_2)-(\overline{\ww_1-\ww_2})\right),\uu)\\&=(\al \left((\ww_1- \ww_2)-(\overline{\ww_1-\ww_2})\right),\uu-\overline{\uu})\\&=(\al(\ww_1- \ww_2),\uu-\overline{\uu})\\&=(\al(\ww_1- \ww_2),\uu),
\end{align*}
 recalling in the last equality that $\overline{\uu}\equiv\textbf{0}$.
}
Lemma \ref{lem} \an{then} entails
$$
\gamma\left(\nabla \uu,\al\nabla \uu\right)\geq 0.
$$
Then, by the Cauchy-Schwarz inequality,
$$
( \al\mathbf{A}\uu,\uu)\leq C\Vert\uu\Vert^2.
$$
In conclusion, we have
\begin{align*}
\sum_{i,j=1}^N\left(\al_{ij}(\psi^\prime(u^1_j)-\psi^\prime(u^2_j)),{u}_i\right)=\sum_{i,j=1}^N\int_\Omega\int_0^1\al_{ij}\psi^{\prime\prime}(s{u}_j^1+(1-s){u}_j^2){u}_j{u}_idsdx=\int_\Omega \uu^T\al \mathbf{C}\uu dx,
\end{align*}
where $$\mathbf{C}=\text{diag}(c_1,\ldots,c_N):=\text{diag}\left(\int_0^1\psi^{\prime\prime}(s{u}_1^1+(1-s){u}_1^2)ds,\ldots, \int_0^1\psi^{\prime\prime}(s{u}_N^1+(1-s){u}_N^2)ds\right),$$
so that $c_i\geq 0$, for any $i=1,\ldots,N$, by assumption (\textbf{E0}). Observing now that $\uu(x,t)\in T\Sigma$ for almost any $(x,t)\in\Omega\times(0,T)$ and, by symmetry, $\uu^T\al\mathbf{C}\uu=\uu^T\mathbf{C}\al\uu$, thanks to Lemma \ref{lem}, for almost any $(x,t)\in\Omega\times(0,T)$, we have
$$
\uu^T\al \mathbf{C}\uu\geq 0,
$$
so that
$$
\sum_{i,j=1}^N\left(\al_{ij}(\psi^\prime(u^1_j)-\psi^\prime(u^2_j)),{u}_i\right)=\int_\Omega \uu^T\al \mathbf{C}\uu dx\geq 0, \quad\text{ a.e. in } (0,T).
$$
Therefore, from \eqref{est4}, we infer
$$
\dfrac 1 2 \ddt \Vert \uu(t)\Vert^2\leq C\Vert \uu(t)\Vert^2,\quad \text{for a.a. }t\in(0,T),
$$
and the Gronwall Lemma gives \eqref{contdep}.

\medskip\noindent
\textbf{Existence of a solution.} \an{Here we give the details of the Galerkin scheme since in previous
related contributions (see, e.g., \cite{EL}) they are not given.} We consider the approximation %
\eqref{approx}. In particular, for each $\varepsilon>0$
sufficiently small, we set
\begin{equation*}
\boldsymbol{\phi}_\varepsilon(\mathbf{y})=\Psi^1_{\varepsilon,\mathbf{y}}(%
\mathbf{y})=\{\psi^\prime_\varepsilon(y_i)\}_{i=1,\ldots,N}, \quad \forall\, \mathbf{y}\in \mathbb{R}^N.
\end{equation*}
We then fix $0<\varepsilon<\varepsilon_0$ and consider the complete system of $N$-dimensional eigenfunctions $\{\mathbf{e}_i\}_i$ of the problem $-\Delta \mathbf{e}_i=\lambda_i\mathbf{e}_i$, with homogeneous Neumann boundary conditions $\partial_\mathbf{n}\mathbf{e}_i=0$ on $\partial\Omega$ ($\lambda_i$ is the eigenvalue corresponding  to $\mathbf{e}_i$)\an{, subject to the constraints $\overline{\mathbf{e}}_i=0$ and $\sum_{j=1}^N(\mathbf{e}_i)_{j}\equiv 0$.} The family $\{\mathbf{e}_i\}_i$ can be tuned to form an orthogonal basis in \an{$\VVV_0$}, orthonormal in \an{$\HHH_0$} (see also \cite[Appendix 8.1]{GGPS}).
\an{We then set $\mmm:=\overline{\uu}_0$ and} introduce the finite-dimensional spaces
$$\mathbf{V}_n:=\text{span}\{\mathbf{e}_i,i=1,\ldots,n\} \quad \forall\, n\geq 1,$$
and look for a function ${\uu}_{n,\varepsilon}\in \mathbf{V}_n$ of the form
\begin{align*}
{\uu}_{n,\varepsilon}(t)=\sum_{i=1}^{n}{\hat{\alpha}}_i(t)\mathbf{e}_i\in \mathbf{V}_n,
\end{align*}
\an{and for ${\ww}_{n,\varepsilon}\in \widetilde{\VVV}_0$ such that
\begin{align*}
{\ww}_{n,\varepsilon}(t)-\overline{{\ww}}_{n,\varepsilon}(t)=\sum_{i=1}^{n}{{\delta}}_i(t)\mathbf{e}_i\in \mathbf{V}_n,
\end{align*}}
solving the equations
\begin{align}
&(\partial_t{\uu}_{n,\varepsilon},\mathbf{v})+(\al(\ww_{n,\varepsilon}-\overline{\ww}_{n,\varepsilon}),\mathbf{v})=0, \label{ac}\\&
({\ww}_{n,\varepsilon}-\an{\overline{\ww}_{n,\varepsilon}},\mathbf{v})\nonumber\\&=\gamma(\nabla {\uu}_{n,\varepsilon},\nabla\mathbf{v})+\left(\textbf{P}(\boldsymbol{\phi}_\varepsilon({\uu}_{n,\varepsilon}+\an{\mmm})-\mathbf{A}({\uu}_{n,\varepsilon}+\an{\mmm}))-\overline{\textbf{P}(\boldsymbol{\phi}_\varepsilon({\uu}_{n,\varepsilon}+\an{\mmm})-\mathbf{A}({\uu}_{n,\varepsilon}+\an{\mmm}))},\mathbf{v}\right),\label{mum}\\&\label{subs}
\an{\overline{\ww}_{n,\varepsilon}=\overline{\textbf{P}(\boldsymbol{\phi}_\varepsilon({\uu}_{n,\varepsilon}+\an{\mmm})-\mathbf{A}({\uu}_{n,\varepsilon}+\an{\mmm}))}},
\\&
{\uu}_{n,\varepsilon}(0)={\uu}_{n,0},\label{inic}
\end{align}
for any $\mathbf{v}\in\mathbf{V}_n$ and for any $t\in[0,T]$ where ${\uu}_{n,0}$ is the $\mathbf{L}^2(\Omega)$-projection on $\mathbf{V}_n$ of \an{${\uu}_{0}-\mmm\in\HHH_0$}.

\an{Let us first notice that the quantity $\overline{\ww}_{n,\varepsilon}$ is necessary to be specified since any test function $\mathbf{v}\in\mathbf{V}_n$ has zero integral mean. Moreover, by construction,
$$
\overline{\uu}_{n,\varepsilon}\equiv \textbf{0},\quad \textbf{P}\uu_{n,\varepsilon}=\uu_{n,\varepsilon},\quad \textbf{P}\ww_{n,\varepsilon}=\ww_{n,\varepsilon}.
$$}

In the sequel we will denote by $C$ a generic positive constant independent of $n$. Any other dependence is explicitly pointed out if necessary.

Recalling that $\psi_\varepsilon'$ is at least $C^1(\R)$, we can locally solve the above Cauchy problem \eqref{ac}-\eqref{mum}, \eqref{inic} in the unknowns $\{{\hat{\alpha}}_i\}_i$ and find a unique maximal solution ${\hat{\alpha}}^{(n)}\in C^1([0, t_{n,\varepsilon}] ; \R^{ n})$, from which we also obtain by comparison a unique ${\delta}^{(n)}\in C^1([0, t_{n,\varepsilon}] ; \R^{ n})$. \an{Then by substitution in \eqref{subs} we immediately obtain the complete quantity $\ww_{n,\varepsilon}$.} \an{It is now standard to test \eqref{ac} by $\textbf{v}={\ww}_{n,\varepsilon}-\overline{{\ww}}_{n,\varepsilon}\in\mathbf{V}_n$ and obtain} the energy identity
\begin{align}
\ddt \mathcal{E}_{n,\varepsilon}+(\al(\ww_{n,\varepsilon}-\overline{\ww}_{n,\varepsilon}),\ww_{n,\varepsilon}-\overline{\ww}_{n,\varepsilon})=0,
\label{enid}
\end{align}
 where
 $$
 \mathcal{E}_{n,\varepsilon}:=\frac \gamma 2 \Vert \nabla {\uu}_{n,\varepsilon}\Vert^2+\int_\Omega \Psi_\varepsilon({\uu}_{n,\varepsilon}\an{+\mmm})dx.
 $$
 Let us observe that, being $\psi_\varepsilon'$ Lipschitz (see \eqref{approx}), and recalling that $\Psi_\varepsilon(\uu_0)\leq \Psi(\uu_0)$, we obtain
 \begin{align}
  \int_\Omega \Psi_\varepsilon({\uu}_{n,\varepsilon}(0)\an{+\mmm})dx&=\int_\Omega \left(\Psi_\varepsilon({\uu}_{n,0}\an{+\mmm})- \Psi_\varepsilon({\uu}_{0})\right)dx+\int_\Omega \Psi_\varepsilon(\uu_0)dx\nonumber\\&\leq C_\varepsilon \Vert {\uu}_{n,0}\an{+\mmm}-\uu_0\Vert+\int_\Omega \Psi(\uu_0)dx.
 \label{ii}
 \end{align}
 Therefore, since clearly $\Vert {\uu}_{n,0}\an{+\mmm}-\uu_0\Vert\to0$ as $n\to \infty$, for any $\varepsilon>0$ there exists $\overline{n}=\overline{n}(\varepsilon)$ such that
 \begin{align}
 \int_\Omega \Psi_\varepsilon({\uu}_{n,\varepsilon}(0)+\an{\mmm})dx\leq C,\quad \forall n>\overline{n}.
 \label{ini}
 \end{align}
 An application of Gronwall's Lemma then gives, thanks to \eqref{ini} and $\Vert \nabla {\uu}_{n,0}\Vert\leq \Vert \nabla \uu_0 \Vert$,
 \begin{align}
 \mathcal{E}_{n,\varepsilon}(t)+\int_0^t(\al(\ww_{n,\varepsilon}(s)-\overline{\ww}_{n,\varepsilon}(s)),\ww_{n,\varepsilon}(s)-\overline{\ww}_{n,\varepsilon}(s))ds\leq C,\quad \forall n>\overline{n}.
 \label{eni}
 \end{align}
 Now, recalling property (vi) of $\psi_\varepsilon$, it is immediate to see that, for any $\varepsilon<\varepsilon_0$,
 $$
 \int_\Omega \Psi_\varepsilon({\uu}_{n,\varepsilon}(t)\an{+\mmm})dx\geq -K,
 $$
 for some $K>0$, so that we can conclude, for any $\varepsilon<\varepsilon_0$,
 \begin{align}
 \Vert{\uu}_{n,\varepsilon}\Vert_{L^\infty(0,T;\mathbf{H}^1(\Omega))}+\Vert \Psi_\varepsilon({\uu}_{n,\varepsilon}\an{+\mmm})\Vert_{L^\infty(0,T;\textbf{L}^1(\Omega))}+\Vert \ww_{n,\varepsilon}-\overline{\ww}_{n,\varepsilon}\Vert_{L^2(0,T;\mathbf{L}^2(\Omega))}\leq C,\quad \forall n>\overline{n}(\varepsilon),
 \label{basis}
 \end{align}
 where we also exploited \eqref{pos}. \an{Clearly $C$ does not depend on $\varepsilon$.}  From this we can easily \an{ deduce that local maximal time $t_{n,\varepsilon}$ is $+\infty$.} Moreover, from these estimates it clearly derives, by comparison, that
 \begin{align}
 \Vert \partial_t {\uu}_{n,\varepsilon}\Vert_{L^2(0,T; \mathbf{L}^2(\Omega))}\leq C,\quad \forall n>\overline{n}.
 \label{dt}
 \end{align}
 These estimates, together with the fact that $\psi^\prime_\varepsilon$ is Lipschitz, give \an{from \eqref{subs}}
 $$
\Vert \overline{\ww}_{n,\varepsilon}\Vert_ {L^\infty(0,T)}\leq C_\varepsilon,\quad \forall n>\overline{n}.
 $$
 \an{Here $C_\varepsilon$ could depend on $\varepsilon$.}
 The obtained bounds are enough to pass to the limit as $n\to\infty$ by standard compactness arguments. However, since we also need to prove the existence of strong solutions, we now assume $\uu_0\in \mathbf{H}^2(\Omega)$ such that $\partial_\mathbf{n}\uu_0=\mathbf{0}$ almost everywhere on $\partial\Omega$, together with $\phi({u}_{0,i})\in L^2(\Omega)$, for any $i=1,\ldots,N$, and find a higher order estimate, before passing to the limit. In particular, we test \eqref{ac} with {$\mathbf{v}=\partial_t({\ww}_{n,\varepsilon}-\overline{
 {\ww}}_{n,\varepsilon})\in \mathbf{V}_n$.} Recalling that $\textbf{P}$ is selfadjoint and $\overline{\partial_t{\uu}_{n,\varepsilon}}\equiv \mathbf{0}$ by construction, we obtain
 \begin{align}
\frac 1 2 \ddt (\al(\ww_{n,\varepsilon}-\overline{\ww}_{n,\varepsilon}),\ww_{n,\varepsilon}-\overline{\ww}_{n,\varepsilon} )+ (\partial_t{\uu}_{n,\varepsilon},\partial_t{{\ww}}_{n,\varepsilon})=0.
\label{i}
 \end{align}
 Using \eqref{mum}, since $\partial_t{\uu}_{n,\varepsilon}\in \mathbf{V}_n$, we find
 \begin{align*}
 (\partial_t{\uu}_{n,\varepsilon},\partial_t{{\ww}}_{n,\varepsilon})=\sum_{i=1}^N\int_\Omega\phi_\varepsilon^\prime({u}_{n,\varepsilon,i}+\an{\mmm})\vert \partial_t{\uu}_{n,\varepsilon}\vert^2dx-(\partial_t{\uu}_{n,\varepsilon},\mathbf{A}\partial_t{\uu}_{n,\varepsilon}) +\Vert \nabla\partial_t {\uu}_{n,\varepsilon}\Vert^2.
 \end{align*}
 Being $\phi_\varepsilon^\prime\geq 0$ by property (iv) of $\psi_\varepsilon$, we have only to treat the term related to the matrix $\mathbf{A}$. This is readily done by comparison from \eqref{ac}: indeed, being $\mathbf{v}=\partial_{t}{\uu}_{n,\varepsilon}\in \mathbf{V}_n$, we get, by the Cauchy-Schwarz, Young, and Poincar\'{e} inequalities,
 \begin{align*}
 \vert(\partial_t{\uu}_{n,\varepsilon},\mathbf{A}\partial_t{\uu}_{n,\varepsilon})\vert\leq C\Vert\partial_t{\uu}_{n,\varepsilon}\Vert^2 \leq C\vert(\al(\ww_{n,\varepsilon}-\overline{\ww}_{n,\varepsilon}),\partial_t{\uu}_{n,\varepsilon})\vert\leq \frac 1 2 \Vert \nabla\partial_t {\uu}_{n,\varepsilon}\Vert^2+C\Vert \ww_{n,\varepsilon}-\overline{\ww}_{n,\varepsilon}\Vert^2.
 \end{align*}
Putting everything together in \eqref{i} and recalling \eqref{pos}, we end up with
\begin{align}
\frac 1 2 \ddt (\al(\ww_{n,\varepsilon}-\overline{\ww}_{n,\varepsilon}),\ww_{n,\varepsilon}-\overline{\ww}_{n,\varepsilon} )+\frac 1 2 \Vert \nabla\partial_t {\uu}_{n,\varepsilon}\Vert^2\leq C(\al(\ww_{n,\varepsilon}-\overline{\ww}_{n,\varepsilon}),\ww_{n,\varepsilon}-\overline{\ww}_{n,\varepsilon} ).
\label{end}
\end{align}
Observe now that, \an{from \eqref{mum},}
\begin{align*}
&\Vert (\al(\ww_{n,\varepsilon}(0)-\overline{\ww}_{n,\varepsilon}(0)),\ww_{n,\varepsilon}(0)-\overline{\ww}_{n,\varepsilon}(0) )\Vert\\&\leq C(\Vert \Delta {\uu}_{n,0}\Vert^2+\Vert \mathbf{A}{\uu}_{n,0}\Vert^2+\sum_{i=1}^N\Vert\phi_\varepsilon({u}_{n,0,i}\an{+\mmm})\Vert^2).
\end{align*}
On the other hand, by the properties of the eigenfunctions, we have
$$
\Vert \Delta {\uu}_{n,0}\Vert^2+\Vert \mathbf{A}{\uu}_{n,0}\Vert^2\leq C(\an{\vert \mmm\vert^2}+\Vert \Delta \uu_0\Vert^2+\Vert \uu_0\Vert^2)\leq C\Vert \uu_0\Vert^2_{\mathbf{H}^2(\Omega)}.
$$
Thus, recalling properties (ii)-(iii) of $\psi_\varepsilon$, we get
\begin{align*}
\sum_{i=1}^N\Vert\phi_\varepsilon({u}_{n,0,i}\an{+\mmm})\Vert^2&\leq 2\sum_{i=1}^N\Vert\phi_\varepsilon({u}_{n,0,i}\an{+\mmm})-\phi_\varepsilon({u}_{0,i})\Vert^2+2\sum_{i=1}^N\Vert\phi_\varepsilon({u}_{0,i})\Vert^2\\&\leq \frac{2}{\varepsilon^2}\Vert {\uu}_{n,0}\an{+\mmm}-{\uu}_{0}\Vert^2+2\sum_{i=1}^N\Vert\phi({u}_{0,i})\Vert^2.
\end{align*}
Therefore, since $\Vert {\uu}_{n,0}\an{+\mmm}-{\uu}_{0}\Vert\to0$ as $n\to\infty$ and by the stronger assumptions on the initial data, we deduce that for any $\varepsilon<\varepsilon_{0}$ there exists $\overline{\overline{n}}=\overline{\overline{n}}(\varepsilon)>0$ such that
\begin{align*}
\sum_{i=1}^N\Vert\phi_\varepsilon({u}_{n,0,i}\an{+\mmm})\Vert^2\leq C+2\sum_{i=1}^N\Vert\phi({u}_{0,i})\Vert^2, \quad \forall n>\overline{\overline{n}}.
\end{align*}
We can thus conclude that, for any $n>n_0(\varepsilon)=\max\{\overline{n},\overline{\overline{n}}\}$, owing to Gronwall's Lemma and \eqref{pos}, it holds
\begin{align}
\label{esss}
 \Vert \ww_{n,\varepsilon}-\overline{\ww}_{n,\varepsilon}\Vert_{L^\infty(0,T;\mathbf{L}^2(\Omega))}+\Vert \partial_t {\uu}_{n,\varepsilon}\Vert_{L^2(0,T;\mathbf{H}^1(\Omega))}\leq C(T),\quad \forall n>n_0,
\end{align}
\an{where $C(T)$ does not depend on $\varepsilon$.}
Furthermore, by comparison (choosing $\mathbf{v}=\partial_t{\uu}_{n,\varepsilon}$ in \eqref{ac}) it also holds
\begin{align}
\Vert \partial_t {\uu}_{n,\varepsilon}\Vert_{L^\infty(0,T;\mathbf{L}^2(\Omega))}\leq C(T),\quad \forall n>n_0.
\label{dtt}
\end{align}
We can now pass to the limit in $n$ for both the situations (according to the regularity of the initial data), to deduce, by standard compactness arguments,
that, for any $\varepsilon
<\varepsilon_0$, there exists a pair $(\mathbf{u}%
_\varepsilon,\ww_\varepsilon)$, defined on $[0,+\infty)$, with $\ww_\varepsilon(t)\in \widetilde{\VVV}_0$ for almost any $t\geq0$, such that (in the case of less regularity on $\uu_0$), for each $T>0$,
\begin{align*}
&\uu_\varepsilon\in L^\infty(0,T;\mathbf{H}^1(\Omega)), \\
& \partial_t\uu_\varepsilon\in L^2(0,T; \mathbf{L}^2(\Omega)),
\\
& \ww_\varepsilon\in L^2(0,T;\mathbf{L}^2(\Omega)).
\end{align*}
and
\begin{align}
\Vert \uu_\varepsilon\Vert_{ L^\infty(0,T;\mathbf{H}^1(\Omega))}+\Vert \partial_t\uu_\varepsilon\Vert_{L^2(0,T; \mathbf{L}^2(\Omega))}+\Vert\ww_\varepsilon-\overline{\ww}_\varepsilon\Vert_{ L^2(0,T;\mathbf{L}^2(\Omega))}\leq C(T),
\label{co}
\end{align}
for some $C(T)>0$ independent of $\varepsilon$, whereas there exists $C_\varepsilon>0$ such that
$$
\Vert\overline{\ww}_\varepsilon\Vert_{L^\infty(0,T)}\leq C_\varepsilon.
$$
If the stronger assumptions hold (see point (2) of Theorem \ref{thm}), then there exists a constant $C>0$, depending on the initial datum and on $T$,
	but independent of $\varepsilon $, such that
	\begin{equation}
	\Vert  \ww_{\varepsilon }-\overline{\ww}_{\varepsilon }\Vert_{L^\infty(0,T;\mathbf{L}^2(\Omega))) \cap L^2(0,T;\mathbf{H}^1(\Omega))) }+\Vert\partial _{t}\uu_\varepsilon\Vert _{L^2(0,T;\mathbf{H}^1(\Omega))}+\Vert\partial _{t}\uu_\varepsilon\Vert _{L^\infty(0,T;\mathbf{L}^2(\Omega))}\leq C,
	\label{est3b}
	\end{equation}%
    where the $L^2(0,T;\mathbf{H}^1(\Omega))$ control on the chemical potential differences is obtained by comparison in \eqref{e1} below.
It is then standard to show that $(\uu_\varepsilon,\ww_\varepsilon)$ satisfies
\begin{align}
\label{e1}
&\partial_t\uu_\varepsilon + \boldsymbol{\alpha}(\ww_\varepsilon - \overline{\ww}_\varepsilon)=0, \quad\text{ a.e. in } \Omega\times (0,T),\\
& (\ww_\varepsilon,\boldsymbol{\eta})=\gamma(\nabla\mathbf{u}%
_\varepsilon,\nabla\boldsymbol{\eta})+(\mathbf{P}(-\mathbf{A}\mathbf{u}%
_\varepsilon+\boldsymbol{\phi}_\varepsilon(\uu_\varepsilon)),%
\boldsymbol{\eta}), \quad\forall\, \boldsymbol{\eta}\in \mathbf{H}^1(\Omega), \quad\text{ a.e. in } (0,T), \label{e2}\\
&\uu_\varepsilon(0) = \uu_0, \quad\text{ a.e. in } \Omega. \label{e1in}
\end{align}
\an{Notice that, to be precise, we find that $\uu_{n,\varepsilon}$ converges in suitable norms to a function $\widetilde{\uu}_{\varepsilon}(t)\in \VVV_0$ (for almost any $t\geq0$) as $n\to \infty$. We then define $\uu_\varepsilon:=\widetilde{\uu}_\varepsilon+\mmm$ to obtain the results above.}
Then, by elliptic regularity, being $\boldsymbol{\phi}_\varepsilon$ Lipschitz, from \eqref{e2} we deduce its strong version, namely $\uu_\varepsilon \in L^2(0,T;\mathbf{H}^2(\Omega))$ and
\begin{align}
&\ww_\varepsilon=-\gamma \Delta \mathbf{u}_\varepsilon + \mathbf{P}(-\mathbf{A}\mathbf{u}_\varepsilon+\boldsymbol{\phi}_\varepsilon(\uu_\varepsilon)), \quad\text{ a.e. in } \Omega\times (0,T), \label{e2bis}\\
& \partial_\textbf{n}\uu_\varepsilon=0, \quad\text{ a.e. in } \partial\Omega\times (0,T). \label{e2bc}
\end{align}
By standard computations (see also \cite{EL} for similar results), we then have

\begin{itemize}
	\item Conservation of mass:
	\begin{equation*}
	\overline{\mathbf{u}}_\varepsilon(t)=\overline{\mathbf{u}}_0, \qquad \forall\, t\geq 0.
	\end{equation*}
	
	\item Conservation of total mass:
	\begin{equation}
	\sum_{i=1}^{N}u_{\varepsilon ,i}(x,t)=1,\qquad \text{ for a.a. } x\in \Omega \text{ and for all } t \in [0,\infty).
	   \label{mass}
	\end{equation}
	
	\item Conservation of chemical potential differences
	\begin{equation}
	\sum_{i=1}^{N}{w}_{\varepsilon ,i}=0,\qquad  \text{ for a.a. } x\in \Omega \text{ and a.a. } t \in (0,\infty).
	 \label{potentia}
	\end{equation}
	\item It holds from \eqref{enid}-\eqref{ii}, by standard arguments, the energy inequality
\begin{align}
 \mathcal{E}_{\varepsilon}(t)+\int_0^t(\al(\ww_{\varepsilon}-\overline{\ww}_{\varepsilon}),\ww_{\varepsilon}-\overline{\ww}_{\varepsilon})\leq  \mathcal{E}_{\varepsilon}(0),
\label{enid2}
\end{align}
for any $t\geq0$,
where
$$
\mathcal{E}_{\varepsilon}:=\frac \gamma 2 \Vert \nabla {\uu}_{\varepsilon}\Vert^2+\int_\Omega \Psi_\varepsilon({\uu}_{\varepsilon})dx.
$$
\end{itemize}
At this point, we can argue as in the proof of \cite[Thm 3.1]{GGPS} (which is based on \cite{Garke}), in order to control
$\overline{\ww}_{\varepsilon }(t)$ which then allows us to control $\Vert \ww_{\varepsilon }(t)\Vert$.
Following the proof of \cite[Lemma 3.3]{Garke}, we define
\begin{equation*}
\ww_{\varepsilon ,0}:=\ww_{\varepsilon }-\boldsymbol{\lambda }%
_{\varepsilon },
\end{equation*}%
where, on account of the boundary conditions,
\begin{equation*}
\boldsymbol{\lambda }_{\varepsilon }:=\overline{\ww}_{\varepsilon }=%
\overline{\mathbf{P}(-\mathbf{A}\uu_\varepsilon+\boldsymbol{\phi }%
	_{\varepsilon }(\uu_\varepsilon))}.
\end{equation*}%
Taking advantage of \eqref{e2}, we have,
\begin{equation}
(\ww_{\varepsilon ,0}+\boldsymbol{\lambda }_{\varepsilon },%
\boldsymbol{\eta })=\gamma (\nabla \uu_\varepsilon,\nabla
\boldsymbol{\eta })+(\mathbf{P}(-\mathbf{A}\uu_\varepsilon+%
\boldsymbol{\phi }_{\varepsilon }(\uu_\varepsilon)),\boldsymbol{%
	\eta }), \quad\forall\,\boldsymbol{\eta }\in \mathbf{H}^{1}(\Omega ),\; \text{ a.e. in } (0,T).
\label{e3}
\end{equation}%
Exploiting the convexity of $\Psi _{\varepsilon }^{1}$, for any $%
\mathbf{k}\in \mathbf{G}$, $\mathbf{G}$ being the Gibbs simplex, because $%
\mathbf{k}-\uu_\varepsilon\in T\Sigma $ almost everywhere in $\Omega\times (0,T)$, we find
\begin{align}
C& \geq \int_{\Omega }\Psi _{\varepsilon }^{1}(\mathbf{k})dx\geq \int_{\Omega
}\Psi _{\varepsilon }^{1}(\uu_\varepsilon)dx+\int_{\Omega }\Psi
_{\varepsilon ,\mathbf{u}}^{1}(\uu_\varepsilon)\cdot (\mathbf{k}-%
\uu_\varepsilon) dx \label{psi} \\
& =\int_{\Omega }\Psi _{\varepsilon }^{1}({\uu}_{\varepsilon
})dx+\int_{\Omega }\mathbf{P}\boldsymbol{\phi }_{\varepsilon }(\mathbf{u}%
_{\varepsilon })\cdot (\mathbf{k}-\uu_\varepsilon)dx,  \notag
\end{align}%
where we used (see property (i) of $\psi _{\varepsilon }$)
\begin{equation*}
\int_{\Omega }\Psi _{\varepsilon }^{1}(\mathbf{k}%
)dx\leq \int_{\Omega }\Psi ^{1}(\mathbf{k})dx\leq \max_{\mathbf{s}\in \lbrack
	0,1]}\vert \Psi^1 (\mathbf{s})\vert =C.
\end{equation*}
Here and in the sequel $C>0$ stands for a generic constant independent of $\varepsilon$.
Recalling that $\Psi
_{\varepsilon ,\mathbf{u}}^{1}(\uu_\varepsilon)=\{\phi
_{\varepsilon }(u_{\varepsilon ,i})\}_{i=1,\ldots ,N}$ and choosing $\boldsymbol{\eta }=\mathbf{k}-\mathbf{u}%
_{\varepsilon }$ in \eqref{e3}, on account of \eqref{psi}, we deduce that%
\begin{align*}
C& \geq \int_{\Omega }\Psi _{\varepsilon }^{1}(\mathbf{k})dx \\
& \geq \int_{\Omega }\Psi _{\varepsilon }^{1}(\uu_\varepsilon)dx+(%
\mathbf{P}(\mathbf{A}\uu_\varepsilon),\mathbf{k}-\mathbf{u}%
_{\varepsilon })dx \\
& +\gamma \Vert \nabla \uu_\varepsilon\Vert ^{2}+(\ww%
_{\varepsilon ,0},\mathbf{k}-\uu_\varepsilon)+(\boldsymbol{\lambda
}_{\varepsilon },\mathbf{k}-\uu_\varepsilon),
\end{align*}%
for almost all $t\in (0,T)$. On the other hand, we have ($\mathbf{k}\in \textbf{G}$ and thus $0\leq \mathbf{k}\leq 1$)
\begin{equation*}
\int_{\Omega }\sum_{i=1}^{N}k_{i}^{2}dx\leq \int_{\Omega }\left(
\sum_{i=1}^{N}k_{i}\right) ^{2}dx=|\Omega |_d.
\end{equation*}%
{\color{black}Then, using Cauchy-Schwarz's and Young's inequalities and recalling recalling property (vi) of $\psi_\varepsilon$ we obtain, }
\begin{align}
& (\boldsymbol{\lambda }_{\varepsilon },\mathbf{k}-{\uu}_{\varepsilon
})+\gamma \Vert \nabla \uu_\varepsilon\Vert ^{2}-K\leq (\boldsymbol{\lambda }_{\varepsilon },\mathbf{k}-{\uu}_{\varepsilon
})+\gamma \Vert \nabla \uu_\varepsilon\Vert ^{2} +\int_{\Omega }\Psi _{\varepsilon }^{1}(\uu_\varepsilon)dx \notag \\
& \leq C-(%
\mathbf{P}(\mathbf{A}\uu_\varepsilon),\mathbf{k}-\mathbf{u}%
_{\varepsilon })-(\ww_{\varepsilon ,0},\mathbf{k}-\mathbf{u}%
_{\varepsilon })  \notag \\
& \leq C\left( 1+\Vert \uu_\varepsilon\Vert +\Vert \uu%
_{\varepsilon }\Vert ^{2}+\Vert  \ww_{\varepsilon,0 }\Vert
(1+\Vert \uu_\varepsilon\Vert )\right) \leq C(1+\Vert
\ww_{\varepsilon,0 }\Vert ),  \label{e4}
\end{align}%
where in the last estimate we have exploited \eqref{co}. By the
conservation of mass and Remark \ref{control} we also deduce that, for all $%
i=1,\ldots ,N$ and all $t\in \lbrack 0,T]$,
\begin{equation*}
0<\delta _{0}<\overline{\mathbf{u}}_{\varepsilon ,i}(t)<1-(N-1)\delta
_{0}<1-\delta _{0}.
\end{equation*}%
Therefore, for any fixed $k,l=1,\ldots ,N$, we choose
\begin{equation*}
\mathbf{k}=\overline{\mathbf{u}}_{\varepsilon }+\delta _{0}\text{sign}(\lambda _{\varepsilon ,k}-\lambda _{\varepsilon ,l})(\pmb%
\zeta_{k}-%
\pmb%
\zeta_{l})\in \mathbf{G}
\end{equation*}%
in \eqref{e4}, where $\pmb\zeta_{j}:=(0,\ldots,\underbrace{1}_{j},\ldots,0)$. Thus, from %
\eqref{e4} we get that
\begin{equation}
|(\lambda _{\varepsilon ,k}-\lambda _{\varepsilon ,l})(t)|\leq \frac{C%
}{\delta _{0}|\Omega |_d}(1+\Vert \ww_{\varepsilon,0 }\Vert ).
\label{s}
\end{equation}%
Integrating $|(\lambda _{\varepsilon ,k}-\lambda _{\varepsilon
	,l})(t)|^{2}$ over $(0,T)$ and using the identity
\begin{equation*}
\boldsymbol{\lambda }_{\varepsilon }=\frac{1}{N}\left( \sum_{l=1}^{N}(\lambda _{\varepsilon ,k}-\lambda _{\varepsilon ,l})\right) _{k=1,\ldots
	,N},
\end{equation*}%
we find, owing to \eqref{co},
\begin{equation*}
\int_{0}^{T}|\boldsymbol{\lambda }_{\varepsilon }(t)|^{2}dt\leq C.
\end{equation*}
This, using again \eqref{co}, gives
\begin{equation}
\Vert \ww_{\varepsilon }\Vert _{L^{2}(0,T;\mathbf{L}^{2}(\Omega
	))}\leq C.  \label{ww1}
\end{equation}%
As a consequence, we deduce from \eqref{s} that
\begin{equation*}
|\boldsymbol{\lambda }_{\varepsilon }(t)|^{2}\leq C(1+\Vert \ww%
_{\varepsilon }(t)-\overline{\ww}%
_{\varepsilon }
(t)\Vert ^{2}),
\end{equation*}%
for almost any $t\in (0,T)$. Therefore, in the case of a more regular initial datum (see \eqref{est3b}), we have
\begin{equation}
\Vert \ww_{\varepsilon }\Vert _{L^{\infty }(0,T;\mathbf{L}%
	^{2}(\Omega ))}+\Vert \ww_{\varepsilon }\Vert _{L^{2}(0,T;\mathbf{H}^{1}(\Omega
	))}\leq C.  \label{mu2}
\end{equation}%
We are now left with some estimates related to ${\phi }_{\varepsilon }(u_{\varepsilon ,i})$. We follow again the proof of \cite[Thm.3.1]{GGPS}.
Being $\phi _{\varepsilon }^{\prime }$ bounded for a fixed $\varepsilon\in (0,\varepsilon _{0})$, we have that
\begin{equation*}
\nabla \phi _{\varepsilon }(u_{\varepsilon ,i})=\phi _{\varepsilon
}^{\prime }(u_{\varepsilon ,i})\nabla {u}_{\varepsilon,i}\in \mathbf{L}^{2}(\Omega ),
\end{equation*}
for almost any $t\in (0,T)$. Thus we can test \eqref{e2} with $\boldsymbol{\eta} =\boldsymbol{\phi }_{\varepsilon }(\uu_\varepsilon(t) )$ to get

\begin{align}
\sum_{i=1}^{N}(w_{\varepsilon ,i},{\phi }_{\varepsilon }(u_{\varepsilon ,i}))& =\sum_{i=1}^{N}\left( \gamma (\nabla u_{\varepsilon ,i},\phi _{\varepsilon }^{\prime }(u_{\varepsilon ,i})\nabla u_{\varepsilon ,i})\right)  \label{e5} \\
& +(\mathbf{P}(-\mathbf{A}\uu_\varepsilon+\boldsymbol{\phi }%
_{\varepsilon }(\uu_\varepsilon)),\boldsymbol{\phi }_{\varepsilon
}(\uu_\varepsilon)).  \notag
\end{align}%
Observe that
\begin{equation*}
(\mathbf{P}(\boldsymbol{\phi }_{\varepsilon }(\uu_\varepsilon)),%
\boldsymbol{\phi }_{\varepsilon }({u}_{\varepsilon
}))=\sum_{k=1}^{N}\int_{\Omega }\left( \phi _{\varepsilon }(u_{\varepsilon ,k})-\frac{1}{N}\sum_{l=1}^{N}\phi _{\varepsilon }(u_{\varepsilon ,l})\right) \phi _{\varepsilon }(u_{\varepsilon ,k})dx,
\end{equation*}%
and
\begin{align*}
& \sum_{k=1}^{N}\left( \phi _{\varepsilon }(u_{\varepsilon ,k})-%
\frac{1}{N}\sum_{l=1}^{N}\phi _{\varepsilon }(u_{\varepsilon ,l})\right) \phi _{\varepsilon }(u_{\varepsilon ,k}) \\
& =\frac{1}{N}\sum_{k,l=1}^{N}\left( \phi _{\varepsilon }(u_{\varepsilon ,k})-\phi _{\varepsilon }(u_{\varepsilon ,l})\right)
\phi _{\varepsilon }(u_{\varepsilon ,k}) \\
& =\frac{1}{N}\sum_{k<l}^{N}\left( \phi _{\varepsilon }(u_{\varepsilon ,k})-\phi _{\varepsilon }(u_{\varepsilon ,l})\right)
\phi _{\varepsilon }(u_{\varepsilon ,k})+\frac{1}{N}%
\sum_{k>l}^{N}\left( \phi _{\varepsilon }(u_{\varepsilon ,k})-\phi
_{\varepsilon }(u_{\varepsilon ,l})\right) \phi _{\varepsilon }(u_{\varepsilon ,k}) \\
& =\frac{1}{N}\sum_{k<l}^{N}\left( \phi _{\varepsilon }(u_{\varepsilon ,k})-\phi _{\varepsilon }(u_{\varepsilon ,l})\right)
(\phi _{\varepsilon }(u_{\varepsilon ,k})-\phi _{\varepsilon }(%
u_{\varepsilon ,l})) \\
& =\frac{1}{N}\sum_{k<l}^{N}\left( \phi _{\varepsilon }(u_{\varepsilon ,k})-\phi _{\varepsilon }(u_{\varepsilon ,l})\right)
^{2}.
\end{align*}%
Thanks to \eqref{mass}, we have
\begin{equation}
u_{\varepsilon ,m}:=\min_{i=1,\ldots ,N}u_{\varepsilon ,i}\leq \frac{1}{N}\leq \max_{i=1,\ldots ,N}u_{\varepsilon ,i}=:%
u_{\varepsilon ,M},  \label{minmax}
\end{equation}%
so that, being $\phi _{\varepsilon }$ monotone, we infer
\begin{align*}
& \frac{1}{N}\sum_{k<l}^{N}\left( \phi _{\varepsilon }(u_{\varepsilon ,k})-\phi _{\varepsilon }(u_{\varepsilon ,l})\right)
^{2}\geq \frac{1}{N}\left( \phi _{\varepsilon }(u_{\varepsilon,m})-\phi _{\varepsilon }(u_{\varepsilon ,M})\right) ^{2} \\
& \geq \frac{1}{N}\max_{i=1,\ldots ,N}\left( \phi _{\varepsilon }(u_{\varepsilon ,i})-\phi _{\varepsilon }\left( \frac{1}{N}\right) \right) ^{2}
\\
& \geq \frac{1}{N}\max_{i=1,\ldots ,N}\left( \frac{1}{2}\phi _{\varepsilon }(%
u_{\varepsilon ,i})^{2}-\phi _{\varepsilon }\left( \frac{1}{N}%
\right) ^{2}\right) \\
& \geq \frac{1}{2N}\max_{i=1,\ldots ,N}\phi _{\varepsilon }(u_{\varepsilon ,i})^{2}-C,
\end{align*}%
owing to the inequality $(a-b)^{2}\geq \tfrac{1}{2}a^{2}-b^{2}$. Notice
that $C$ is independent of $\varepsilon $ provided that we choose $%
\varepsilon $ sufficiently small. Indeed, since we have the pointwise
convergence $\phi _{\varepsilon }(\tfrac{1}{N})\rightarrow \phi (\tfrac{1}{N})$
as $\varepsilon \rightarrow 0$, then there exists $C>0$, independent of $\varepsilon$,
such that $|\phi _{\varepsilon }(\tfrac{1}{N})|\leq C$ for any $%
\varepsilon \in (0,\varepsilon _{0})$, with $\varepsilon_0>0$ sufficiently small. Then we get
$$
\sum_{i=1}^{N}(w_{\varepsilon ,i},{\phi }_{\varepsilon }(u_{\varepsilon ,i})) \leq \sum_{i=1}^{N}\Vert w_{\varepsilon,i}\Vert \Vert {\phi }_{\varepsilon }
(u_{\varepsilon ,i})\Vert
\leq C\Vert \ww_{\varepsilon }\Vert^2 +\frac{1}{8N}\int_{\Omega
}\max_{i=1,\ldots ,N}\phi _{\varepsilon }(u_{\varepsilon ,i})^{2}dx,
$$
and (see \eqref{co})
\begin{align*}
|(\mathbf{P}(-\mathbf{A}\uu_\varepsilon,\boldsymbol{\phi }%
_{\varepsilon }(\uu_\varepsilon))|& \leq C\Vert \mathbf{u}%
_{\varepsilon }\Vert ^{2}+\frac{1}{8N}\int_{\Omega }\max_{i=1,\ldots ,N}\phi
_{\varepsilon }(u_{\varepsilon ,i})^{2}dx \\
& \leq C+\frac{1}{8N}\int_{\Omega }\max_{i=1,\ldots ,N}\phi _{\varepsilon }(%
u_{\varepsilon ,i})^{2}dx.
\end{align*}%
Therefore, on account of the above inequalities and recalling that $\phi
_{\varepsilon }^{\prime }\geq 0$, we deduce from \eqref{e5} that
\begin{equation}
\frac{1}{4N}\int_{\Omega }\max_{i=1,\ldots ,N}\phi _{\varepsilon }(u_{\varepsilon ,i})^{2}dx\leq C\left( 1+\Vert \ww_{\varepsilon }\Vert
^{2}\right) ,  \label{l2}
\end{equation}%
which yields (see \eqref{ww1})
\begin{equation}
\Vert \boldsymbol{\phi }_{\varepsilon }(\uu_\varepsilon)\Vert
_{L^{2}(0,T;\boldsymbol{L}^{2}(\Omega ))}\leq C(T).  \label{phi}
\end{equation}
From this result, together with \eqref{co} and \eqref{ww1}, by elliptic
regularity, we infer from \eqref{e2bis}-\eqref{e2bc} that
\begin{equation}
\Vert \uu_\varepsilon\Vert _{L^{2}(0,T;\mathbf{H}^{2}(\Omega
	))}\leq C(T).  \label{H2}
\end{equation}%
Moreover, from \eqref{l2}, assuming a more regular initial datum, we infer (see \eqref{mu2})
\begin{equation}
\Vert \boldsymbol{\phi }_{\varepsilon }(\uu_\varepsilon)\Vert
_{L^{\infty}(0,T;\boldsymbol{L}^{2}(\Omega ))}\leq C(T),  \label{phi2}
\end{equation}
as well as
\begin{equation}
\Vert \uu_\varepsilon\Vert _{L^{\infty}(0,T;\mathbf{H}^{2}(\Omega
	))}\leq C(T).  \label{H2b}
\end{equation}%
We have obtained all the bounds we need to pass to the limit as $\varepsilon \rightarrow 0$. Being this step
standard (see, e.g., \cite{Garke}), we only present a sketch of the argument. By
compactness we immediately deduce that, up to subsequences,
\begin{align}
& \uu_\varepsilon\rightharpoonup \mathbf{u}\quad \text{ weakly* in }%
L^{\infty }(0,T;\mathbf{H}^{1}(\Omega )),  \notag \\
& \uu_\varepsilon\rightharpoonup \mathbf{u}\quad \text{ weakly in }%
L^{2}(0,T;\mathbf{H}^{2}(\Omega )),  \notag \\
&\partial _{t}\uu_\varepsilon\rightharpoonup%
\partial _{t}\mathbf{u}\quad \text{ weakly in }L^{2}(0,T;\mathbf{L}^{2}(\Omega )),  \notag \\
& \uu_\varepsilon\rightarrow \mathbf{u}\quad \text{ strongly in }%
L^{2}(0,T;\boldsymbol{L}^{2}(\Omega )),  \label{convv} \\
& \uu_\varepsilon\rightarrow \mathbf{u}\quad \text{ a.e. in }\Omega \times (0,T),  \notag \\
& \ww_{\varepsilon }\rightharpoonup \ww\quad \text{ weakly in }%
L^{2}(0,T;\mathbf{L}^{2}(\Omega )). \notag
\end{align}%
Then, arguing as in \cite[Section 6]{Garke} and exploiting \eqref{phi}, we infer that
\begin{align}
& \phi _{\varepsilon }(u_{\varepsilon ,k})\rightarrow \phi (u_k) \quad \text{ a.e. in }\Omega \times (0,T),  \label{ae} \\
& \phi _{\varepsilon }(u_{\varepsilon ,k})\rightharpoonup \phi (u_k) \quad \text{ weakly in }L^{2}(0,T;L^{2}(\Omega )),
\label{phi3}
\end{align}
for any $k=1,\ldots ,N$.
Thus the pair $(\mathbf{u},\ww)$ satisfies \eqref{eq1bis}-\eqref{eq1}. The energy inequality \eqref{ident0} is then retrieved by standard lower semicontinuity arguments. If the initial datum is more regular, then, up to subsequences,we also have the convergences
\begin{align}
& \uu_\varepsilon\rightharpoonup \mathbf{u}\quad \text{ weakly* in }%
L^{\infty }(0,T;\mathbf{H}^{2}(\Omega )),  \notag \\
& \uu_\varepsilon\rightharpoonup \mathbf{u}\quad \text{ weakly in }%
L^{2}(0,T;\mathbf{H}^{2}(\Omega )),  \notag \\
& \partial _{t}\uu_\varepsilon\rightharpoonup
\partial _{t}\mathbf{u}\quad \quad \text{ weakly* in } L^{\infty}(0,T;\mathbf{L}^{2}(\Omega )\ \text{and weakly in }L^{2}(0,T;\mathbf{H}^{1}(\Omega )),  \notag\\
& \ww_{\varepsilon }\rightharpoonup \ww\quad \text{ weakly* in } L^{\infty}(0,T;\mathbf{L}^{2}(\Omega )) \text{ and weakly in }L^{2}(0,T;\mathbf{H}^{1}(\Omega )),\notag\\&
\notag
\phi _{\varepsilon }(u_{\varepsilon ,k})\rightharpoonup \phi (u_k) \quad \text{ weakly* in }L^{\infty}(0,T;L^{2}(\Omega )),\quad \forall k=1,\ldots,N,
\end{align}%
which ensure the regularity of Theorem \ref{thm}, point (2).
The energy identity \eqref{ident0} can be recovered since $t \mapsto \Vert \nabla \uu(t)\Vert^2$ is absolutely continuous in  $[0,T]$ \an{and because of $\Psi^1(\uu)\in H^1(0,T;\LLL^1(\Omega))$ entailing that the function $t \mapsto \int_\Omega  \Psi(\uu(t))dx$ is absolutely continuous in  $[0,T]$ as well. Indeed, $\Vert \partial_t\Psi^1(\uu)\Vert_{\LLL^1(\Omega)}\leq\Vert \Psi^1_{,\pphi }(\uu)\Vert\Vert\partial_t\uu\Vert\leq C$ from the regularity above.}
This concludes the proof of the existence part of Theorem \ref{thm}.

\medskip\noindent
\textbf{Strict separation property of strong solutions.}
We recall that (\textbf{M1}) is in force. Let us now introduce the following notation: we
define $\pp^{s}$, with $s=1,\ldots,N-1$ and $\sigma\in \N$, \an{as} any possible subset of $s$ (non repeated) indices from $1,\ldots,N$.
\an{Note that $\sigma$ indicates the choice of the subset, and }$\sigma=1,\ldots,\binom{N}{s}.$ In case $s=N-1$ we define the only index not belonging to $\pp^{N-1}$ by $j_\sigma$.

\medskip\noindent
\textbf{Step 1. Case $N-1$.}
 Let us then start from $s=N-1$, having fixed $\sigma$. We consider the vector $\textbf{e}_\sigma^{N-1}$ as, for $i=1,\ldots,N$,
 $$
 (\textbf{e}_\sigma^{N-1})_i=\begin{cases}
 1 \text{ if }i\in \pp^{N-1},\\
 0\text{ if }i\not\in \pp^{N-1}.
 \end{cases}
 $$
Then  we take $\boldsymbol{\eta}=\eta\textbf{e}_\sigma^{N-1}$, for $\eta\in H^1(\Omega)$, in \eqref{phi1}. This gives
 \begin{align*}
 \left(\partial_t\left(\sum_{i\in\pp^{N-1}}{u}_i\right),\eta\right)+\left(\sum_{i\in\pp^{N-1}}\sum_{j=1}^N\al_{ij}(w_j-\overline{w}_j),\eta\right)=0,
 \end{align*}
 for almost any $t\in[0,T]$.
 We now fix $\delta>0$ (to be chosen later on) and consider $u_{\sigma,\delta}^{N-1}=\left(\sum_{i\in\pp^{N-1}}{u}_i-\delta\right)^-$.
 Setting $\eta=-u_{\sigma,\delta}^{N-1}$ and integrating by parts, we find
 \begin{align}
 &\dfrac 1 2 \dfrac {d}{dt}\left\Vert u_{\sigma,\delta}^{N-1}\right\Vert^2-\sum_{i\in\pp^{N-1}}\sum_{j=1}^N\int_\Omega \alpha_{ij}\psi^\prime({u}_j) u_{\sigma,\delta}^{N-1}dx-\gamma\sum_{i\in\pp^{N-1}}\sum_{j=1}^N\int_\Omega\alpha_{ij} \nabla {u}_j\cdot \nabla u_{\sigma,\delta}^{N-1}dx\nonumber\\&=-\sum_{i\in\pp^{N-1}}\sum_{j,k=1}^N\int_\Omega\alpha_{ij}\textbf{A}_{jk}{u}_ku_{\sigma,\delta}^{N-1}dx-\sum_{i\in\pp^{N-1}}\sum_{j=1}^N\int_\Omega \alpha_{ij}\overline{w}_ju_{\sigma,\delta}^{N-1}dx,
 \label{p}
 \end{align}
 where we used the property that, given any vector $\boldsymbol{\zeta}\in\R^N$, $\al\mathbf{P}\boldsymbol{\zeta}=\al\boldsymbol{\zeta}$.
 Now notice that, being $\alpha_{ii}=A>0$ for any $i=1,\ldots,N$, we have
 \begin{align*}
 &-\sum_{i\in\pp^{N-1}}\sum_{j=1}^N\int_\Omega\alpha_{ij} \nabla {u}_j\cdot \nabla u_{\sigma,\delta}^{N-1}dx\\&=-\sum_{i\in\pp^{N-1}}\int_\Omega\left(\alpha_{ii}\nabla {u}_i\cdot \nabla u_{\sigma,\delta}^{N-1}\right)dx-\sum_{i\in\pp^{N-1}}
\sum_{j\not=i,j=1}^N\int_\Omega\alpha_{ij} \nabla {u}_j\cdot \nabla u_{\sigma,\delta}^{N-1} dx
\\&=\int_\Omega\left(A\nabla u_{\sigma,\delta}^{N-1}\cdot \nabla u_{\sigma,\delta}^{N-1}\right)dx
-\sum_{i\in\pp^{N-1}}
\left(\sum_{j\not=i,j=1}^N\int_\Omega\alpha_{ij} \nabla {u}_j\cdot \nabla u_{\sigma,\delta}^{N-1} dx\right).
\end{align*}
Since $\alpha_{ij}=B<0$ for any $i\not=j$ (clearly we have $A+(N-1)B=0$), we see that the second summand becomes
\begin{align*}
&-\sum_{i\in\pp^{N-1}}
\left(\sum_{j\not=i,j=1}^N\int_\Omega\alpha_{ij} \nabla {u}_j\cdot \nabla u_{\sigma,\delta}^{N-1} dx\right)\\&=-B\sum_{i\in\pp^{N-1}}\sum_{j\not=i,j\in\pp^{N-1}}\int_\Omega\nabla {u}_j\cdot \nabla u_{\sigma,\delta}^{N-1} dx-B\sum_{i\in\pp^{N-1}}\sum_{j\not\in\pp^{N-1}}\int_\Omega \nabla {u}_j\cdot \nabla u_{\sigma,\delta}^{N-1} dx\\&
=-B(N-2)\sum_{j\in\pp^{N-1}}\int_\Omega \nabla {u}_j\cdot \nabla u_{\sigma,\delta}^{N-1} dx-B(N-1)\int_\Omega \nabla {u}_{j_\sigma}\cdot \nabla u_{\sigma,\delta}^{N-1} dx\\&
=B(N-2)\int_\Omega \nabla u_{\sigma,\delta}^{N-1}\cdot \nabla u_{\sigma,\delta}^{N-1}dx -B(N-1)\int_\Omega \nabla {u}_{j_\sigma}\cdot \nabla u_{\sigma,\delta}^{N-1} dx.
\end{align*}
 Notice also that, recalling ${u}_{j_\sigma}=1-\sum_{\an{j}\in \pp^{N-1}}{u}_j
$ and that $A=-B(N-1)$, it holds
$$
-B(N-1)\int_\Omega \nabla {u}_{j_\sigma}\cdot \nabla u_{\sigma,\delta}^{N-1} dx=-B(N-1)\int_\Omega \nabla u_{\sigma,\delta}^{N-1}\cdot \nabla u_{\sigma,\delta}^{N-1} dx=A\int_\Omega \nabla u_{\sigma,\delta}^{N-1}\cdot \nabla u_{\sigma,\delta}^{N-1} dx,
$$
\an{where we used the fact that, when $u_{\sigma,\delta}^{N-1}\leq \delta$, it holds
\begin{align}
\nonumber&\nabla u_{j_\sigma}=-\nabla\sum_{\an{j}\in \pp^{N-1}}{u}_j= -\nabla\left(\sum_{\an{j}\in \pp^{N-1}}{u}_j-\delta\right)\\&=-\nabla\left(\sum_{\an{j}\in \pp^{N-1}}{u}_j-\delta\right)^++\nabla\left(\sum_{\an{j}\in \pp^{N-1}}{u}_j-\delta\right)^-=\nabla\left(\sum_{\an{j}\in \pp^{N-1}}{u}_j-\delta\right)^-.
    \label{context}
\end{align}}
 Therefore, in the end we get
 \begin{align*}
 -\sum_{i\in\pp^{N-1}}\sum_{j=1}^N\int_\Omega\alpha_{ij} \nabla {u}_j\cdot \nabla u_{\sigma,\delta}^{N-1}dx=\left(2A+B(N-2)\right)\int_\Omega\vert\nabla u_{\sigma,\delta}^{N-1}\vert^2 dx\geq 0,
 \end{align*}
 recalling that $2A+B(N-2)=A-B\geq0$.

 Concerning the terms related to $\psi^\prime(u_j)$, we first observe that, on account of \eqref{general1}, we have
 $$
 \sum_{i\in \pp^{N-1}}\alpha_{ij}=-\alpha_{j_\sigma j}, \quad j=1,\ldots,N.
 $$
 Then we write
 \begin{align*}
& -\sum_{i\in\pp^{N-1}}\sum_{j=1}^N\int_\Omega \alpha_{ij}\psi^\prime({u}_j) u_{\sigma,\delta}^{N-1}dx\\&=-\sum_{i\in\pp^{N-1}}\sum_{j\in\pp^{N-1}}\int_\Omega \alpha_{ij}\psi^\prime({u}_j)u_{\sigma,\delta}^{N-1}dx-\sum_{i\in\pp^{N-1}}\sum_{j\not\in\pp^{N-1}}\int_\Omega \alpha_{ij}\psi^\prime({u}_j) u_{\sigma,\delta}^{N-1}dx\\&=\sum_{j\in\pp^{N-1}}\alpha_{j_\sigma j}\int_{E_{N-1}(t)} \psi^\prime({u}_j) u_{\sigma,\delta}^{N-1}dx+\alpha_{j_\sigma j_\sigma}\int_{E_{N-1}(t)} \psi^\prime({u}_{j_\sigma}) u_{\sigma,\delta}^{N-1}dx,
 \end{align*}
 where
 $$E_{N-1}(t)=\left\{x\in\Omega: \sum_{i\in\pp^{N-1}}{u}_i(x,t)\leq\delta\right\}.$$
 Observe that, in ${E_{N-1}(t)}$, it holds
 $$
1\geq \sum_{j\not\in \pp^{N-1}}{u}_{j}(t)={u}_{j_\sigma}(t)=1-\sum_{i\in\pp^{N-1}}{u}_i(t)>1-\delta.
 $$
 Thus, for $\delta\leq \tfrac 1 2$, we deduce
 \begin{align}
 \vert \psi^\prime(u_{j_\sigma})\vert\leq \max\left\{\left \vert \psi^\prime\left(\frac 1 2\right)\right\vert,\left \vert \psi^\prime\left(1\right)\right\vert\right\}\leq C,
 \label{easy}
 \end{align}
 being $\psi^\prime$ monotone increasing. Moreover, in $E_{N-1}(t)$ it also holds, being $\alpha_{j_\sigma j}=B=-\vert B\vert$ for $j\in\pp^{N-1}$ (recall that $B\leq0$),
 $$
 -\vert B\vert \psi^\prime({u}_j(t)) u_{\sigma,\delta}^{N-1}(t)\geq  -\vert B\vert \psi^\prime(\delta) u_{\sigma,\delta}^{N-1}(t),\quad\forall j\in\pp^{N-1},
 $$
 since we have $0\leq{u}_i(t)\leq \delta$ for any $i\in \pp^{N-1}$ and
 $$
 -\psi^\prime({u}_j)\geq -\psi^\prime(\delta),\quad \forall j\in\pp^{N-1}.
 $$
 Concerning the other terms in \eqref{p}, we have, clearly, being $0\leq{u}_k\leq 1$ for $k=1,\ldots,N$, that
 $$
 -\sum_{i\in\pp^{N-1}}\sum_{j,k=1}^N\int_\Omega\alpha_{ij}\textbf{A}_{jk}{u}_ku_{\sigma,\delta}^{N-1}dx\leq C\int_\Omega u_{\sigma,\delta}^{N-1}dx,
 $$
 and observing that (see \eqref{mu}) $\overline{\ww}\in L^\infty(0,T)$,
 we have, similarly,
 $$
 -\sum_{i\in\pp^{N-1}}\sum_{j=1}^N\int_\Omega \alpha_{ij}\overline{w}_ju_{\sigma,\delta}^{N-1}dx\leq C(T)\int_\Omega u_{\sigma,\delta}^{N-1}dx.
 $$
 Coming back to \eqref{p} and collecting all these results we end up with
 \begin{align*}
 &\dfrac 1 2 \dfrac {d}{dt}\left\Vert u_{\sigma,\delta}^{N-1}\right\Vert^2+\gamma\left(A-B\right)\int_\Omega\vert\nabla u_{\sigma,\delta}^{N-1}\vert^2 dx-(N-1)\vert B\vert \psi^\prime(\delta)\int_\Omega u_{\sigma,\delta}^{N-1}dx\\&\leq C(T)\int_\Omega u_{\sigma,\delta}^{N-1}dx-\alpha_{j_\sigma j_\sigma}\int_ {E_{N-1}(t)}\psi^\prime({u}_{j_\sigma}) u_{\sigma,\delta}^{N-1}dx\\&\leq C(T)\int_\Omega u_{\sigma,\delta}^{N-1}dx,
 \end{align*}
 so that, assuming $\delta$ sufficiently small to satisfy  (see assumption (\textbf{E1}))
 $$
 -(N-1)\vert B\vert \psi^\prime(\delta)-C(T)\geq0,
 $$
 we get, for almost any $t\in[0,T]$,
  \begin{align*}
 &\dfrac 1 2 \dfrac {d}{dt}\left\Vert u_{\sigma,\delta}^{N-1}\right\Vert^2\leq 0.
 \end{align*}
 Hence, having assumed the initial datum strictly separated, i.e., there exists $0<\delta_0\leq \tfrac 1 N$ such that
 \begin{align}
 {u}_{i,0}\geq \delta_0,\quad \forall i=1,\ldots,N,
 \label{initial}
 \end{align}
 we can choose $\delta\leq \delta_0$ in such a way that $u_{\sigma,\delta}^{N-1}(0)\equiv0$ and Gronwall's Lemma yields
 $$
\Vert u_{\sigma,\delta}^{N-1}(t)\Vert\equiv 0\quad \forall\, t\in[0,T].
 $$
 Notice now that the choice of the set $\pp^{N-1}$ is completely arbitrary, thus we infer that there exists $\delta_{N-1}$, such that $\delta_0\geq \delta_{N-1}>0$ and, for any possible $\pp^{N-1}$, with $\sigma=1,\ldots, N$,
 \begin{align}
 \sum_{i\in \pp^{N-1}}{u}_i(t)\geq \delta>0\quad \text{ in }\Omega,\quad\forall\, t\in[0,T],\quad \forall\, \delta\in (0,\delta_{N-1}).
 \label{sums1}
 \end{align}
 \begin{remark}
 	We point out that in the case $N=2$ the proof is ended. This means that $\textbf{(E2)}$ is not necessary in this case, consistently with \cite[Thm.3.5]{GGP}. \label{N2}
 	\end{remark}
\noindent
 \textbf{Step 2. Case $N-2$.}
 If $N=2$ we are done. Otherwise we need to consider the sets $\pp^{N-2}$, $\sigma=1,\ldots, \frac{N(N-1)}{2}$. Let us fix $\sigma$ and $0<\delta\leq \delta_{N-1}$ (to be chosen later on). Then, we set $$u_{\sigma,\delta}^{N-2}=\left(\sum_{i\in \pp^{N-2}}{u}_i-2\delta^2\right)^-$$
 and define the vector $\textbf{e}_\sigma^{N-2}$ as
 $$
 (\textbf{e}_\sigma^{N-2})_i=\begin{cases}
 1 \text{ if }i\in \pp^{N-2},\\
 0\text{ if }i\not\in \pp^{N-2}.
 \end{cases}
 $$
 for $i=1,\ldots,N$.

 We make a crucial observation: in the set  $$E_{N-2}(t)=\left\{x\in\Omega: \sum_{i\in\pp^{N-2}}{u}_i(x,t)\leq 2\delta^2\right\}$$
we infer from \eqref{sums1} that
 \begin{align}
 {u}_j(t)\geq \delta-2\delta^2,\quad \forall j\not\in \pp^{N-2}.
 \label{pp}
 \end{align}
 Recall that $\delta<\tfrac1 2$ and $0<2\delta^2< \delta\leq \delta_{N-1}<1$.
 Then  we take in \eqref{phi1}, as in Step 1, the test function $\boldsymbol{\eta}=\eta\textbf{e}_\sigma^{N-2}$, for $\eta\in H^1(\Omega)$ and we get
 \begin{align*}
 \left(\partial_t\left(\sum_{i\in\pp^{N-2}}{u}_i\right),\eta\right)+\left(\sum_{i\in\pp^{N-2}}\sum_{j=1}^N\alpha_{ij}({w}_j-\overline{w}_j),\eta\right)=0.
 \end{align*}
 Choosing in the equation above $\eta=-u_{\sigma,\delta}^{N-2}$ and integrating by parts, we find
 \begin{align}
 &\dfrac 1 2 \dfrac {d}{dt}\left\Vert u_{\sigma,\delta}^{N-2}\right\Vert^2-\sum_{i\in\pp^{N-2}}\sum_{j=1}^N\int_\Omega \alpha_{ij}\psi^\prime({u}_j) u_{\sigma,\delta}^{N-2}dx-\gamma\sum_{i\in\pp^{N-2}}\sum_{j=1}^N\int_\Omega\alpha_{ij} \nabla {u}_j\cdot \nabla u_{\sigma,\delta}^{N-2}dx\nonumber\\&=-\sum_{i\in\pp^{N-2}}\sum_{j,k=1}^N\int_\Omega\alpha_{ij}\textbf{A}_{jk}{u}_ku_{\sigma,\delta}^{N-2}dx-\sum_{i\in\pp^{N-2}}\sum_{j=1}^N\int_\Omega \alpha_{ij}\overline{w}_ju_{\sigma,\delta}^{N-2}dx.
 \label{p2}
 \end{align}
Recalling once more that $\alpha_{ii}=A>0$ for any $i=1,\ldots,N$, and arguing exactly as in Step 1, we find
 \begin{align*}
 &-\sum_{i\in\pp^{N-2}}\sum_{j=1}^N\int_\Omega\alpha_{ij} \nabla {u}_j\cdot \nabla u_{\sigma,\delta}^{N-2}dx\\&=-\sum_{i\in\pp^{N-2}}\int_\Omega\left(\alpha_{ii}\nabla {u}_i\cdot \nabla u_{\sigma,\delta}^{N-2}\right)dx-\sum_{i\in\pp^{N-2}}
 \sum_{j\not=i,j=1}^N\int_\Omega\alpha_{ij} \nabla {u}_j\cdot \nabla u_{\sigma,\delta}^{N-2} dx
 \\&=\int_\Omega\left(A\nabla u_{\sigma,\delta}^{N-2}\cdot \nabla u_{\sigma,\delta}^{N-2}\right)dx
 -\sum_{i\in\pp^{N-2}}
 \left(\sum_{j\not=i,j=1}^N\int_\Omega\alpha_{ij} \nabla {u}_j\cdot \nabla u_{\sigma,\delta}^{N-2} dx\right).
 \end{align*}
 Since $\alpha_{ij}=B<0$ for any $i\not=j$, the second summand becomes
 \begin{align*}
 &-\sum_{i\in\pp^{N-2}}
 \left(\sum_{j\not=i,j=1}^N\int_\Omega\alpha_{ij} \nabla {u}_j\cdot \nabla u_{\sigma,\delta}^{N-2} dx\right)\\&=-B\sum_{i\in\pp^{N-2}}\sum_{j\not=i,j\in\pp^{N-2}}\int_\Omega\nabla {u}_j\cdot \nabla u_{\sigma,\delta}^{N-2} dx-B\sum_{i\in\pp^{N-2}}\sum_{j\not\in\pp^{N-2}}\int_\Omega \nabla {u}_j\cdot \nabla u_{\sigma,\delta}^{N-2} dx\\&
 =-B(N-3)\sum_{j\in\pp^{N-2}}\int_\Omega \nabla {u}_j\cdot \nabla u_{\sigma,\delta}^{N-2} dx-B(N-2)\sum_{j\not\in\pp^{N-2}}\int_\Omega \nabla {u}_{j}\cdot \nabla u_{\sigma,\delta}^{N-2} dx\\&
 =B(N-3)\int_\Omega \nabla u_{\sigma,\delta}^{N-2}\cdot \nabla u_{\sigma,\delta}^{N-2}dx -B(N-2)\sum_{j\not\in\pp^{N-2}}\int_\Omega \nabla {u}_{j}\cdot \nabla u_{\sigma,\delta}^{N-2} dx.
 \end{align*}
Recall now that $\sum_{j\not\in\pp^{N-2}}{u}_{j}=1-\sum_{i\in\pp^{N-2}}{u}_i $ and $A=-B(N-1)$. Then, it holds
 \begin{align*}
 &-B(N-2)\sum_{j\not\in\pp^{N-2}}\int_\Omega \nabla {u}_{j}\cdot \nabla u_{\sigma,\delta}^{N-2} dx=-B(N-2)\int_\Omega \nabla u_{\sigma,\delta}^{N-2}\cdot \nabla u_{\sigma,\delta}^{N-2} dx\\&=(A+B)\int_\Omega \nabla u_{\sigma,\delta}^{N-2}\cdot \nabla u_{\sigma,\delta}^{N-2} dx.
 \end{align*}
This entails
 \begin{align*}
 -\sum_{i\in\pp^{N-2}}\sum_{j=1}^N\int_\Omega\alpha_{ij} \nabla {u}_j\cdot \nabla u_{\sigma,\delta}^{N-2}dx=\left(A+B(N-3)+A+B\right)\int_\Omega\vert\nabla u_{\sigma,\delta}^{N-2}\vert^2 dx\geq 0,
 \end{align*}
 since $A+B(N-3)+A+B=A-B\geq0$.

 The terms $\psi^\prime(u_j)$ can be handled as above. Indeed, observing that
 $$
 \sum_{i\in \pp^{N-2}}\alpha_{ij}=-\sum_{l\not\in \pp^{N-2}}\alpha_{l j}=-2B\geq 0, \quad \forall\,j\in \pp^{N-2},
 $$
 we obtain
 \begin{align*}
 & -\sum_{i\in\pp^{N-2}}\sum_{j=1}^N\int_\Omega \alpha_{ij}\psi^\prime({u}_j) u_{\sigma,\delta}^{N-2}dx\\&=-\sum_{i\in\pp^{N-2}}\sum_{j\in\pp^{N-2}}\int_\Omega \alpha_{ij}\psi^\prime({u}_j)u_{\sigma,\delta}^{N-2}dx-\sum_{i\in\pp^{N-2}}\sum_{j\not\in\pp^{N-2}}\int_\Omega \alpha_{ij}\psi^\prime({u}_j) u_{\sigma,\delta}^{N-2}dx\\&=\sum_{l\not\in \pp^{N-2}}\sum_{j\in\pp^{N-2}}\alpha_{l j}\int_{E_{N-2}(t)} \psi^\prime({u}_j) u_{\sigma,\delta}^{N-2}dx-B\sum_{i\in\pp^{N-2}}\sum_{j\not\in\pp^{N-2}}\int_{E_{N-2}(t)} \psi^\prime({u}_j) u_{\sigma,\delta}^{N-2}dx.
 \end{align*}
Thanks to \eqref{pp} we know that in $E_{N-2}(t)$, for $\delta$ sufficiently small, it holds
$$
\vert \psi^\prime({u}_j)\vert\leq -\psi^\prime(\delta-2\delta^2),\quad \forall j\not\in \pp^{N-2},
$$
being $\psi^\prime$ monotone increasing. This entails that
\begin{align*}
B\sum_{i\in\pp^{N-2}}\sum_{j\not\in\pp^{N-2}}\int_{E_{N-2}(t)} \psi^\prime({u}_j) u_{\sigma,\delta}^{N-2}dx&\leq -\vert B\vert\sum_{i\in\pp^{N-2}}\sum_{j\not\in\pp^{N-2}}\int_{E_{N-2}(t)} \psi^\prime(\delta-2\delta^2)u_{\sigma,\delta}^{N-2}dx\\&\leq -2\vert B\vert(N-2)\int_\Omega \psi^\prime(\delta-2\delta^2)u_{\sigma,\delta}^{N-2}dx.
\end{align*}
 Moreover, in $E_{N-2}(t)$ it also holds, being $0\geq B=-\vert B\vert$ and $u_{\sigma,\delta}^{N-2}\geq 0$,
 \begin{align*}
\sum_{l\not\in \pp^{N-2}}\sum_{j\in\pp^{N-2}}\alpha_{l j} \psi^\prime({u}_j)u_{\sigma,\delta}^{N-2}(t)&=-\vert B\vert \sum_{l\not\in \pp^{N-2}}\sum_{j\in\pp^{N-2}} \psi^\prime({u}_j)u_{\sigma,\delta}^{N-2}(t)\\&\geq  -2(N-2)\vert B\vert \psi^\prime(2\delta^2) u_{\sigma,\delta}^{N-2}(t),
 \end{align*}
 since $0\leq{u}_i(t)\leq 2\delta^2$ for any $i\in \pp^{N-2}$ and thus
 $$
 -\psi^\prime({u}_{j})\geq -\psi^\prime(2\delta^2),\quad \forall j\in\pp^{N-2}.
 $$
 Concerning the other terms in \eqref{p}, we have (recall that $0\leq{u}_k\leq 1$ for $k=1,\ldots,N$)
 $$
 -\sum_{i\in\pp^{N-2}}\sum_{j,k=1}^N\int_\Omega\alpha_{ij}\textbf{A}_{jk}{u}_ku_{\sigma,\delta}^{N-2}dx\leq C\int_\Omega u_{\sigma,\delta}^{N-2}dx,
 $$
 and, arguing similarly (see \eqref{mu}), we find
 $$
 -\sum_{i\in\pp^{N-2}}\sum_{j=1}^N\int_\Omega \alpha_{ij}\overline{w}_ju_{\sigma,\delta}^{N-2}dx\leq C(T)\int_\Omega u_{\sigma,\delta}^{N-2}dx.
 $$
 Combining \eqref{p2} with the obtained estimates, we end up with
 \begin{align*}
 &\dfrac 1 2 \dfrac {d}{dt}\left\Vert u_{\sigma,\delta}^{N-2}\right\Vert^2+\gamma\left(A-B\right)\int_\Omega\vert\nabla u_{\sigma,\delta}^{N-1}\vert^2 dx-2\vert B\vert (N-2)\psi^\prime(2\delta^2)\int_\Omega u_{\sigma,\delta}^{N-1}dx\\&\leq C(T)\int_\Omega u_{\sigma,\delta}^{N-1}dx-2\vert B\vert (N-2)\int_ {E_{N-2}(t)}\psi^\prime(\delta-2\delta^2) u_{\sigma,\delta}^{N-1}dx,
 \end{align*}
 that is,
  \begin{align*}
 &\dfrac 1 2 \dfrac {d}{dt}\left\Vert u_{\sigma,\delta}^{N-2}\right\Vert^2+\gamma\left(A-B\right)\int_\Omega\vert\nabla u_{\sigma,\delta}^{N-1}\vert^2 dx\\&+\left[2\vert B\vert (N-2)(-\psi^\prime(2\delta^2)+\psi^\prime(\delta-2\delta^2))-C(T)\right]\int_\Omega u_{\sigma,\delta}^{N-1}dx\leq 0.
 \end{align*}
Therefore, on account of (\textbf{E2}), for $0<\delta\leq \delta_{N-1}\leq \delta_0\leq \tfrac 1 N$, $\delta$ sufficiently small, we can ensure
 $$
 2\vert B\vert (N-2)(-\psi^\prime(\delta^2)+\psi^\prime(\delta-\delta^2))-C(T)\geq0.
 $$
Recalling that $A-B\geq0$, we deduce, for almost any $t\in[0,T]$,
 \begin{align*}
 &\dfrac 1 2 \dfrac {d}{dt}\left\Vert u_{\sigma,\delta}^{N-2}\right\Vert^2\leq 0.
 \end{align*}
 Then, thanks to \eqref{initial} and to the choice $\delta\leq \delta_{N-1}$ (entailing also $2\delta^2\leq \delta_{N-1}$), we get $ u_{\sigma,\delta}^{N-2}(0)\equiv 0$. Therefore, by Gronwall's Lemma, we get
 $$
 \Vert u_{\sigma,\delta}^{N-2}(t)\Vert \equiv 0\quad \forall t\in[0,T].
 $$
 Again the choice of the set $\pp^{N-2}$ is completely arbitrary, meaning that there exists a $0<\delta_{N-2}\leq \delta_{N-1}$ such that, for any possible $\pp^{N-2}$, with $\sigma=1,\ldots, \frac{N(N-1)}{2}$,
 \begin{align}
 \sum_{i\in \pp^{N-2}}{u}_i(t)\geq \delta>0\quad \text{ in }\Omega,\quad\forall t\in[0,T],\quad \forall 0<\delta\leq\delta_{N-2}.
 \label{sums2}
 \end{align}
 \subsection{Step 3. Iterative procedure and conclusion}
 If $N=3$ we are done. Otherwise we consider the sets $\pp^{N-3}$, $\sigma=1,\ldots, \binom{N}{N-3}$. Let us fix $\sigma$ and $\delta\leq \delta_{N-2}$ (to be chosen later on), introduce as before $$u_{\sigma,\delta}^{N-3}=\left(\sum_{i\in \pp^{N-3}}{u}_i-2\delta^2\right)^-$$
and define the vector $\textbf{e}_\sigma^{N-3}$ as, for $i=1,\ldots,N$,
$$
(\textbf{e}_\sigma^{N-3})_i=\begin{cases}
1 \text{ if }i\in \pp^{N-3},\\
0\text{ if }i\not\in \pp^{N-3}.
\end{cases}
$$
The essential observation is again the following: in the set  $$E_{N-3}(t)=\left\{x\in\Omega: \sum_{i\in\pp^{N-3}}{u}_i(x,t)\leq2\delta^2\right\}$$
from \eqref{sums2}, being $0<2\delta^2< \delta\leq \delta_{N-2}\leq \tfrac 1 N$, we deduce that
\begin{align}
{u}_j(t)\geq \delta-2\delta^2,\quad \forall j\not\in \pp^{N-3}.
\label{pp2}
\end{align}
This implies that in $E_{N-3}(t)$, for $\delta>0$ sufficiently small, it holds
$$
\vert \psi^\prime({u}_j(t))\vert\leq -\psi^\prime(\delta-2\delta^2)\quad\forall\, j\not \in \pp^{N-3},
$$
and
$$
-\psi^\prime({u}_{i}(t))\geq -\psi^\prime(2\delta^2)\quad \forall\, i\in \pp^{N-3}.
$$
We can now argue as in Step 2 and conclude that there exists a $\delta_{N-3}\in (0,\delta_{N-2}]$ such that, for any possible $\pp^{N-3}$, with $\sigma=1,\ldots, \binom{N}{N-3}$,
\begin{align}
\sum_{i\in \pp^{N-3}}{u}_i(t)\geq \delta>0\quad \text{ in }\Omega,\quad\forall\, t\in[0,T],\quad \forall\, \delta\in(0,\delta_{N-3}].
\label{sums3}
\end{align}
Repeating iteratively these arguments, we reach a generic step $m$ and we find $\delta_{N-m}\in (0,\delta_{N-m+1}]$ such that, for any $\pp^{N-m}$
with $\sigma=1,\ldots, \binom{N}{N-m}$, we have
\begin{align}
\sum_{i\in \pp^{N-m}}{u}_i(t)\geq \delta>0\quad \text{ in }\Omega,\quad\forall\, t\in[0,T],\quad \forall\, \delta\in (0,\delta_{N-m}].
\label{sumsm}
\end{align}

Therefore, we can continue the procedure until $N-m=1$, which entails in the end that
there exists a $0<\delta\leq \delta_0\leq \tfrac 1 N$ such that, for any $i=1,\ldots,N$,
\begin{align}
{u}_i(t)\geq \delta>0\quad \text{ in }\Omega,\quad\forall\, t\in[0,T],
\label{separ}
\end{align} 	
i.e., the strict separation property holds. This concludes the proof of Theorem \ref{thm}.

\subsection{Proof of Theorem \protect\ref{exp}}
Let us take $\boldsymbol{\eta } =\mathbf{u}(t)-\overline{\mathbf{u}}(t)$ in equation \eqref{mu1}. This gives
\begin{equation}
(\mathbf{P}\pmb\phi (\mathbf{u}),\mathbf{u}-\overline{\mathbf{u}})+\Vert
\nabla \mathbf{u}\Vert ^{2}=(\ww-\overline{\ww},\mathbf{u}-\overline{\mathbf{u}})+(%
\mathbf{Au},\mathbf{u}-\overline{\mathbf{u}}).
\label{la}
\end{equation}%
Moreover, by convexity of $\Psi ^{1}$ (recall that $\mathbf{u}-\overline{\mathbf{u}}\in T\Sigma $), we have
\begin{equation*}
(\mathbf{P}\pmb\phi (\mathbf{u}),\mathbf{u}-\overline{\mathbf{u}})=(\pmb\phi
(\mathbf{u}),\mathbf{u}-\overline{\mathbf{u}})\geq \int_{\Omega }\Psi ^{1}(%
\mathbf{u})dx-\int_{\Omega }\Psi ^{1}(\overline{\mathbf{u}})dx,
\end{equation*}%
but, being $\overline{\mathbf{u}}\equiv \overline{\mathbf{u}}_{0}$, it holds
\begin{equation*}
|\Psi ^{1}(\overline{\mathbf{u}})|\leq C,
\end{equation*}%
where $C>0$ depends on $\overline{\mathbf{u}}_{0}$. Applying standard inequalities, from \eqref{la} we infer that
$$
\int_{\Omega }\Psi ^{1}(\mathbf{u})dx+\Vert \nabla \mathbf{u}\Vert^2\leq C+C\Vert \nabla \mathbf{u}\Vert \Vert \ww-\overline{\ww}\Vert+(\mathbf{Au},\mathbf{u})-(\overline{\mathbf{u}},\mathbf{Au})
$$
and using \eqref{pos} we get
\begin{align*}
& \int_{\Omega }\Psi ^{1}(\mathbf{u})dx-\frac{1}{2}(\mathbf{Au},\mathbf{u})+%
\frac{1}{4}\Vert \nabla \mathbf{u}\Vert ^{2} \\
& \leq C(\pmb\alpha (\ww-\overline{\ww}),\ww-\overline{\ww})+\frac{1}{2}(\mathbf{%
	Au},\mathbf{u})-(\overline{\mathbf{u}},\mathbf{Au}) \\
& \leq C(1+(\pmb\alpha (\ww-\overline{\ww}),\ww-\overline{\ww}))+C\Vert \mathbf{u%
}\Vert ^{2} \\
& \leq C(1+(\pmb\alpha (\ww-\overline{\ww}),\ww-\overline{\ww}))+\frac{1}{2}%
\int_{\Omega }\Psi (\mathbf{u})dx,
\end{align*}%
where in the last step we applied property (vi) of the potential $\psi_\varepsilon$ (recall that these estimates must be obtained
in an approximating scheme, so for $\varepsilon $ sufficiently
small, see above). Therefore, we obtain
\begin{equation}
\frac{1}{4}\Vert \nabla \mathbf{u}\Vert ^{2}+\frac{1}{2}\int_{\Omega }\Psi (%
\mathbf{u})dx\leq C(1+(\pmb\alpha (\ww-\overline{\ww}),\ww-\overline{\ww})).
\label{inbis}
\end{equation}%
Combining \eqref{ident} with \eqref{inbis} (multiplied by $\epsilon >0$ sufficiently small), we end up with
\begin{equation*}
\frac{d}{dt}\mathcal{E}(t)+\frac{\epsilon }{2}\mathcal{E}(t)\leq \frac{d}{dt}%
\mathcal{E}(t)+\frac{\epsilon }{2}\mathcal{E}(t)+(1-\epsilon C)(\pmb\alpha (\ww-\overline{\ww}),\ww-\overline{\ww})\leq C,
\end{equation*}%
and the result follows from Gronwall's Lemma.
\subsection{Proof of Theorem \protect\ref{thm2}}

\label{pr2}

\textbf{Instantaneous regularization of weak solutions.} Thanks to Theorem \ref{thm}, point (1), for any $\tau>0$ we can find $\widetilde{\tau} \leq \tau$ such that $\uu(\widetilde{\tau})\in \mathbf{H}^2(\Omega)$ and $\partial_\textbf{n}\uu(\widetilde{\tau})=0$ on $\partial\Omega$ such that the solution starting from $\widetilde{\tau}$ is more regular. Having assumed (\textbf{M1}), this solution coincides with the weak one (generated from $\uu_0$) in $[\widetilde{\tau},+\infty)$ , by uniqueness, whence its instantaneous regularization and the validity of properties \eqref{regg0}-\eqref{reggtris}. Concerning the global bounds, for the sake of brevity, here we simply show the formal
estimates. A rigorous argument can be performed within an approximation scheme like the previous one. First, we observe that \eqref{ident} entails
\begin{equation}
\Vert \mathbf{u}\Vert _{L^{\infty }(0,\infty ;\mathbf{H}^{1}(\Omega
	))}+\Vert  \ww-\overline{\ww}\Vert _{L^{2}(t,t+1;\mathbf{L}^{2}(\Omega ))}\leq
C, \quad \forall\, t\geq 0.  \label{first}
\end{equation}%
Notice that the constant $C>0$ only depends on the initial energy $\mathcal{E}(0)$.
Then, arguing as in \eqref{end}, we obtain
\begin{equation*}
\frac{1}{2}\frac{d}{dt}(\al (\ww-\overline{\ww}),\ww-\overline{\ww})+\frac 1 2\Vert \nabla
\partial _{t}\mathbf{u}\Vert ^{2}\leq C(\al (\ww-\overline{\ww}),\ww-\overline{\ww}).
\end{equation*}%
Due to \eqref{first}, we can apply the uniform Gronwall's Lemma (see, e.g., \cite{Temam}, by choosing, e.g., $r=\tfrac{\tau }{2}$) to deduce, for any given $\tau >0$,
\begin{equation*}
\Vert  \ww-\overline{\ww}\Vert _{L^{\infty }(\tau ,\infty ;\mathbf{L}%
	^{2}(\Omega ))}+\Vert \nabla \partial _{t}\mathbf{u}\Vert _{L^{2}(t,t+1;%
	\mathbf{L}^{2}(\Omega ))}\leq C,\quad \forall\, t\geq \tau .
\end{equation*}%
From now on we can argue as in the proof of Theorem \ref{thm}, to get
\begin{equation*}
\Vert  \ww\Vert _{L^{\infty }(\tau ,\infty ;\mathbf{L}%
	^{2}(\Omega ))} \leq C,\quad \forall\, t\geq \tau .
\end{equation*}
where $C>0$, now and in the sequel, stands for a generic constant depending on $\Omega$, $\al$, $\Psi$, $\overline{\mathbf{u}}_{0}$, and $\mathcal{E}(0)$.
This allows us to deduce
\begin{equation}
\label{H2globbound}
\Vert \boldsymbol{\phi }_{ }(\uu )\Vert
_{L^{\infty}(\tau,\infty;\boldsymbol{L}^{2}(\Omega ))} + \Vert \uu\Vert _{L^{\infty}(\tau,\infty;\mathbf{H}^{2}(\Omega))}\leq C, \quad \forall\, t\geq \tau.
\end{equation}
Also, by comparison in \eqref{eq1}, we find
\begin{equation*}
\Vert  \ww\Vert _{L^{2}(t,t+1 ;\mathbf{H}^{1}(\Omega ))} \leq C,\quad \forall\, t\geq \tau.
\end{equation*}
The proof is finished.

\medskip\noindent
\textbf{Instantaneous strict separation.}
We are in the case $\Vert\uu_0\Vert_{\mathbf{L}^\infty(\Omega)}\leq 1$, that is, $\uu_0$ is not necessarily strictly separated like in Section \ref{exi}. Therefore we need to adapt the proof we performed in Section \ref{exi}. In order to do that, we perform a De Giorgi-type iterative scheme at each step.

The basic tool is the following
\begin{lemma}
	\label{conv}
	Let $\{y_n\}_{n\in\N\cup \{0\}}\subset \R^+$ satisfy the recursive inequalities
	\begin{align}
	y_{n+1}\leq Cb^ny_n^{1+\varepsilon},
	\label{ineq}\qquad \forall n\geq 0,
	\end{align}
	for some $C>0$, $b>1$ and $\varepsilon>0$. If
	\begin{align}
	\label{condition}
	y_0\leq \theta:= C^{-\frac{1}{\varepsilon}}b^{-\frac{1}{\varepsilon^2}},
	\end{align}
	then
	\begin{align}
	y_n\leq \theta b^{-\frac{n}{\varepsilon}},\qquad \forall n\geq 0,
	\label{yn}
	\end{align}
	and consequently $y_n\to 0$ for $n\to \infty$.
\end{lemma}
Lemma \eqref{conv} can be found, e.g., in \cite[Ch. I, Lemma 4.1]{DiBenedetto} (see also \cite[Ch.2, Lemma 5.6]{Lady}) and can be easily proven by induction \an{(see, e.g., \cite[Lemma 3.8]{P}).}
Since the iterative argument in which we sum up some components of $\uu$ (in decreasing number at each step) is exactly the same as in the case treated in Section \ref{exi}, we directly assume to be at Step $m>1$ and show the differences with respect to estimate \eqref{p} (Step 1 is even easier, as we have seen in Section \ref{exi}, thanks to the relation \eqref{easy}, thus it can be easily adapted following the analysis of the other steps). We assume to know, \an{for an arbitrary $\tau>0$}, that there exists $0<\delta_{N-m+1}\leq \tfrac 1 N$ such that, for any $\sigma$,
\begin{align}
\sum_{i\in \pp^{N-m+1}}{u}_i\geq \delta,\quad\text{ in }\overline{\Omega}\times\left[\frac \tau 2+\frac{m\tau}{2N},+\infty\right),\quad \forall \delta\leq \delta_{N-m+1},
\label{sumsmm}
\end{align}
with the same notation as in Section \ref{exi}. 
Notice that the upper bound $\delta_{N-m+1}\leq \tfrac 1 N$ is set since clearly in the end the necessary condition for the separation will be that $\delta\leq \tfrac 1 N$.
We now consider the set of indices $\pp^{N-m}$ for a certain $\sigma$. Then, for $i=1,\ldots,N$, we set
$$
(\textbf{e}_\sigma^{N-m})_i=\begin{cases}
1 \text{ if }i\in \pp^{N-m},\\
0\text{ if }i\not\in \pp^{N-m},
\end{cases}
$$
We can now perform De Giorgi's scheme.
Let us set $\delta$ sufficiently small such that $\delta\leq {\delta}_{N-m+1}$ and fix $\widetilde{\tau}$ such that
\begin{equation}
2\widetilde{\tau}+\frac{\tau}{2}+\frac{m\tau}{2N}= \frac \tau 2+\frac{(m+1)\tau}{2N},
\label{tt}
\end{equation}
i.e., $\widetilde{\tau}= \frac{\tau}{4N}$. Choose now $T>0$ such that $T-3\widetilde{\tau}= \frac{\tau}{2}+\frac{m\tau}{2N}\geq \frac \tau 2$, i.e., $T=\frac \tau 2+\dfrac{3+2m}{4N}\tau$.
Notice that condition \eqref{tt} implies
\begin{align}
T-\widetilde{\tau}=\frac \tau 2 +\frac{(m+1)\tau}{2N}.
\label{ttt}
\end{align}
Let us define the sequence
\begin{align}
k_n={\delta^2}+\frac{\delta^2}{2^{n}}, \quad \forall n\geq 0,
\label{kn}
\end{align}
where
\begin{align}
\delta^2< k_{n+1}<k_{n}<2\delta^2,\qquad \forall n\geq 1,\qquad k_n\to \delta^2\qquad \text{as }n\to \infty,
\label{kn1}
\end{align}
and the sequence of times
\begin{align}
\begin{cases}
t_{-1}=T-3\widetilde{\tau},\\
t_n=t_{n-1}+\frac{\widetilde{\tau}}{2^n},\qquad n\geq 0,
\end{cases}
\end{align}
which satisfies
$$
t_{-1}<t_n<t_{n+1}< T-\widetilde{\tau},\qquad \forall n\geq 0.
$$
Then, introduce a cutoff function $\eta_n\in C^1(\R)$ by setting
\begin{align}
\eta_n(t):=\begin{cases}
0,\quad t\leq t_{n-1},\\
1,\quad t\geq t_{n},
\end{cases}\text{ and }\quad \vert \eta^\prime_n(t)\vert\leq \frac{2^{n+1}}{\widetilde{\tau}},
\label{cutoff}
\end{align}
on account of the above definition of $\{t_n\}_n$, and set
\begin{align}
u_{\sigma,\delta}^{N-m,n}(x,t):=\left(\sum_{i\in \pp^{N-m}}{u}_i-k_n\right)^-.
\label{phik0}
\end{align}
Also, for any $n\geq 0$, let us introduce the interval $I_n=[t_{n-1},T]$ and the set
$$
A_n(t):=\left\{x\in \Omega: \sum_{i\in \pp^{N-m}}{u}_i(x,t)-k_n\leq 0\right\},\quad \forall t\in I_n,
$$
so that on $A_n(t)$ it holds (see \eqref{sumsmm}), being $0<2\delta^2< \delta\leq \delta_{N-
	m+1}\leq \tfrac 1 N$,
\begin{align}
{u}_j(t)\geq \delta-2\delta^2,\quad \forall j\not\in \pp^{N-m}.
\label{ppm}
\end{align}
This means that, on $A_{n}(t)$, for $\delta>0$ sufficiently small, we have
\begin{align}
\vert \psi^\prime({u}_j(t))\vert\leq -\psi^\prime(\delta-2\delta^2)\quad\forall j\not \in \pp^{N-m},
\label{1a}
\end{align}
and
\begin{align}
-\psi^\prime({u}_{i}(t))\geq -\psi^\prime(2\delta^2)\quad \forall i\in \pp^{N-m}.
\label{1b}
\end{align}
Observe now that
$$
I_{n+1}\subseteq I_n,\qquad \forall n\geq 0,$$
$$A_{n+1}(t)\subseteq A_n(t),\qquad \forall n\geq 0,\qquad \forall t\in I_{n+1},
$$
and set
$$
y_n=\int_{I_n}\int_{A_n(s)}1dxds,\qquad \forall n\geq0.
$$
For any $n\geq 0$, we take the test function $\boldsymbol{\eta }=\textbf{e}_\sigma^{N-m}u_{\sigma,\delta}^{N-m,n}\eta_n^2$ in \eqref{phi1}, and integrate over $[t_{n-1},t]$, $t_n\leq t\leq T$. After an integration by parts, we get
\begin{align}
&\dfrac 1 2 \eta_n^2(t)\left\Vert u_{\sigma,\delta}^{N-m,n}(t)\right\Vert^2\underbrace{-\sum_{i\in\pp^{N-m}}\sum_{j=1}^N\int_{t_{n-1}}^t\int_\Omega \alpha_{ij}\psi^\prime({u}_j) u_{\sigma,\delta}^{N-m,n}\eta_n^2dxds}_{I_1}\nonumber\\&\underbrace{-\gamma\sum_{i\in\pp^{N-m}}\sum_{j=1}^N\int_{t_{n-1}}^t\int_\Omega\alpha_{ij} \nabla {u}_j\cdot \nabla u_{\sigma,\delta}^{N-m,n}\eta_n^2dxds\nonumber}_{I_2}=\underbrace{-\sum_{i\in\pp^{N-m}}\sum_{j,k=1}^N\int_{t_{n-1}}^t\int_\Omega\alpha_{ij}\textbf{A}_{jk}{u}_ku_{\sigma,\delta}^{N-m,n}\eta_n^2dxds}_{I_3}\\&\underbrace{-\sum_{i\in\pp^{N-m}}\sum_{j=1}^N\int_{t_{n-1}}^t\int_\Omega \alpha_{ij}\overline{w}_ju_{\sigma,\delta}^{N-m,n}\eta_n^2dxds}_{I_4}+\underbrace{\int_{t_{n-1}}^t\partial_t\eta_n\eta_n\Vert u_{\sigma,\delta}^{N-m,n}\Vert^2ds}_{I_5},
\label{p1}
\end{align}
where we used
$$
\dfrac 1 2 \eta_n^2(t)\left\Vert u_{\sigma,\delta}^{N-m,n}(t)\right\Vert^2=\int_{t_{n-1}}^t<\partial_tu_{\sigma,\delta}^{N-m,n},u_{\sigma,\delta}^{N-m,n}\eta_n^2>ds-
\int_{t_{n-1}}^t\partial_t\eta_n\eta_n\Vert u_{\sigma,\delta}^{N-m,n}\Vert^2ds.$$
As in Section \ref{exi}, recalling that $\alpha_{ii}=A>0$ for any $i=1,\ldots,N$, we obtain
\begin{align*}
&I_2=-\gamma\sum_{i\in\pp^{N-m}}\sum_{j=1}^N\int_{t_{n-1}}^t\int_\Omega\alpha_{ij} \nabla {u}_j\cdot \nabla u_{\sigma,\delta}^{N-m,n}\eta_n^2dxds\\&=-\gamma\sum_{i\in\pp^{N-m}}\int_{t_{n-1}}^t\int_\Omega\left(\alpha_{ii}\nabla {u}_i\cdot \nabla u_{\sigma,\delta}^{N-m,n}\right)\eta_n^2dxds\\&-\gamma\sum_{i\in\pp^{N-m}}
\sum_{j\not=i,j=1}^N\int_{t_{n-1}}^t\int_\Omega\alpha_{ij} \nabla {u}_j\cdot \nabla u_{\sigma,\delta}^{N-m,n}\eta_n^2 dxds
\\&=\gamma\int_{t_{n-1}}^t\int_\Omega\left(A\nabla u_{\sigma,\delta}^{N-m,n}\cdot \nabla u_{\sigma,\delta}^{N-m,n}\right)\eta_n^2dxds
\\&-\gamma\sum_{i\in\pp^{N-m}}
\left(\sum_{j\not=i,j=1}^N\int_{t_{n-1}}^t\int_\Omega\alpha_{ij} \nabla {u}_j\cdot \nabla u_{\sigma,\delta}^{N-m,n}\eta_n^2 dxds\right).
\end{align*}
Since $\alpha_{ij}=B<0$ for any $i\not=j$, the second summand becomes
\begin{align*}
&-\sum_{i\in\pp^{N-m}}
\left(\sum_{j\not=i,j=1}^N\int_\Omega\alpha_{ij} \nabla {u}_j\cdot \nabla u_{\sigma,\delta}^{N-m,n} dx\right)\\&=-B\sum_{i\in\pp^{N-m}}\sum_{j\not=i,j\in\pp^{N-m}}\int_\Omega\nabla {u}_j\cdot \nabla u_{\sigma,\delta}^{N-m,n} dx-B\sum_{i\in\pp^{N-m}}\sum_{j\not\in\pp^{N-2}}\int_\Omega \nabla {u}_j\cdot \nabla u_{\sigma,\delta}^{N-m} dx\\&
=B(N-m-1)\sum_{j\in\pp^{N-m}}\int_\Omega \nabla {u}_j\cdot \nabla u_{\sigma,\delta}^{N-m,n} dx-B(N-m)\sum_{j\not\in\pp^{N-m}}\int_\Omega \nabla {u}_{j}\cdot \nabla u_{\sigma,\delta}^{N-m,n} dx.
\end{align*}
On the other hand, observe that $\sum_{j\not\in\pp^{N-m}}{u}_{j}=1-\sum_{i\in\pp^{N-m}}{u}_{i}$. Hence, we get
$$
-B(N-m)\sum_{j\not\in\pp^{N-m}}\int_\Omega \nabla {u}_{j}\cdot \nabla u_{\sigma,\delta}^{N-m,n} dx=-B(N-m)\int_\Omega \nabla u_{\sigma,\delta}^{N-m,n}\cdot \nabla u_{\sigma,\delta}^{N-m,n} dx.
$$
This entails (recall that $A=-B(N-1)$)
\begin{align*}
I_2&=\gamma\left(A+B(N-m-1)-B(N-m)\right)\int_{t_{n-1}}^t\int_\Omega\vert\nabla u_{\sigma,\delta}^{N-2}\vert^2\eta_n^2 dxds\\&=\gamma N\vert B\vert\int_{t_{n-1}}^t\int_\Omega\vert\nabla u_{\sigma,\delta}^{N-m,n}\vert^2\eta_n^2 dxds,
\end{align*}
since $A+B(N-m-1)-B(N-m)=A-B=N\vert B\vert>0$.\\
Concerning $I_1$, recalling that, for any $j\in \pp^{N-m}$, $\sum_{i\in \pp^{N-m}}\alpha_{ij}=-\sum_{l\not\in \pp^{N-m}}\alpha_{l j}$,
we can write
\begin{align*}
& -\sum_{i\in\pp^{N-m}}\sum_{j=1}^N\int_\Omega \alpha_{ij}\psi^\prime({u}_j) u_{\sigma,\delta}^{N-m,n}dx\\&=-\sum_{i\in\pp^{N-m}}\sum_{j\in\pp^{N-m}}\int_\Omega \alpha_{ij}\psi^\prime({u}_j)u_{\sigma,\delta}^{N-m,n}dx-\sum_{i\in\pp^{N-m}}\sum_{j\not\in\pp^{N-m}}\int_\Omega \alpha_{ij}\psi^\prime({u}_j) u_{\sigma,\delta}^{N-m,n}dx\\&=\sum_{l\not\in \pp^{N-m}}\sum_{j\in\pp^{N-m}}\alpha_{l j}\int_{A_n(t)} \psi^\prime({u}_j) u_{\sigma,\delta}^{N-m,n}dx-B\sum_{i\in\pp^{N-m}}\sum_{j\not\in\pp^{N-m}}\int_{A_n(t)} \psi^\prime({u}_j) u_{\sigma,\delta}^{N-m,n}dx.
\end{align*}
Thus, by \eqref{1a}, we deduce
\begin{align*}
B\sum_{i\in\pp^{N-m}}\sum_{j\not\in\pp^{N-m}}\int_{A_n(t)}\psi^\prime({u}_j) u_{\sigma,\delta}^{N-m,n}dx&\leq -\vert B\vert\sum_{i\in\pp^{N-m}}\sum_{j\not\in\pp^{N-m}}\int_{A_n(t)} \psi^\prime(\delta-2\delta^2)u_{\sigma,\delta}^{N-m,n}dx\\&\leq -m\vert B\vert(N-m)\int_\Omega \psi^\prime(\delta-2\delta^2)u_{\sigma,\delta}^{N-m,n}dx.
\end{align*}
Moreover, being $0\geq B=-\vert B\vert$ and thanks to \eqref{1b}, in $A_{n}(t)$ it also holds
\begin{align*}
\sum_{l\not\in \pp^{N-m}}\sum_{j\in\pp^{N-m}}\alpha_{l j} \psi^\prime({u}_j)u_{\sigma,\delta}^{N-m,n}&=-\vert B\vert \sum_{l\not\in \pp^{N-m}}\sum_{j\in\pp^{N-m}} \psi^\prime({u}_j)u_{\sigma,\delta}^{N-m,n}\\&\geq  -m(N-m)\vert B\vert \psi^\prime(2\delta^2)u_{\sigma,\delta}^{N-m,n}.
\end{align*}
About the other terms in \eqref{p}, recalling that $0\leq{u}_k\leq 1$ for $k=1,\ldots,N$, we clearly have
$$
I_3=-\sum_{i\in\pp^{N-m}}\sum_{j,k=1}^N\int_{t_{n-1}}^t\int_\Omega\alpha_{ij}\textbf{A}_{jk}{u}_ku_{\sigma,\delta}^{N-m,n}\eta_n^2dxds\leq C\int_{t_{n-1}}^t\int_\Omega u_{\sigma,\delta}^{N-m,n}\eta_n^2dxds,
$$
and by \eqref{reggbis} on $(\tfrac \tau 2,+\infty)$
we have, similarly,
$$
I_4=-\sum_{i\in\pp^{N-m}}\sum_{j=1}^N\int_{t_{n-1}}^t\int_\Omega \alpha_{ij}\overline{w}_ju_{\sigma,\delta}^{N-m,n}\eta_n^2dxds\leq C\int_{t_{n-1}}^t\int_\Omega u_{\sigma,\delta}^{N-m,n}\eta_n^2dxds.
$$
We are left with $I_5$, which is not present when $\uu_0$ is strictly separated.
Note that, since $\sum_{i=1}^N{u}_i=1$ and  $0\leq {u}_i\leq 1$ almost everywhere in $\Omega\times[0,+\infty)$, for any $i=1,\ldots,N$, we have
$$
0\leq \sum_{i\in \pp^{N-m}}{u}_i,\quad\text{ a.e. in }\Omega\times[0,+\infty),
$$
and thus
\begin{align}
\label{delti}
0\leq u_{\sigma,\delta}^{N-m,n}\leq 2\delta^2\quad\text{ a.e. in }\Omega\times[0,+\infty).
\end{align}
Then, thanks to the above inequality, we infer
\begin{align}
\nonumber
&I_5=\int_{t_{n-1}}^t\Vert u_{\sigma,\delta}^{N-m,n}(s)\Vert^2\eta_n\partial_t\eta_nds=\int_{t_{n-1}}^t\int_\Omega  \left(u_{\sigma,\delta}^{N-m,n}(s)\right)^2\eta_n\partial_t\eta_ndxds\\&=\int_{t_{n-1}}^t\int_{A_n(s)}  \left(u_{\sigma,\delta}^{N-m,n}(s)\right)^2\eta_n\partial_t\eta_ndxds
\leq
\int_{t_{n-1}}^t\int_{A_n(s)} 4\delta^4\frac{2^{n+1}}{\widetilde{\tau}}dxds\leq \frac{2^{n+3}\delta^4}{\widetilde{\tau}}y_n.
\label{del}
\end{align}
Therefore, collecting all the above results, we end up with
\begin{align}
&\dfrac 1 2 \eta_n^2(t)\left\Vert u_{\sigma,\delta}^{N-m,n}(t)\right\Vert^2-m(N-m)\vert B\vert\psi^\prime(2\delta^2){\int_{t_{n-1}}^t\int_\Omega u_{\sigma,\delta}^{N-m,n}\eta_n^2dxds}\nonumber+{\gamma m\vert B\vert \int_{t_{n-1}}^t\int_\Omega\vert\nabla u_{\sigma,\delta}^{N-m,n}\vert^2\eta_n^2dxds\nonumber}\\&\leq (C-m\vert B\vert(N-m) \psi^\prime(\delta-2\delta^2))\int_{t_{n-1}}^t\int_\Omega u_{\sigma,\delta}^{N-m,n}\eta_n^2dxds+\frac{2^{n+3}\delta^4}{\widetilde{\tau}}y_n,
\label{p1f}
\end{align}
for any $t\in[t_{n},T]$, i.e.,
\begin{align}
&\nonumber\dfrac 1 2 \eta_n^2(t)\left\Vert u_{\sigma,\delta}^{N-m,n}(t)\right\Vert^2+\left[m(N-m)\vert B\vert(-\psi^\prime(2\delta^2)+\psi^\prime(\delta-2\delta^2))-C\right]{\int_{t_{n-1}}^t\int_\Omega u_{\sigma,\delta}^{N-m,n}\eta_n^2dxds}\\&+{\gamma m\vert B\vert \int_{t_{n-1}}^t\int_\Omega\vert\nabla u_{\sigma,\delta}^{N-m,n}\vert^2\eta_n^2dxds}\leq \frac{2^{n+3}\delta^4}{\widetilde{\tau}}y_n,
\label{p1f2}
\end{align}
We now recall assumption (\textbf{E2}), so that
\begin{align}
\psi^\prime(\delta-2\delta^2)-\psi^\prime(2\delta^2)\to +\infty\quad\text{ as }\delta\to 0^+.
\label{assumF0}
\end{align}
Therefore, for $0<\delta\leq \delta_{N-m+1}$ sufficiently small, we have
$$
m\vert B\vert (N-m)(-\psi^\prime(2\delta^2)+\psi^\prime(\delta-2\delta^2))-C\geq0.
$$
This entails
\begin{align}
\max_{t\in I_{n+1}}\Vert u_{\sigma,\delta}^{N-m,n} (t)\Vert^2\leq X_n,\qquad  {\color{black}}2\gamma m\vert B\vert \int_{I_{n+1}}\Vert \nabla u_{\sigma,\delta}^{N-m,n} \Vert^2ds \leq X_n,
\label{est}
\end{align}
where
$$
X_n:= \frac{2^{n+4}\delta^4}{\widetilde{\tau}}y_n.
$$
On the other hand, for any $t\in I_{n+1}$ and for almost any $x\in A_{n+1}(t)$, we get
\begin{align*}
& u_{\sigma,\delta}^{N-m,n} (x,t)=\delta^2+\frac{\delta^2}{2^{n}}-\sum_{i\in \pp^{N-m}}{u}_i(x,t)\\&=
\underbrace{-\sum_{i\in \pp^{N-m}}{u}_i(x,t)+\left[\delta^2+\frac{\delta^2}{2^{n+1}}\right]}_{ u_{\sigma,\delta}^{N-m,n+1}(x,t)\geq 0}+\delta^2\left[\frac{1}{2^{n}}-\frac{1}{2^{n+1}}\right]\geq \frac{\delta^2}{2^{n+1}},
\end{align*}
which implies
\begin{align*}
\int_{I_{n+1}}\int_{\Omega}\left\vert u_{\sigma,\delta}^{N-m,n} \right\vert^3dxds&\geq \int_{I_{n+1}}\int_{A_{n+1}(s)}\left\vert u_{\sigma,\delta}^{N-m,n} \right\vert^3dxds\\&\geq \left(\frac{\delta^2}{2^{n+1}}\right)^3\int_{I_{n+1}}\int_{A_{n+1}(s)}dxds\\&=\left(\frac{\delta^2}{2^{n+1}}\right)^3y_{n+1}.
\end{align*}
Then, for $d=2,3$, we find
\begin{align}
&\nonumber\left(\frac{\delta^2}{2^{n+1}}\right)^3y_{n+1}\leq \int_{I_{n+1}}\int_{\Omega}\vert u_{\sigma,\delta}^{N-m,n} \vert^3dxds\\&=\int_{I_{n+1}}\int_{A_n(s)}\vert u_{\sigma,\delta}^{N-m,n} \vert^3dxds \leq \left(\int_{I_{n+1}}\int_{\Omega}\vert u_{\sigma,\delta}^{N-m,n} \vert^{\frac{2d+ 4} {d}}dxds\right)^{\frac{3d}{2d+4}}\left(\int_{I_{n+1}}\int_{A_n(s)}dxds\right)^{\frac{4-d}{2d+4}}.
\label{est2}
\end{align}
Thanks to Sobolev-Gagliardo-Nirenberg's inequality and Poincar\'{e}'s inequality, we have
\begin{align}
\label{sobolev}
&\| u\|_{L^\frac{2d+4}{d}(\Omega)}\leq \tilde{C} \left(\|u-\overline{u}\|^{\frac{2}{d+2}}\|\nabla u\|^{\frac{d}{d+2}}+\vert \overline{u}\vert\right), \quad &&\forall \, u \in H^1(\Omega),
\end{align}
so we get
\begin{align*}
&\int_{I_{n+1}}\int_{\Omega}\vert u_{\sigma,\delta}^{N-m,n} \vert^{\frac{2d+4}{d}}dxds\leq \tilde{C} \int_{I_{n+1}}\left(\Vert u_{\sigma,\delta}^{N-m,n} \Vert^{\frac{4d+8}{d(d+2)}}\Vert \nabla u_{\sigma,\delta}^{N-m,n} \Vert^2+\vert \overline{ u_{\sigma,\delta}^{N-m,n}}\vert^{\frac{2d+4}d}\right)ds\\&\leq \hat{C}\int_{I_{n+1}}\left(\Vert u_{\sigma,\delta}^{N-m,n} \Vert^{\frac{4d+8}{d(d+2)}}\Vert \nabla u_{\sigma,\delta}^{N-m,n}\Vert^{2}+\Vert u_{\sigma,\delta}^{N-m,n} \Vert^{\frac{2d+4}d}\right)ds.
\end{align*}
On the other hand, by \eqref{est}, we obtain
\begin{align*}
&\hat{C}\int_{I_{n+1}}\Vert\nabla u_{\sigma,\delta}^{N-m,n} \Vert^2\Vert u_{\sigma,\delta}^{N-m,n} \Vert^{\frac{4d+8}{d(d+2)}}ds\leq \hat{C}\max_{t\in I_{n+1}}\Vert u_{\sigma,\delta}^{N-m,n} (t)\Vert^{^{\frac{4d+8}{d(d+2)}}}\int_{I_{n+1}}\Vert\nabla u_{\sigma,\delta}^{N-m,n} \Vert^2ds\\&
\leq \frac{{\color{black}}\hat{C}}{2\gamma m\vert B\vert}X_n^{{\frac{2d+4}{d(d+2)}}}2\gamma m\vert B\vert\int_{I_{n+1}}\Vert\nabla u_{\sigma,\delta}^{N-m,n} \Vert^2ds\leq \frac{\hat{C}}{2\gamma m\vert B\vert}X_n^{{\frac{d+2}{d}}}\leq
\frac{{\color{black}}\hat{C}}{2\gamma m\vert B\vert }\dfrac{2^{ \frac{d+2}{d}n +\frac{4(d+2)}{d}}\delta^{\frac {4(d+2)} {d}}}{\widetilde{\tau}^{\frac  {d+2}d}}y_n^{\frac{d+2}{d}}.
\end{align*}
Similarly, using \eqref{est} once more, we have
$$
\hat{C}\int_{I_{n+1}}\Vert u_{\sigma,\delta}^{N-m,n} \Vert^{\frac{2d+4}d}ds\leq 2\hat{C}\widetilde{\tau}X_n^{\frac {d+2} d}=\hat{C}\dfrac{2^{\frac {d+2} dn +\frac{5d+8}{d}}\delta^{\frac {4d+8} d}}{\widetilde{\tau}^{\frac 2 d}}y_n^{\frac{d+2}{d}}.
$$
Therefore, we infer from \eqref{est2} that
\begin{align}
&\nonumber\left(\frac{\delta^2}{2^{n+1}}\right)^3y_{n+1}\leq \left(\int_{I_{n+1}}\int_{\Omega}\vert u_{\sigma,\delta}^{N-m,n} \vert^{\frac{2d+ 4} {d}}dxds\right)^{\frac{3d}{2d+4}}\left(\int_{I_{n+1}}\int_{A_n(s)}dxds\right)^{\frac{4-d}{2d+4}}\\&\leq \delta^{6}\frac{2^{{\frac{3}{2}n+{6}}}\hat{C}^{\frac{3d}{2d+4}}}{\widetilde{\tau}^{\frac 3{2}}}\left(\dfrac{1}{2\gamma m\vert B\vert}+{2\widetilde{\tau}}\right)^{\frac {3d}{2d+4}}y_n^{\frac{5+d}{2+d}}.
\label{est3}
\end{align}
In conclusion, we end up with
\begin{align}
y_{n+1}\leq \frac{2^{{\frac{9}{2}n+{9}}}\hat{C}^{\frac{3d}{2d+4}}}{\widetilde{\tau}^{\frac 3{2}}}\left(\dfrac{1}{2\gamma m\vert B\vert}+{2\widetilde{\tau}}\right)^{\frac {3d}{2d+4}}y_n^{\frac{5+d}{2+d}},\qquad \forall n\geq 0.
\label{last0}
\end{align}
Thus we can apply Lemma \ref{conv}. In particular, we have $b=2^\frac{9}{2}>1$, $C=\frac{2^{{6}}\hat{C}^{\frac{3d}{2d+4}}}{\widetilde{\tau}^{\frac 3{2}}}\left(\dfrac{1}{2\gamma m\vert B\vert}+{2\widetilde{\tau}}\right)^{\frac {3d}{2d+4}}>0$, $\varepsilon=\frac{3}{d+2}$, to get that ${y}_n\to 0$, as long as
$$
{y}_0\leq C^{-\frac{d+2}{3}}b^{-\frac{(d+2)^2}{9}},
$$
i.e.,
\begin{align}
y_0\leq \dfrac{2^{-\left[2(d+2)-\frac{(d+2)^2}{2}\right]}\widetilde{\tau}^\frac {d+2} 2}{\hat{C}^{\frac{d}{2}}\left(\dfrac{1}{2\gamma m\vert B\vert}+{2\widetilde{\tau}}\right)^{\frac d{2}}}.
\label{last}
\end{align}
On the other hand, owing to \eqref{reggtris}, we know that $\Vert \psi^\prime({u}_j)\Vert_{L^\infty(\frac{\tau}{2},\infty;L^1(\Omega))}\leq C(\tau)$ for any $j=1,\ldots,N$ and $\psi^\prime$ is monotone in a neighborhood of $0^+$. Then, we get, for $\delta$ sufficiently small,
\begin{align*}
&y_0=\int_{I_0}\int_{A_0(s)}1dxds\leq\int_{I_0}\int_{\{x\in\Omega:\ \sum_{i\in \pp^{N-m}}{u}_i(x,t) \leq 2\delta^2\}}1dxds\\&\leq \int_{I_0}\int_{A_0(s)}\dfrac{1}{N-m}\sum_{i\in \pp^{N-m}}\frac{\vert \psi^\prime({u}_i)\vert}{-\psi^\prime(2\delta^2)} dxds\leq -\frac{3C(\tau)\widetilde{\tau}}{{\psi^\prime(2\delta^2)}(N-m)}.
\end{align*}
Hence, if we ensure that
$$
-\frac{3C(\tau)\widetilde{\tau}}{{\psi^\prime(2\delta^2)}(N-m)}\leq \dfrac{2^{-\left[2(d+2)-\frac{(d+2)^2}{2}\right]}\widetilde{\tau}^\frac {d+2} 2}{\hat{C}^{\frac{d}{2}}\left(\dfrac{1}{2\gamma m\vert B\vert}+{2\widetilde{\tau}}\right)^{\frac d{2}}},
$$
then \eqref{last} holds. This is equivalent to
\begin{align}
\label{delta}
\frac{3C(\tau)2^{-\left[2(d+2)-\frac{(d+2)^2}{2}\right]}\hat{C}^{\frac{d}{2}}\left(\dfrac{1}{2\gamma m\vert B\vert}+{2\widetilde{\tau}}\right)^{\frac d{2}}}{\widetilde{\tau}^{\frac d 2}(N-m)}\leq -{\psi^\prime(2\delta^2)}.
\end{align}
Having fixed $\widetilde{\tau}$ such that \eqref{tt} holds,
we obtain the result by choosing $\delta$ sufficiently small, since $-\psi^\prime(2\delta^2)\to +\infty$ as $\delta\to0$ by assumption (\textbf{E1}). Notice that $\delta>0$ is fixed and not infinitesimal.

In the end, passing to the limit in $y_n$ as $n\to\infty$, we have obtained that
$$
\left\Vert\left(\sum_{i\in \pp^{N-m}}{u}_i-\delta^2\right)^-\right\Vert_{L^\infty(\Omega\times({T}-\widetilde{\tau},{T}))}=0,
$$
\an{by uniqueness of the limit, since, as $n\to\infty$,
\begin{align*}
y_n\to \left\vert\left\{(x,t)\in \Omega\times[T-\widetilde{\tau},T]: \sum_{i\in \pp^{N-m}}{u}_i\leq \delta^2\right\}\right\vert_d,
\end{align*}
and, on the other hand, $y_n\to0$.}
Notice that, due to the choice of $T$, we have (see \eqref{ttt}) $T-\widetilde{\tau}=\tfrac\tau 2 +\tfrac{(m+1)\tau}{2N}\leq \tau$, therefore we can repeat the same procedure on the interval $({T},T+\widetilde{\tau})$ (the new starting time will be $t_{-1}={T}-2\widetilde{\tau}\geq \tfrac{\tau}{2}$) and so on, reaching eventually the entire interval $\left[\frac\tau 2 +\frac{(m+1)\tau}{2N},+\infty\right)$. Clearly $\delta$ and $\widetilde{\tau}$ do not change,  since the estimates are independent of $T$. Therefore, being $\sigma=1,\ldots,\binom{N}{N-m}$ arbitrary, we have obtained that there exists a $0<\delta_{N-m}\leq \delta_{N-m+1}\leq \tfrac 1 N$ such that, for any possible $\pp^{N-m}$, with $\sigma=1,\ldots, \binom{N}{N-m}$,
\begin{align}
\sum_{i\in \pp^{N-m}}{u}_i(t)\geq \delta>0\quad\text{ a.e. in }\Omega\times\left[\frac\tau 2 +\frac{(m+1)\tau}{2N},+\infty\right),\quad \forall\, \delta\in (0,\delta_{N-m}].
\label{sumsm2}
\end{align}
Recalling Remark \ref{sepimp}, we can deduce that \eqref{sumsm2} actually holds everywhere in $\overline{\Omega}\times\left[\frac\tau 2 +\frac{(m+1)\tau}{2N},+\infty\right)$.
We can thus repeat the procedure increasing $m$, for a finite number of times, until each set $\mathcal{P}_\sigma$ is a singleton (as in the case discussed in Section \ref{exi}). This entails that
there exists $0<\delta\leq \tfrac 1 N$ such that, for any $i=1\ldots,N$,
\begin{align}
{u}_i\geq \delta>0\quad\text{ a.e. in }\Omega\times[\tau,+\infty),
\label{separ1}
\end{align} 	
concluding the proof. Notice that the quantity $\delta$ depends on the initial datum only through the initial energy $\mathcal{E}(0)$ and $\overline{\uu}_0$ since all the estimates involved in this proof are the ones mentioned in Theorem  \ref{thm2}.

\section{Proof of Theorem \protect\ref{global}}

\label{pr3}

By Remark \ref{trick}, we only need to show the existence of a compact
absorbing set. From Theorem \ref{exp}, we deduce that, for any $\mathbf{u}%
_0\in \mathcal{V}_\textbf{M}$, there exist constants
$C_3, C_4>0$ such that
\begin{align*}
\Vert S(t)\mathbf{u}_0\Vert_{\mathcal{V_\textbf{M}}}^2\leq C_3e^{-\omega t}\Vert
\mathbf{u}_0\Vert_{\mathcal{V}_\textbf{M}}^2+C_4\quad \forall t\geq 0.
\end{align*}
Indeed, being $\Psi$ bounded on $[0,1]$ and $0\leq \uu_0\leq1$, it holds
\begin{align}
\frac 1 2\Vert
\mathbf{u}_0\Vert_{\mathcal{V}_\textbf{M}}^2 -C\leq \mathcal{E}(0)\leq \frac 1 2\Vert
\mathbf{u}_0\Vert_{\mathcal{V}_\textbf{M}}^2 +C,
\label{energ0}
\end{align}
for some $C>0$ independent of the initial datum $\uu_0$.
This means that the set
\begin{equation*}
\widetilde{\mathcal{B}}_0:=\left\{\mathbf{u}\in \mathcal{V}_\textbf{M}:\
\Vert \mathbf{u}\Vert_{\mathcal{V}_\textbf{M}}\leq \sqrt{\frac {C_3}{2}+C_4}%
:=R_0\right\}
\end{equation*}
is an absorbing set, i.e., for any bounded set $B\subset \mathcal{V}_\textbf{%
	M}$ there exists $t_e>0$ depending on {\color{black}$\sup_{\uu_0\in B}\Vert\uu_0\Vert_{\mathcal{V}_\textbf{M}}$} such that $S(t)B\subset \widetilde{%
	\mathcal{B}}_0$ for any $t\geq t_e$.

On account of \eqref{H2globbound} and \eqref{separ1}, we can find $\delta=\delta(R_0)>0$ and a bounded set
\begin{equation}
\mathcal{B}_0:=\left\{\mathbf{u}\in \widetilde{\mathcal{B}}_0:\ \Vert
\mathbf{u}\Vert_{\mathbf{H}^{2}(\Omega)}\leq C_0,\quad \uu\geq \delta\text{ in }\overline{\Omega},\quad \partial_\textbf{n}\uu=0\quad \text{ a.e. on }\partial\Omega \right\},
\label{compact}
\end{equation}
for some $C_0=C_0(R_0)>0$, and
a time $t_{R_0}$, depending only on $R_0$, such that $S(t)\widetilde{\mathcal{B}}_0\subset \mathcal{B}_0$ for any $t\geq t_{R_0}$. Note that we can state \textit{for any }$t\geq t_{R_0}$ instead of \textit{for almost any} $t$ (see Remark \ref{sepimp}). This clearly implies that $\mathcal{B}_0$ is a compact absorbing set and ends the proof.
\section{Proof of Theorem \ref{expatt}}
We need some preliminary lemmas. First, recalling \eqref{compact}, we know that there exists $\widetilde{t}=\widetilde{t}(R_0,\mathbf{M})>0$ (with $\mathbf{M}$ fixed) such that $S(t)\mathcal{B}_0\subset \mathcal{B}_0$, for any $t\geq \widetilde{t}$. We then introduce the set
\begin{align*}
\mathbb{B}:=\overline {\bigcup_{t\geq \widetilde{t}}S(t)\mathcal{B}_0}^{\mathcal{V}_\textbf{M}},
\end{align*}
which is compact, positively invariant and absorbing. Let us prove the following
\begin{lemma}
	For any $T\geq 0$ there exists $C=C(T)>0$ such that, given $\uu_{0,1},\uu_{0,2}\in \mathbb{B}$,  it holds
	\begin{align}
	\Vert S(t)\uu_{0,1}-S(t)\uu_{0,2} \Vert_{\mathcal{V}_\mathbf{M}}^2+\int_0^t\Vert \partial_sS(s)\uu_{0,1}- \partial_sS(s)\uu_{0,2}\Vert^2ds\leq C(T)\Vert \uu_{0,1}-\uu_{0,2} \Vert_{\mathcal{V}_\mathbf{M}}^2,\quad \forall t\in[0,T],
	\label{expon}
	\end{align}
	and
		\begin{align}
	\Vert S(t)\uu_{0,1}-S(t)\uu_{0,2} \Vert_{\mathbf{H}^2(\Omega)}^2\leq \frac{C(T)(1+t^2)}{t^2}\Vert \uu_{0,1}-\uu_{0,2} \Vert_{\mathcal{V}_\mathbf{M}}^2,\quad \forall t\in(0,T].
	\label{expo2}
	\end{align}
\end{lemma}
\begin{proof}
	The following computations are formal, but they can be performed within a suitable approximating scheme as the one used in the proof of Theorem \ref{thm}. In particular, leaning on the strict separation property, which holds uniformly (depending only $R_0$ and $\mathbf{M}$, this last one being fixed, see Remark \ref{dependence2}) if the initial data belong to $\mathbb{B}$ (see Theorem \ref{thm}), then we are able to interpret, by uniqueness, the solutions to problem \eqref{syst} as the solutions to a similar problem where $\psi$ is
	replaced by a suitable regular potential (i.e. obtained by extending $\psi$ outside $[\delta,1-(N-1)\delta]$ in a smooth way).
	
	We start by observing that there exists $\delta>0$ (possibly smaller than the one in the definition of $\mathbb{B}$) such that (see \eqref{sep} and \eqref{delt})
	\begin{align}
	S(t)\uu_0\geq \delta\quad \text{ in }\overline{\Omega },\quad \forall\, t\geq 0,\quad \forall\, \uu_0\in\mathbb{B}.
	\label{separa}
	\end{align}
	Set now $\uu^i=S(t)\uu_{0,i}$, with $\uu_{0,i}\in\mathbb{B}$, $i=1,2$. Then, taking the difference between the equations satisfied by $\uu^1$ and $\uu^2$, multiplying it by $\partial_t\uu$, where $\uu=\uu^1-\uu^2$, and integrating over $\Omega$, after an integration by parts, we get
	\begin{align*}
	\frac 1 2\ddt(\al \nabla \uu,\nabla\uu)-\sum_{i,j=1}^N(\alpha_{ij}(\psi^\prime(u_j^1)-\psi^\prime(u_j^2)),\partial_t u_i)+\sum_{i,j=1}^N(\alpha_{ij}(\textbf{Au})_j,\partial_t u_i)+\Vert\partial_t\uu\Vert^2=0,
	\end{align*}
	where we exploited the following facts: $\overline{\partial_t\uu}\equiv0$, $\mathbf{P}(\partial_t\uu)=\partial_t\uu$, and the property $\al(\mathbf{P}\boldsymbol{\xi})=\al\boldsymbol{\xi}$ for any $\boldsymbol{\xi}\in\R^N$.
	
	Thanks to \eqref{separa}, we have $\Vert\psi^{\prime\prime}(su_j^1+(1-s)u_j^2)\Vert_{L^\infty(\Omega)}\leq C$, for any $j=1,\ldots,N$, so that, by standard inequalities,
	\begin{align}
	\nonumber\sum_{i,j=1}^N(\alpha_{ij}(\psi^\prime(u_j^1)-\psi^\prime(u_j^2)),\partial_t u_i)&\nonumber=\sum_{i,j=1}^N\int_\Omega\int_0^1\psi^{\prime\prime}(su_j^1+(1-s)u_j^2)(u_j^1-u_j^2)\alpha_{ij}\partial_t u_idsdx\\&\leq C\Vert \uu\Vert\Vert\partial_t\uu\Vert\leq C\Vert \uu\Vert^2+\frac{1}{4}\Vert \partial_t\uu\Vert^2\leq C\Vert \nabla \uu\Vert^2+\frac{1}{4}\Vert \partial_t\uu\Vert^2.
	\label{psip}
	\end{align}
	Then, similarly,
	\begin{align*}
	\sum_{i,j=1}^N(\alpha_{ij}(\textbf{Au})_j,\partial_t u_i)\leq C\Vert \uu\Vert^2+\frac{1}{4}\Vert\partial_t\uu\Vert^2,
	\end{align*}
	so that, owing to Poincar\'{e}'s inequality, we obtain
		\begin{align}
	\frac 1 2\ddt(\al \nabla \uu,\nabla\uu)+\frac 1 4\Vert \partial_t\uu\Vert^2\leq C(\al \nabla \uu,\nabla\uu),\quad\text{ for a.a. } t\in[0,T].
	\label{e}
	\end{align}
	where we exploited the fact that ($f_{,k}:=\partial_{x_k}f$)
	$$
	(\al \nabla \uu,\nabla\uu)=\sum_{k=1}^d\sum_{i,j=1}^N(\alpha_{ij}u_{i,k},u_{j,k})=\sum_{k=1}^d(\al \uu_{,k},\uu_{,k})\geq C\Vert \nabla \uu\Vert^2,
	$$
	by \eqref{pos} (recall that $\mathbf{P}\uu_{,k}=\uu_{,k}$).
	Thus \eqref{expon} follows from  \eqref{e} owing to Gronwall's Lemma and Poincar\'{e}'s inequality. Notice that the constant $C$, thanks to \eqref{separa}, does not depend on the specific $\uu_{0,i}\in \mathbb{B}$.
	
	Concerning \eqref{expo2}, we write \eqref{phi1} for the difference (defined as $\uu$) between $\uu^1$ and $\uu^2$ and we differentiate the resulting equation with respect to time.
    Then, we multiply it by $\partial_t\uu$ and integrate over $\Omega$. This gives, after an integration by parts, the identity
	$$
	\frac 1 2 \ddt \Vert \partial_t\uu\Vert^2+\sum_{i,j=1}^N(\alpha_{ij}(\psi^{\prime\prime}(u_j^1)\partial_t u_j^1-\psi^{\prime\prime}(u_j^2)\partial_tu_j^2),\partial_t u_i)-\sum_{i,j=1}^N(\alpha_{ij}(\textbf{A}\partial_t\uu)_j,\partial_tu_i)+(\al\nabla\partial_t\uu,\nabla\partial_t\uu)=0,
	$$
	 where we exploited $\overline{\partial_t\uu}\equiv0$, $\mathbf{P}\partial_t\uu=\partial_t\uu$, and the properties of $\al$.
	 Using now \eqref{separa} once more, standard inequalities, and on account of assumption $\psi\in C^3(0,1]$, we get
	 \begin{align*}
	 &\left\vert\sum_{i,j=1}^N(\alpha_{ij}(\psi^{\prime\prime}(u_j^1)\partial_t u_j^1-\psi^{\prime\prime}(u_j^2)\partial_tu_j^2),\partial_t u_i)\right\vert\\&=
	 \left\vert\sum_{i,j=1}^N(\alpha_{ij}(\psi^{\prime\prime}(u_j^1)-\psi^{\prime\prime}(u_j^2))\partial_tu_j^1),\partial_t u_i)\right\vert\\&+\left\vert\sum_{i,j=1}^N(\alpha_{ij}\psi^{\prime\prime}(u_j^2)(\partial_t u_j^1-\partial_tu_j^2),\partial_t u_i)\right\vert\\&\leq \left\vert\sum_{i,j=1}^N\int_\Omega\int_0^1 \alpha_{ij}\psi^{\prime\prime\prime}(su_j^1+(1-s)u_j^2)(u_j^1-u_j^2)\partial_tu_j^1\partial_tu_idsdx\right\vert+C\Vert \partial_t\uu\Vert^2\\&\leq C\Vert \uu\Vert_{\mathbf{L}^4(\Omega)}\Vert\partial_t \uu^1 \Vert\Vert \partial_t\uu\Vert_{\mathbf{L}^4(\Omega)}+C\Vert \partial_t\uu\Vert^2\\&\leq
	 C\Vert \uu\Vert_{\mathcal{V_\textbf{M}}}\Vert \nabla\partial_t\uu\Vert+\Vert \partial_t\uu\Vert^2\leq C(\Vert \uu\Vert_{\mathcal{V_\textbf{M}}}^2+\Vert \partial_t\uu\Vert^2)+\frac{1}{2}(\al\nabla\partial_t\uu,\nabla\partial_t\uu),
	 \end{align*}
	 where we exploited the embedding $\mathbf{H}^1(\Omega)\hookrightarrow \mathbf{L}^4(\Omega)$, the bound $\Vert \partial_t\uu^1\Vert_{L^\infty(0,T;\mathbf{L}^2(\Omega))}\leq C$ with $C$ depending only on $R_0$ (see point (2) of Theorem \ref{thm}, Theorem \ref{thm2}, and \eqref{compact}), Poincar\'{e}'s inequality and the fact that $(\al\nabla\partial_t\uu,\nabla\partial_t\uu)\geq C\Vert \nabla\partial_t\uu\Vert^2$. This last estimate comes from \eqref{pos}, since we have
	 $$
	 (\al\nabla\partial_t\uu,\nabla\partial_t\uu)=\sum_{k=1}^d\sum_{i,j=1}^N(\alpha_{ij}\partial_tu_{i,k},\partial_t u_{j,k})=\sum_{k=1}^d(\al\partial_t\uu_{,k},\partial_t\uu_{,k})\geq C\Vert\nabla\partial_t\uu\Vert^2,
	 $$
	 with $\mathbf{P}\partial_t\uu_{,k}=\partial_t\uu_{,k}$.
	 In conclusion, we have
	 $$
	 \left\vert\sum_{i,j=1}^N(\alpha_{ij}(\textbf{A}\partial_t\uu)_j,\partial_tu_i)\right\vert\leq C\Vert \partial_t\uu\Vert^2.
	 $$
	 We thus end up with
	 $$
	 	\frac 1 2 \ddt \Vert \partial_t\uu\Vert^2+\frac 1 2(\al\nabla\partial_t\uu,\nabla\partial_t\uu)\leq C(\Vert \uu\Vert_{\mathcal{V_\textbf{M}}}^2+\Vert \partial_t\uu\Vert^2),
	 $$
	 and, multiplying both sides by $s^2\in[0,T^2]$, we obtain
	 	 $$
	 \frac 1 2 \ddt s^2\Vert \partial_t\uu\Vert^2+\frac {s^2} 2(\al\nabla\partial_t\uu,\nabla\partial_t\uu)\leq C(s^2\Vert \uu\Vert_{\mathcal{V_\textbf{M}}}^2+(s^2+s)\Vert \partial_t\uu\Vert^2).
	 $$
	 Integrating over $(0,t)$, recalling \eqref{expon}, and dividing by $t^2$, we deduce
	 \begin{align}
	 \Vert \partial_t \uu(t)\Vert\leq \frac{C(T)}{t^2}\Vert \uu_{0,1}-\uu_{0,2}\Vert_{\mathcal{V_\textbf{M}}},\quad \forall t\in (0,T].
	 \label{pt}
	 \end{align}
	We now multiply the equation for $\uu$ by $-\Delta\uu$ and integrate over $\Omega$. We get
	\begin{align}
	-(\partial_t\uu,\Delta\uu)+(\al\Delta\uu,\Delta\uu)-\sum_{i,j=1}^N(\alpha_{ij}(\psi^\prime(u_i^1)-\psi^\prime(u_i^2)),\Delta u_j)+\sum_{i,j=1}^N(\alpha_{ij}(\textbf{A}\uu)_j,\Delta u_j)=0,
	\label{ll}
	\end{align}
	where we used $\overline{\Delta\uu}=0$ and the properties of $\al$.
	Now, being $\mathbf{P}\Delta\uu=\Delta\uu$, we have (see \eqref{pos})
	$$
	(\al\Delta\uu,\Delta\uu)\geq C\Vert \Delta\uu\Vert^2.
	$$
	Moreover, similarly to \eqref{psip} we have
	$$
	\left\vert\sum_{i,j=1}^N(\alpha_{ij}(\psi^\prime(u_i^1)-\psi^\prime(u_i^2)),\Delta u_j)\right\vert\leq C\Vert\uu\Vert^2+\frac 1 4(\al\Delta\uu,\Delta\uu),
	$$
    and Cauchy-Schwarz and Young's inequalities yield
    $$
    \left\vert\sum_{i,j=1}^N(\alpha_{ij}(\textbf{A}\uu)_j,\Delta u_j)\right\vert\leq  C\Vert\uu\Vert^2+\frac 1 4(\al\Delta\uu,\Delta\uu).
    $$
    Therefore, from \eqref{ll} and Poincar\'{e}'s inequality we deduce
    \begin{align*}
    C\Vert \Delta\uu\Vert^2\leq C\Vert \nabla \uu\Vert^2+C\Vert \partial_t\uu\Vert^2,
    \end{align*}
    and combining it with \eqref{expon} and \eqref{pt}, we infer \eqref{expo2}.
\end{proof}
We can now continue the proof of Theorem \ref{expatt} following \cite{Zelik}. By \eqref{ddt}, given $\uu(t)=S(t)\uu_0$, with $\uu_0\in\mathbb{B}$, we have, for any given $T>0$,
\begin{align}
\label{lipt}
\Vert \uu(t)-\uu(s)\Vert_{\mathcal{V_\textbf{M}}}\leq \int_s^t\Vert \partial_t\uu(\tau)\Vert_{\mathcal{V_\textbf{M}}} d\tau\leq \vert t-s\vert^{\frac 1 2}\left(\int_s^t\Vert \partial_t\uu(\tau)\Vert_{\mathcal{V_\textbf{M}}}^2 d\tau\right)^\frac 1 2\leq C(T)\vert t-s\vert^{\frac 1 2}
\end{align}
for any $s,t\in[0,T]$, i.e., $t \mapsto S(t)\uu_0$ is $\tfrac 12$-H\"{o}lder continuous in $[0,T]$, with $C(T)$ depending only on $R_0$. Let us now fix $t^*>0$. Thanks to the smoothing property \eqref{expo2}, valid at $t=t^*>0$, the discrete dynamical system generated by the iterations of $S(t^*)$ possesses an exponential attractor $\mathcal{M}^*\subset \mathbb{B}$ (see, e.g., \cite[Thm.3.7]{Zelik}). Moreover \eqref{expon} and \eqref{lipt} entail
$$
S:[0,t^*]\times \mathbb{B}\to \mathbb{B}, \quad S(t,\uu_0):=S(t)\uu_0,
$$
is H\"{o}lder continuous, when $\mathbb{B}$ is endowed with the $\mathcal{V}_\mathbf{M}$ topology. Therefore, we can define
$$
\mathcal{M}:=\bigcup_{t\in[0,t^*]}S(t)\mathcal{M}^*\subset \mathbb{B},
$$
and, following \cite{Zelik}, show that $\mathcal{M}$ is an exponential attractor for $S(t)$ on $\mathbb{B}$. Since $\mathbb{B}$ is also a compact absorbing set, the basin of exponential attraction of $\mathcal{M}$ is the whole phase space $\mathcal{V}_{\textbf{M}}$. This means that $\mathcal{M}$ is an exponential attractor on $\mathcal{V_\mathbf{M}}$. The proof is finished.

\bigskip
\paragraph{{\bf Acknowledgments.}} \an{The authors are grateful to the anonymous referees for their careful reading of the manuscript as well as for their valuable comments.} The authors have been partially funded by MIUR-PRIN research grant no.~2020F3NCPX ``Mathematics for Industry 4.0 (Math4I4)''.
The authors are also members of Gruppo Nazionale per l'Analisi Ma\-te\-ma\-ti\-ca, la Probabilit\`{a} e le loro Applicazioni (GNAMPA), Istituto Nazionale di Alta Matematica (INdAM).
\an{The present research is part of the activities of ``Dipartimento di Eccellenza 2023-2027''.}


\bigskip
\begin{thebibliography}{}

\bibitem{Abels2009} H. Abels, \textit{On a diffuse interface model for two-phase flows of viscous, incompressible fluids
	with matched densities}, Arch. Ration. Mech. Anal. \textbf{194} (2009), 463-506.

\an{
\bibitem{ATT2023} S. Aihara, N. Takada, T. Takaki, \textit{Highly conservative Allen–Cahn-type multi-phase-field
model and evaluation of its accuracy}, Theor. Comput. Fluid Dyn. (2023), https://doi.org/10.1007/s00162-023-00655-0.
}

\bibitem{AC79} S.M. Allen, J.W. Cahn, \textit{A microscopic theory for antiphase boundary motion and its application to antiphase domain coarsening},
Acta Metall, \textbf{27} (1979), 1085-1095.

\bibitem{BGSS2013} L. Blank, H. Garcke, L. Sarbu, V. Styles, \textit{Nonlocal Allen–Cahn systems: analysis and a primal–dual active set method} IMA J. Numer. Anal. \textbf{33} (2013), 1126-1155.

\bibitem{BW2005} T. Blesgen, U. Weikard, \textit{Multi-component Allen-Cahn equation for elastically stressed solids}, Electron. J. Differential Equations \textbf{89} (2005), 17 pp.

\bibitem{BTP2015} C.P. Brangwynne, P. Tompa, R.V. Pappu, \textit{Polymer physics of intracellular phase transitions}, Nat. Phys. \textbf{11}, 899-904.

\bibitem{BS97} L. Bronsard, B. Stoth, \textit{Volume-preserving mean curvature flow as a limit of a nonlocal Ginzburg-Landau equation}, SIAM J. Math. Anal. \textbf{28} (1997), 769-807.

\bibitem{CH1} J.W. Cahn, J.E. Hilliard, \textit{Free energy of a nonuniform system I. Interfacial free energy},
J. Chem. Phys. \textbf{2} (1958), 258-267.

\bibitem{CH2} J.W. Cahn, J.E. Hilliard, \textit{On spinodal decomposition}, Acta Metall. \textbf{9} (1961), 795-801.

\bibitem{CHL2010} X. Chen, D. Hilhorst, E. Logak, \textit{Mass conserving Allen-Cahn equation and volume preserving mean curvature flow},
Interfaces Free Bound. \textbf{12} (2010), 527-549.
	
\bibitem{DiBenedetto} E. DiBenedetto, \textit{Degenerate parabolic equations}, Springer, New York, 1993.

\bibitem{D1} E. Dolgin, \textit{What lava lamps and vinaigrette can teach us about cell
biology}, Nat. \textbf{555} (2018), 300-302.

\bibitem{D2} E. Dolgin, \textit{The shape-shifting blobs that rule biology}, Nat. \textbf{611} (2022), 24-27.

\bibitem{E} C.M. Elliott, \textit{The Cahn-Hilliard model for the kinetics of phase
separation} in Mathematical models for phase change problems (\'{O}bidos,
1988), Internat. Ser. Numer. Math. \textbf{88}, J. F. Rodrigues, ed., Birkh\"{a}user,
Basel, 1989, 35-73.
	
\bibitem{EL} C.M. Elliott, S. Luckhaus, \textit{A generalized diffusion
		equation for phase separation of a multi component mixture with interfacial
		free energy}, IMA Preprint Series \# 887, 1991.

\bibitem{FLWZ2023} M. Fei, F. Lin, W. Wang, Z. Zhang, \textit{Matrix-valued Allen–Cahn equation
and the Keller–Rubinstein–Sternberg problem}, Invent. Math. (2023), https://doi.org/10.1007/s00222-023-01183-8.

\bibitem{Fu2017} G. Fusco, \textit{Layered solutions to the vector Allen-Cahn equation in R2. Minimizers and heteroclinic connections},
Commun. Pure Appl. Anal. \textbf{16} (2017), 1807-1841.

\bibitem{GGG} C.G. Gal, A. Giorgini, M. Grasselli, \textit{The separation
		property for 2D Cahn-Hilliard equations: local, nonlocal and fractional
		energy cases}, Discrete Contin. Dyn. Syst. \textbf{43} (2023), 2270-2304.
	
\bibitem{GGG2} C.G. Gal, A. Giorgini, M. Grasselli, \textit{The nonlocal
		Cahn-Hilliard equation with singular potential: well-posedness, regularity
		and separation property}, J. Diff. Eqns. 263 (2017), 5253-5297.

\bibitem{GGP} C.G. Gal, M. Grasselli, A. Poiatti, \textit{Allen-Cahn-Navier-Stokes-Voigt systems with moving contact lines,}  ResearchGate preprint 10.13140/RG.2.2.16038.86086, 06 2021.

\bibitem{GGPS} C.G. Gal, M. Grasselli, A. Poiatti,  J. Shomberg, \textit{Multi-component Cahn-Hilliard systems with singular potentials: Theoretical results}, Appl. Math. Optim. \textbf{88}(73), 46 pp.
	
\bibitem{Garke} H. Garcke, \textit{On a Cahn-Hilliard model for phase
		separation with elastic misfit}, Ann. Inst. H. Poincar\'{e} Anal. Non Lin\'{e}aire \textbf{22} (2020), 165-185.

\bibitem{GNSW2008} H. Garcke, B. Nestler, B. Stinner, F. Wendler, \textit{Allen-Cahn systems with volume constraints},  Math. Models Methods Appl. Sci. \textbf{18} (2008), 1347-1381.
		
\bibitem{GGW} A. Giorgini, M. Grasselli, H. Wu, \textit{ On the mass-conserving Allen-Cahn approximation for incompressible binary fluids}, J. Funct. Anal. \textbf{283} (2022), 109631.

\bibitem{G97} D. Golovaty, \textit{The volume-preserving motion by mean curvature as an asymptotic limit of reaction-diffusion equations}, Quart. Appl. Math. \textbf{55} (1997),  243-298.

\bibitem{GPS2006} M. Grasselli, H. Petzeltov\'{a}, G. Schimperna, \textit{Long time behavior of solutions to the Caginalp system with singular potential}, Z. Anal. Anwend. \textbf{25} (2006), 51-72.

\bibitem{KL2017} J. Kim, H. G. Lee, \textit{A new conservative vector-valued Allen–Cahn equation and its fast numerical method}, Comput. Phys. Commun. \textbf{221} (2017), 102-108.

\bibitem{KK2006} R. Kornhuber, R. Krause, \textit{Robust multigrid methods for vector-valued Allen-Cahn equations with logarithmic free energy}, Comput. Vis. Sci. \textbf{9} (2006), 103-116.

\bibitem{KSZ} A. Kostianko, C. Sun, S. Zelik, \textit{Reaction-diffusion systems with supercritical nonlinearities revisited}, Math. Ann. \textbf{384} (2022), 1-45.

\bibitem{Lady} O.A. Lady\v{z}enskaja, V. A. Solonnikov, N. Ural'ceva, \textit{ Linear and quasilinear equations of parabolic type},
Translations of Mathematical Monographs \textbf{23}, American Mathematical Society, Providence, R.I. 1968.

\bibitem{LOCY2022} P. Liu, Z. Ouyang, C. Chen, X. Yang, \textit{A novel fully-decoupled, linear, and unconditionally energy-stable scheme of the
conserved Allen-Cahn phase-field model of a two-phase incompressible flow system with variable density and viscosity},
Commun. Nonlinear Sci. Numer. Simul. \textbf{107} (2022), Paper No. 106120, 21 pp.

\bibitem{M} A. Miranville, \textit{The Cahn-Hilliard Equation: Recent Advances and Applications},
CBMS-NSF Regional Conf. Ser. in Appl. Math., SIAM, Philadelphia, PA., 2019.

\bibitem{Zelik} A. Miranville, S. Zelik, \textit{Attractors for dissipative
	partial differential equations in bounded and unbounded domains}, in
Evolutionary Equations (C.M. Dafermos, M. Pokorn\'{y}, Eds.), Chap. \textbf{4}, 103-200,
Handb. Differ. Equ., Elsevier/North-Holland, Amsterdam, 2008.

\bibitem{MHVB2018} P.G. Moerman, P.C. Hohenberg, P.C., E. Vanden-Eijndenc, J. Brujica, \textit{Emulsion patterns in the wake of a liquid-liquid phase separation front},
Proc. Natl. Acad. Sci. USA \textbf{115} (2018), 3599-3604.

\an{
\bibitem{M2023} M. Moser, \textit{Convergence of the Scalar- and Vector-Valued Allen-Cahn Equation to Mean Curvature Flow with 90°-Contact Angle in Higher Dimensions, part I: Convergence result}, Asymptot. Anal. \textbf{131} (2023), 297-383.
}

\bibitem{NC92} A. Novick-Cohen, \textit{On the viscous Cahn-Hilliard equation}, in Material Instabilities in Continuum
Mechanics and Related Mathematical Problems, J. Ball, ed., Oxford Scientific, Oxford, UK, 1988, 329-342.

\bibitem{P} A. Poiatti, \textit{The 3D strict separation property for the nonlocal Cahn--Hilliard
	equation with singular potential}, 2022, to appear in Anal. PDE, \href{https://arxiv.org/abs/2303.07745v1}{Preprint arXiv:2303.07745}.


\bibitem{RB2023} U. Rana, K. Xu, A. Narayanan, M.T. Walls, A.Z. Panagiotopoulos, J.L. Avalos, C.P. Brangwynne, \textit{Asymmetric oligomerization state and sequence patterning can tune multiphase condensate miscibility}, bioRxiv preprint doi: https://doi.org/10.1101/2023.03.11.532188


\bibitem{RS92} J. Rubinstein, P. Sternberg, \textit{Nonlocal reaction-diffusion equations and nucleation}
IMA J. Appl. Math. \textbf{48} (1992), 249-264.

\bibitem{ST2009} M. Saez Trumper, \textit{Existence of a solution to a vector-valued Allen-Cahn equation with a three well potential},
Indiana Univ. Math. J. \textbf{58} (2009), 213-267.

\bibitem{Temam} R. Temam, \textit{Infinite-Dimensional Dynamical Systems in
	Mechanics and Physics}, Springer-Verlag, New York, 1997.

\bibitem{WS2021} J. Wang, Z. Shi, \textit{Multi-Reconstruction from Points Cloud by Using a Modified
Vector-Valued Allen–Cahn Equation}, Mathematics \textbf{9} (2021), 1326. https://doi.org/10.3390/math9121326.

\bibitem{Y2021} X. Yang, \textit{Efficient, second-order in time, and energy stable scheme for a new hydrodynamically coupled three components volume-conserved Allen-Cahn phase-field model},
Math. Models Methods Appl. Sci. \textbf{31} (2021), 753-787.

\bibitem{YH2022} X. Yang, X. He, \textit{A fully-discrete decoupled finite element method for the conserved
Allen–Cahn type phase-field model of three-phase fluid flow system}, Comput. Methods Appl. Mech. Engrg. \textbf{389} (2022), 114376.

\bibitem{ZL2022} D. Zwicker, L. Laan, \textit{Evolved interactions stabilize many coexisting phases in multicomponent liquids}, Proc. Natl. Acad. Sci. USA \textbf{119} (2022), e2201250119, 8 pp.



\end{thebibliography}
\end{document}